\documentclass[12pt,a4,reqno]{amsart}

\usepackage{amscd}

\newcommand{\N}{{\mathbb N}}
\newcommand{\Z}{{\mathbb Z}}

\newcommand{\C}{{\mathbb C}}
\newcommand{\R}{{\mathbb R}}
\renewcommand{\P}{{\mathbb P}}

\newcommand{\GG}{{\mathcal G}}
\newcommand{\OO}{{\mathcal O}}
\newcommand{\QQ}{{\mathcal Q}}
\newcommand{\KK}{{\mathcal K}}
\newcommand{\BB}{{\mathcal B}}
\newcommand{\DD}{{\mathcal D}}
\newcommand{\PP}{{\mathcal P}}
\newcommand{\RR}{{\mathcal R}}
\newcommand{\TT}{{\mathcal T}}
\newcommand{\UU}{{\mathcal U}}
\newcommand{\aaa}{{\bf a}}
\newcommand{\mmm}{{\bf m}}
\newcommand{\qqq}{{\bf q}}
\newcommand{\tm}{{\mathcal T}_M}
\newcommand{\modmu}{{\rm mod}_\mu}
\renewcommand{\Im}{{\rm Im}}

\DeclareMathOperator{\specan}{Specan}
\DeclareMathOperator{\Der}{Der}
\DeclareMathOperator{\lie}{Lie}
\DeclareMathOperator{\sym}{Sym}
\DeclareMathOperator{\embdim}{embdim}
\DeclareMathOperator{\Ann}{Ann}
\DeclareMathOperator{\codim}{codim}
\DeclareMathOperator{\rank}{rank}
\DeclareMathOperator{\ad}{ad}
\DeclareMathOperator{\mult}{mult}

\begin{document}

\theoremstyle{plain}
\newtheorem{lemma}{Lemma}[section]
\newtheorem{theorem}[lemma]{Theorem}
\newtheorem{proposition}[lemma]{Proposition}
\newtheorem{corollary}[lemma]{Corollary}

\theoremstyle{definition}
\newtheorem{definition}[lemma]{Definition}
\newtheorem{withouttitle}[lemma]{}
\newtheorem{remark}[lemma]{Remark}
\newtheorem{remarks}[lemma]{Remarks}
\newtheorem{example}[lemma]{Example}
\newtheorem{examples}[lemma]{Examples}

\title{Multiplication on the tangent bundle} 

\author{
Claus Hertling}

\address{Claus Hertling\\
Mathematisches Institut der Universit\"at Bonn\\
Beringstra{\ss}e 1, 53115 Bonn, Germany}

\email{hertling\char64 math.uni-bonn.de}

\date{September 26, 1999}

\begin{abstract}
Manifolds with a commutative and associative multiplication on the
tangent bundle are called F-manifolds if a unit field exists
and the multiplication satisfies a natural integrability condition. 
They are studied here. They are closely related to
discriminants and Lagrange maps. Frobenius manifolds are F-manifolds.
As an application a conjecture of Dubrovin on Frobenius manifolds
and Coxeter groups is proved.
\end{abstract}

\maketitle


\tableofcontents

\setcounter{section}{-1}

\section{Introduction}
\setcounter{equation}{0}

In \cite{HM}\cite{Man}(I.5) the notion of an F-manifold was introduced.
It is a triple $(M,\circ,e)$ where $M$ is a manifold, $\circ$ is an
$\OO_M$-bilinear commutative and associative multiplication on the
tangent sheaf $\tm$, $e$ is a global unit field, and the multiplication
satisfies the integrability
condition
\begin{eqnarray}
\lie_{X\circ Y}(\circ) = X\circ \lie_Y(\circ)+Y\circ\lie_X(\circ)
\end{eqnarray}
for any two local vector fields $X$ and $Y$. 
Frobenius manifolds are F-manifolds (\cite{HM}\cite{Man}(I.5)); 
this motivated the definition of this notion.

This paper is devoted to the local structure of complex F-manifolds.
It turns out to be closely related to singularity theory.
Discriminants and Lagrange maps will play a key role.

\medskip
Distinguished (germs of) F-manifolds with many typical and some special
properties are the base spaces of semiuniversal unfoldings of 
isolated hypersurface singularities and of boundary singularities
(chapters 16 and 17). Here the tangent space at each parameter  is 
canonically isomorphic to the sum of the Jacobi algebras of the 
singularities above this parameter. Many of the general results on 
F-manifolds have been known in another guise in the hypersurface singularity
case and all should be compared with it.

One reason why the integrability condition (0.1) is natural is the following:
Let $(M,p)$ be the germ of an F-manifold $(M,\circ,e)$.
The algebra $T_pM$ decomposes uniquely into a sum of (irreducible) 
local algebras which annihilate one another (Lemma 1.1).
Now the integrability condition (0.1) ensures that this infinitesimal
decomposition extends to a unique decomposition of the germ 
$(M,p)$ into a product of germs of F-manifolds (Theorem 4.2).

If the multiplication at $T_pM$ is semisimple, that is, if $T_pM$
decomposes into 1-dimensional algebras, then this provides
canonical coordinates $u_1,...,u_n$ on $(M,p)$ with 
${\partial \over \partial u_i}\circ 
{\partial \over \partial u_j} =
\delta_{ij}{\partial \over \partial u_i}$.
In fact, at points with semisimple multiplication the integrability 
condition (0.1) is equivalent to the existence of such canonical 
coordinates. In the hypersurface case, the decomposition of the germ
$(M,p)$ for some parameter $p$ is the unique decomposition into a 
product of base spaces of semiuniversal unfoldings of the singularities
above $p$.

Another reason why (0.1) is natural is the relation to the potentiality of
Frobenius manifolds.
There exist F-manifolds such that not all tangent spaces are Frobenius
algebras. They cannot be Frobenius manifolds. But if all tangent 
spaces are Frobenius algebras then the integrability condition (0.1)
is related to a version  of potentiality which requires a metric on 
$M$ which is multiplication invariant, but not necessarily flat.
See chapter 5 for details.

This allows a viewpoint on Frobenius manifolds which was already 
stressed in \cite{HM}: In the construction of some Frobenius manifolds,
e.g. in the hypersurface case, the multiplication comes first and is 
canonical, and the most difficult property of a metric is not 
multiplication invariance or potentiality, but flatness.

\smallskip
The most important geometric object which is attributed to an $n$-dimensional
manifold $M$ with multiplication $\circ$ on the tangent sheaf $\tm$ and 
unit field $e$ (with or without (0.1)) is the analytic spectrum
$L:=\specan (\tm)\subset T^{*}M$ (see chapter 2).
The projection $\pi:L\to M$ is flat and finite of degree $n$.
The fiber $\pi^{-1}(p)\subset L$ above $p\in M$ consists of the
characters of $(T_pM,\circ)$; they correspond 1-1 to the 
irreducible subalgebras of $(T_pM,\circ)$ (see Lemma 1.1).
The multiplication on $\tm$ can be recovered from $L$, because the map
\begin{eqnarray}
\aaa:\tm \to \pi_{*}\OO_L,
\ X\mapsto \alpha(\widetilde{X})|_L
\end{eqnarray}
is an isomorphism of $\OO_M$-algebras; 
here $\widetilde{X}$ is any lift of $X$ to $T^{*}M$
and $\alpha$ is the canonical 1--form on $T^{*}M$. 
The values of the function $\aaa(X)$ on $\pi^{-1}(p)$ are the eigenvalues
of $X\circ:T_pM\to T_pM$.

The analytic spectrum $L$ is a reduced variety iff the multiplication is 
generically semisimple. Then the manifold with multiplication
$(M,\circ, e)$ is called
{\it massive}. Now, a third reason why the integrability condition 
(0.1) is natural is this:
$L\subset T^{*}M$ is a Lagrange variety iff $(M,\circ,e)$ is a massive
F-manifold (Theorem 6.2).

The main body of this paper is devoted to study germs of massive F-manifolds
at points where the multiplication is not semisimple.

We will make use of the theory of singular Lagrange varieties and their
Lagrange maps, which has been worked out by Givental in \cite{Gi2}.
In fact, the notion of an irreducible germ (with respect to the 
decomposition above) of a massive F-manifold is equivalent to
Givental's notion of a miniversal germ of a flat Lagrange map 
(Theorem 9.4). Via this equivalence Givental's paper contains many results
on massive F-manifolds and will be extremely useful.

Locally the canonical 1--form $\alpha$ on $T^{*}M$ can be integrated on the
analytic spectrum $L$ of a massive F-manifold $(M,\circ,e)$ to a 
{\it generating function} $F:L\to \C$ which is continuous on $L$ and 
holomorphic on $L_{reg}$. It depends on a property of $L$,
which is weaker than normality or maximality of the complex 
structure of $L$, whether $F$ is holomorphic on $L$ (see chapter 7).

If $F$ is holomorphic on $L$ then it corresponds via (0.2) to
an {\it Euler field} $E=\aaa^{-1}(F)$ of weight 1, that is, a vector field
on $M$ with $\lie_E(\circ)=\circ$ (Theorem 6.3).

In any case, a generating function $F:L\to \C$ gives rise to a 
Lyashko-Looijenga map $\Lambda:M\to \C^n$ (see chapters 8 and 10) 
and a discriminant $\DD=\pi(F^{-1}(0))\subset M$.

If $F$ is holomorphic and an Euler field $E=\aaa^{-1}(F)$ exists then the
discriminant $\DD$ is the hypersurface of points where the multiplication
with $E$ is not invertible. Then it is a free divisor with 
logarithmic fields $\Der_M(\log \DD) = E\circ \tm$ (Theorem 13.1).
This generalizes a result of K. Saito for the hypersurface case.

From the unit field $e$ and a discriminant $\DD\subset M$ one can 
reconstruct everything. One can read off the multiplication on $TM$ 
in a very nice elementary way (Corollary 11.6):
The $e$-orbit of a generic point $p\in M$ intersects $\DD$ in $n$ points.
One shifts the $n$ tangent hyperplanes with the flow of $e$ to $T_pM$.
Then there exist unique vectors $e_1(p),...,e_n(p)\in T_pM$ such
that $\sum_{i=1}^ne_i(p)=e(p)$ and $\sum_{i=1}^n\C\cdot e_i(p)=T_pM$
and such that the subspaces $\sum_{i\neq j}\C\cdot e_i(p)$, 
$j=1,...,n$, are the shifted hyperplanes. 
The multiplication on $T_pM$ is given by 
$e_i(p)\circ e_j(p)=\delta_{ij}e_i(p)$.

In the case of hypersurface singularities and boundary singularities, 
the classical discriminant in the base space of a semiuniversal unfolding
is such a discriminant.
The critical set $C$ in the total space of the unfolding is canonically
isomorphic to the analytic spectrum $L$; this isomorphism identifies
the map $\aaa$ in (0.2) with a Kodaira-Spencer map
$\aaa_C:\tm \to (\pi_C)_{*}\OO_C$ and a generating function 
$F:L\to \C$ with the restriction of the unfolding function to the 
critical set $C$. This Kodaira-Spencer map $\aaa_C$ is the source of the 
multiplication on $\tm$ in the hypersurface singularity case.
The multiplication on $\tm$ had been defined in this way first by 
K. Saito.

Critical set and analytic spectrum are smooth in the hypersurface singularity
case. By work of Arnold on Lagrange maps and singularities 
a beautiful correspondence holds (Theorem 16.6):
each irreducible germ of a massive F-manifold with smooth analytic spectrum
comes from an isolated hypersurface singularity, and this singularity 
is unique up to stable right equivalence.

By work of Nguyen huu Duc and Nguyen tien Dai the same correspondence 
holds for boundary singularities and irreducible germs of massive
F-manifolds whose analytic spectrum is isomorphic
to $(\C^{n-1},0)\times (\{(x,y)\in \C^2\ |\ xy=0\},0)$ with ordered
components (Theorem 17.6).

The complex orbit space $M:=\C^n/W\cong \C^n$ of an irreducible 
Coxeter group $W$ carries an (up to some rescaling) canonical 
structure of a massive F-manifold:
A generating system $P_1,...,P_n$ of $W$-invariant homogeneous polynomials
induces coordinates $t_1,...,t_n$ on $M$. Precisely one polynomial,
e.g. $P_1$, has highest degree. The field 
${\partial \over \partial t_1}$ is up to a scalar independent of any choices.
This field
$e:={\partial \over \partial t_1}$ as unit field and the classical discriminant
$\DD\subset M$, the image of the reflection hyperplanes, determine
in the elementary way described above the structure of a massive
F-manifold. This follows from \cite{Du1}\cite{Du2}(Lecture 4)
as well as from \cite{Gi2}(Theorem 14).

Dubrovin established the structure of a Frobenius manifold on the complex
orbit space $M=\C^n/W$, with this multiplication, with K. Saito's 
flat metric on $M$, and with a canonical Euler field with positive weights
(see Theorem 19.1). At the same place he conjectured that these Frobenius
manifolds and products of them are (up to some well understood rescalings)
the only massive Frobenius manifolds with an Euler field with positive 
weights. We will prove this conjecture (Theorem 19.3).

Crucial for the proof is Givental's result \cite{Gi2}(Theorem 14).
It characterizes the germs $(M,0)$ of F-manifolds of irreducible
Coxeter groups by geometric properties (see Theorem 18.4).
We obtain from it the following intermediate result (Theorem 18.3):
An irreducible germ $(M,p)$ of a {\it simple} F-manifold such that
$T_pM$ is a Frobenius algebra is isomorphic to the germ at 0 of the 
F-manifold of an irreducible Coxeter group.

A massive F-manifold $(M,\circ,e)$ is called {\it simple} if the
germs $(M,p)$, $p\in M$, of F-manifolds 
are contained in finitely many isomorphism classes.
Theorem 18.3 complements in a nice way the relation of irreducible 
Coxeter groups to the simple hypersurface singularities
$A_n, D_n,E_n$ and the simple boundary singularities $B_n, C_n ,F_4$.

In dimension 1 and 2, up to isomorphism all the irreducible germs of 
massive F-manifolds come from the irreducible Coxeter groups 
$A_1$ and $I_2(m)$ ($m\geq 3$) 
with $I_2(3)=A_2,\ I_2(4)=B_2,\ I_2(5)=:H_2,\ I_2(6)=G_2$.
But already in dimension 3 the classification is vast (see chapter 20).

\medskip
This paper is organized as follows. In chapters 1 -- 5 classical facts on 
finite dimensional algebras and the analytic spectrum are presented, 
as well as the definition of F-manifolds, the local decomposition and the 
relation to potentiality.
In chapters 6 -- 10 the Lagrange property of the analytic spectrum of a
massive F-manifold and related objects as generating functions,
Lyashko-Looijenga maps, and Euler fields are discussed.
The chapters 11 -- 15 contain a mixture of results on discriminants,
2-dimensional F-manifolds, logarithmic vector fields,
modality for germs of F-manifolds and different embeddings of analytic
spectra.
The long final chapters 16 -- 20 treat hypersurface singularities, boundary
singularities, finite Coxeter groups and 3-dimensional massive
F-manifolds.

\medskip
The starting point for this work was the common paper \cite{HM}
with Yu. Manin, whose interest was very encouraging and helpful.
I also want to thank E. Brieskorn, my teacher in singularity theory.
This paper builds on work of him, Arnold, Givental, 
Looijenga, Lyashko, K. Saito, O.P. Shcherbak, Teissier,
and many others.
The major part of this paper was written at the mathematics department
of the University Paul Sabatier in Toulouse. I thank the department
and especially J.-F. Mattei for their hospitality.

\section{Finite-dimensional algebras}
\setcounter{equation}{0}

In this chapter $(Q,\circ,e)$ is a $\C$-algebra of finite dimension
($\geq 1$) with commutative and associative multiplication and with
unit $e$. The next lemma gives precise information on the decomposition
of $Q$ into irreducible algebras. The statements are well known and
elementary. They can be deduced directly in the given order or from
more general results in commutative algebra. Algebra homomorphisms
are always supposed to map the unit to the unit.

\begin{lemma} Let $(Q,\circ, e)$ be as above. As the endomorphisms
$x\circ:Q\to Q$, $x\in Q$, commute, there is a unique simultaneous
decomposition $Q=\bigoplus_{k=1}^l Q_k$ into generalized eigenspaces
$Q_k$ (with $\dim_\C Q_k \geq 1$). Define $e_k\in Q_k$ by
$e=\sum_{k=1}^l e_k$. Then
\begin{list}{}{}
\item[i)] 
$Q_j\circ Q_k=0$ for $j\neq k$; $e_k\neq 0$; 
$e_j\circ e_k=\delta_{jk}e_k$; $e_k$ is the unit of the algebra
$Q_k=e_k\circ Q$.
\item[ii)] 
The function $\lambda_k:Q\to \C$ which associates to $x\in Q$ the 
eigenvalue of $x\circ$ on $Q_k$ is an algebra homomorphism; 
$\lambda_j\neq \lambda_k$ for $j\neq k$.
\item[iii)]
$(Q_k,\circ,e_k)$ is an irreducible and a local algebra with 
maximal ideal $\mmm_k=Q_k\cap \ker (\lambda_k)$.
\item[iv)]
The subsets $\ker (\lambda_k) = 
\mmm_k\oplus \bigoplus_{j\neq k} Q_j$, $k=1,..,l$,
are the maximal ideals of the algebra $Q$; the complement
$Q-\bigcup_k \ker (\lambda_k)$ is the group of invertible elements
of $Q$.
\item[v)] $\{\lambda_1,...,\lambda_l\} = {\rm Hom}_{\C -alg}(Q,\C)$.
\item[vi)] the localization $Q_{\ker (\lambda_k)}$ is isomorphic to $Q_k$.
\end{list}
\end{lemma}

We call this decomposition the {\it eigenspace decomposition} of $(Q,\circ,e)$.
The set $L:=\{\lambda_1,...,\lambda_l\}\subset Q^{*}$ carries a natural
complex structure $\OO_L$ such that $\OO_L(L)\cong Q$ and 
$\OO_{L,\lambda_k}\cong Q_k$. More details on this will be given in 
chapter 2.

The algebra (or its multiplication) is called {\it semisimple} if $Q$ 
decomposes into 1-dimensional subspaces,
$Q\cong \bigoplus_{k=1}^{\dim Q} Q_k = \bigoplus_{k=1}^{\dim Q} \C\cdot e_k$.

An irreducible algebra $Q=\C\cdot e\oplus \mmm$ with maximal ideal $\mmm$ is a
Gorenstein ring if the socle ${\rm Ann}_Q(\mmm)$ has dimension 1.

An algebra $Q=\bigoplus_{k=1}^l Q_k$ is a Frobenius algebra if each 
irreducible subalgebra is a Gorenstein ring.
 
The next (also wellknown) lemma gives equivalent conditions and
additional information. We remark that this classical definition of a 
Frobenius algebra is slightly weaker than Dubrovin's: he calls an algebra
$(Q,\circ,e)$ together with a {\it fixed} bilinear form $g$ as in 
Lemma 1.2 a) iii) a Frobenius algebra.

\begin{lemma}
a) The following conditions are equivalent.
\begin{list}{}{}
\item[i)] $(Q,\circ,e)$ is a Frobenius algebra.
\item[ii)] ${\rm Hom}(Q,\C) \cong Q$ as a $Q$-module.
\item[iii)] There exists a bilinear form $g:Q\times Q \to \C$ 
which is symmetric, nondegenerate and multiplication invariant,
i.e. $g(a\circ b,c)=g(a,b\circ c)$.
\end{list}
b) Let $Q=\bigoplus_{k=1}^l Q_k$ be a Frobenius algebra and 
$Q_k=\C\cdot e_k\oplus \mmm_k$. The generators of ${\rm Hom}(Q,\C)$ as a
$Q$-module are the linear forms $f:Q\to \C$ with 
$f({\rm Ann}_{Q_k}(\mmm_k))=\C$ for all $k$.

One obtains a 1-1 correspondence between these linear forms and the 
bilinear forms $g$ as in a) iii) by putting 
$g(x,y):= f(x\circ y)$.
\end{lemma}

{\bf Proof:} 
a) Any of the conditions i), ii), iii) in a) is satisfied for $Q$ iff it is
satisfied for each irreducible subalgebra $Q_k$. One checks this with
$Q_j\circ Q_k =0$ for $j\neq k$. So we may suppose that $Q$ is irreducible.

i)$\iff$ii) A linear form $f\in {\rm Hom}(Q,\C)$ generates
${\rm Hom}(Q,\C)$ as a $Q$-module iff the linear form
$(x\mapsto f(y\circ x))$ is nontrivial for any $y\in Q-\{0\}$, that is, iff
$f(y\circ Q)=\C$ for any $y\in Q-\{0\}$.

The socle ${\rm Ann}_Q(\mmm)$ is the set of the common eigenvectors of all
endomorphisms $x\circ:Q\to Q$, $x\in Q$. If $\dim {\rm Ann}_Q(\mmm)\geq 2$
then for any linear form $f$ an element 
$y\in (\ker f\cap {\rm Ann}_Q(\mmm))-\{0\}$ satisfies $y\circ Q = \C\cdot y$ 
and $f(y\circ Q)=0$; so $f$ does not generate 
${\rm Hom}(Q,\C)$. If $\dim {\rm Ann}_Q(\mmm)=1$ then it is contained in any 
nontrivial ideal, because any such ideal contains a common eigenvector of 
all endomorphisms. 
The set $y\circ Q$ for $y\in Q-\{0\}$ is an ideal.
So, then a linear form $f$ with 
$f({\rm Ann}_Q(\mmm))=\C$ generates ${\rm Hom}(Q,\C)$ as a $Q$-module.

i) $\Rightarrow$ iii) Choose any linear form $f$ with $f({\rm Ann}_Q(\mmm))=\C$
and define $g$ by $g(x,y):=f(x\circ y)$. It rests to show that $g$ is 
nondegenerate. But for any $x\in Q-\{0\}$ there exists a $y\in Q$ with
$\C\cdot x\circ y={\rm Ann}_Q(\mmm)$, because ${\rm Ann}_Q(\mmm)$ is contained
in the ideal $x\circ Q$.

iii) $\Rightarrow$ i) $g(\mmm,{\rm Ann}_Q(\mmm))=g(e,\mmm\circ{\rm Ann}_Q(\mmm))
=g(e,0)=0$ implies $\dim {\rm Ann}_Q(\mmm)=1$.

b) Starting with a bilinear form $g$, the corresponding linear form $f$
is given by $f(x)=g(x,e)$. The rest is clear form the preceding 
discussion. \hfill $\qed$

\bigskip
The semisimple algebra $Q\cong \bigoplus_{k=1}^{\dim Q}\C\cdot e_k$ 
is a Frobenius algebra. 

A classical result is that the complete intersections 
$\OO_{\C^m,0}/(f_1,...,f_m)$ are Gorenstein. But these are not all, e.g.
$\C\{ x,y,z\}/(x^2,y^2,xz,yz,xy-z^2)$ is Gorenstein, but not a 
complete intersection.

Finally, in the next chapter vector bundles with multiplication will be 
considered. Condition ii) of Lemma 1.2 a) shows that there the points whose 
fibers are Frobenius algebras form an open set in the base.

\section{Vector bundles with multiplication}
\setcounter{equation}{0}

Now we consider a holomorphic vector bundle $Q\to M$ on a complex 
manifold $M$ with multiplication on the fibers: 
The sheaf $\QQ$ of holomorphic sections of the bundle $Q\to M$ is 
equipped with an $\OO_M$-bilinear commutative and associative
multiplication $\circ $ and with a global unit section $e$.

The set $\bigcup_{p\in M} {\rm Hom}_{\C -alg}(Q(p),\C)$ of 
algebra homomorphisms from the single fibers $Q(p)$ to $\C$ 
(which map the unit to $1\in \C$) is a subset of the dual bundle
$Q^{*}$ and has a natural complex structure. It is the analytic 
spectrum $\specan (\QQ)$. We sketch the definition \cite{Ho}(ch. 3):

The $\OO_m$-sheaf $\sym_{\OO_m}\QQ$ can be identified with the 
$\OO_m$-sheaf of holomorphic functions on $Q^{*}$ which are polynomial
in the fibers. The canonical $\OO_m$-algebra homomorphism
$\sym_{\OO_m}\QQ \to \QQ$ which maps the multiplication in 
$\sym_{\OO_m}\QQ$  to the multiplication $\circ$ in $\QQ$ is 
surjective. The kernel generates an ideal ${\mathcal I}$ in 
$\OO_{Q^{*}}$. 

One can describe the ideal locally explicitly: suppose $U\subset M$
is open and $\delta_1,...,\delta_n\in \QQ (U)$ is a base of sections of 
the restriction $Q|_U \to U$ with $\delta_1=e$ and 
$\delta_i\circ \delta_j = \sum_k a_{ij}^k \delta_k$; denote by 
$y_1,...,y_n$ the fiberwise linear functions on $Q^{*}|U$ which are 
induced by $\delta_1,...,\delta_n$; then the ideal ${\mathcal I}$ is 
generated in $Q^{*}|U$ by
\begin{eqnarray}
y_1 -1 \mbox{ \ and \ } y_iy_j - \sum_k a_{ij}^k y_k\ .
\end{eqnarray}
The support of $\OO_{Q^{*}} / {\mathcal I}$ with the 
restriction of $\OO_{Q^{*}} / {\mathcal I}$ as structure sheaf is 
$\specan (\QQ) \subset Q^{*}$. We denote the natural projections
by $\pi_{Q^{*}} : Q^{*}\to M$ and $\pi : \specan (\QQ)\to M$.
A part of the following theorem is already clear from 
the discussion. A complete proof and 
thorough discussion can be found in \cite{Ho}(ch. 3).

\begin{theorem}
The set $\bigcup_{p\in M} {\rm Hom}_{\C -alg}(Q(p),\C)$
is the support of the analytic spectrum $\specan (\QQ) =:L$.
The composition of maps
\begin{eqnarray}
\aaa: \QQ \hookrightarrow (\pi_{Q^{*}})_{*} \OO_{Q^{*}} 
               \to \pi_{*} \OO_L
\end{eqnarray}
is an isomorphism of $\OO_M$-algebras and of free $\OO_M$-modules of 
rank $n$, here $n$ is the fiber dimension of $Q\to M$. The projection
$\pi:L\to M$ is finite and flat of degree $n$.

Consider a point $p\in M$ with eigenspace decomposition
$Q(p)=\bigoplus_{k=1}^{l(p)} Q_k(p)$ and 
$L\cap \pi^{-1}(p) = \{ \lambda_1,...,\lambda_{l(p)}\}$. 
The restriction of the isomorphism
\begin{eqnarray}
\QQ_p \stackrel{\aaa}{\longrightarrow} (\pi_{*}\OO_L)_p 
       \cong \bigoplus_{k=1}^{l(p)} \OO_{L,\lambda_k}
\end{eqnarray}
to the fiber over $p$ yields isomorphisms
\begin{eqnarray}
Q_k(p) \cong \OO_{L,\lambda_k} \otimes_{\OO_{M,p}} \C\ .
\end{eqnarray}
\end{theorem}

\begin{corollary}
In a sufficiently small neighborhood $U$ of a point $p\in M$, the 
eigenspace decomposition $Q(p)=\bigoplus_{k=1}^{l(p)} Q_k(p)$ 
of the fiber $Q(p)$ extends uniquely to a decomposition of the 
bundle $Q|_U \to U$ into multiplication invariant holomorphic
subbundles.
\end{corollary}

{\bf Proof of Corollary 2.2:} The $\OO_{M,p}$-free submodules
$\OO_{L,\lambda_k}$ in the decomposition in (2.3) of 
$(\pi_{*}\OO_L)_p$ are obviously multiplication invariant. 
Via the isomorphism $\aaa$ one obtains locally a decomposition
of the sheaf $\QQ$ of sections of $Q\to M$ into multiplication 
invariant free $\OO_M$-submodules.    \hfill $\qed$

\bigskip
Of course, the induced decomposition of $Q(q)$ for a point $q$ in the
neighborhood of $p$ may be coarser than the eigenspace decomposition
of $Q(q)$.

The base is naturally stratified with respect to the numbers and 
dimensions of the components of the eigenspace decompositions of the 
fibers of $Q\to M$. To make this precise we introduce a partial
ordering $\succ$  on the set ${\mathcal P}$ of partitions of $n$:
$${\mathcal P}:= \{\beta =(\beta_1,...,\beta_{l(\beta)})\ |\ \beta_i\in \N,\ 
\beta_i\geq \beta_{i+1},\ \sum_{i=1}^{l(\beta)} \beta_i=n\};$$
for 
$\beta,\gamma\in {\mathcal P}$ 
define
$$\beta \succ \gamma :\iff \exists \ 
\sigma:\{1,...,l(\gamma)\}\to \{1,...,l(\beta)\}
\mbox{ s.t. } \  
\beta_j= \sum_{i\in \sigma^{-1}(j)} \gamma_i \ .$$
The 
partition $P(p)$ of a fiber $Q(p)$ is the partition of $n=\dim Q(p)$ by 
the dimensions of the subspaces of the eigenspace decomposition.

\begin{proposition}
Fix a partition $\beta\in {\mathcal P}$.
The subset $\{p\in M\ |\ P(p)\succ \beta\}$ 
is empty or an analytic subset of $M$.
\end{proposition}

{\bf Proof:}
The partition $P(f)$ of a polynomial of degree $n$ is the partition
of $n$ by the multiplicity of the zeros of $f$.

{\bf Fact:} {\it The space 
$\{ a\in \C^n\ |\ P(z^n+\sum_{i=1}^n a_iz^{n-i}) \succ \beta\}$ is an
algebraic subvariety of $\C^n$ with normalization isomorphic to
$\C^{l(\beta)}$.}

For the proof one has just to regard the finite map 
$\C^n\to \C^n$, $u\mapsto((-1)^i\sigma_i(u))_{i=1,..,n}$
($\sigma_i (u)$ is the $i$-th symmetric polynomial).

A section $X\in \QQ(U)$, $U\subset M$ open, induces via the 
coefficients of the characteristic polynomial $p_{ch,X\circ}$ 
of multiplication by $X$ a holomorphic map $U\to \C^n$. Hence the
set $\{q\in U\ |\ P(p_{ch,X\circ })\succ \beta\}$ is analytic. 
The intersection of such analytic sets for a basis of sections in 
$U$ is $\{q\in U\ |\ P(q)\succ \beta\}$.  \hfill $\qed$

\bigskip
We suppose that $M$ is connected. Then there is a unique partition
$\beta_0$ such that $\{p\in M\ |\ P(p)=\beta_0\}$ is open in $M$.
The complement $\KK := \{p\in M\ |\ P(p)\neq \beta_0\}$ will be called
the {\it caustic} in $M$; this name is motivated by the 
Lagrange maps (chapters 6, 8, 9) and the hypersurface 
singularities (chapter 16). The multiplication is generically 
semisimple iff $\beta_0 = (1,...,1)$.

\begin{proposition}
If the multiplication is generically semisimple then the caustic 
$\KK$ is a hypersurface or empty.
\end{proposition}

{\bf Proof:}
Suppose $\dim (\KK,p)\leq \dim M -2$ for some point $p\in \KK$. Then in
a neighborhood $U$ the complement $U-\KK$ is simply connected.
There is no monodromy for the locally defined idempotent sections
$e_1,...,e_n$ on $U-\KK$. They are sections on $U-\KK$ and extend
to a basis of sections on $U$, with multiplication
$e_i\circ e_j = \delta_{ij}e_i$. Hence $\KK\cap U= \emptyset$.
\hfill $\qed$

\section{Definition of F-manifolds}
\setcounter{equation}{0}

An F-manifold is a manifold $M$ with a multiplication on the tangent
bundle $TM$ which harmonizes with the Lie bracket in the most natural 
way. A neat formulation of this property requires the Lie derivative
of tensors.

\bigskip
\begin{remark}
Here the sheaf of $(k,l)$-tensors ($k,l\in \N_0$) on a manifold $M$ is 
the sheaf of $\OO_M$-module homomorphisms
${\rm Hom}_{\OO_M} (\bigotimes_{i=1}^k\tm , \bigotimes_{i=1}^l\tm )$.
A $(0,l)$-tensor $T:\OO_M \to \bigotimes_{i=1}^l \tm$ can be 
identified with $T(1)$. Vector fields are $(0,1)$-tensors, a 
(commutative) multiplication on $\tm$ is a (symmetric) $(2,1)$-tensor.

The Lie derivative $\lie_X$ with respect to a vector field $X$ is an
$\OO_M$-derivation on the sheaf of $(k,l)$-tensors. 
It is $\lie_X(f)=X(f)$ for functions $f$, $\lie_X (Y)=[X,Y]$ for vector
fields $Y$, 
$\lie_X (Y_1\otimes ... \otimes Y_l) = 
\sum_i Y_1\otimes .. [X,Y_i] .. \otimes Y_l$ for $(0,l)$-tensors, and 
$(\lie_X T)(Y)= \lie_X (T(Y)) - T(\lie_X(Y))$ for $(k,l)$-tensors $T$.
One can always write it explicitly with Lie brackets. Because of the
Jacobi identity the Lie derivative satisfies
$\lie_{[X,Y]}= [\lie_X,\lie_Y]$.
\end{remark}

\begin{definition}
a) An {\it F-manifold} is a triple $(M,\circ,e)$ where $M$ is a complex 
connected manifold of dimension $\geq 1$, $\circ$ is a commutative
and associative $\OO_M$-bilinear multiplication
$\tm\times \tm \to \tm$, $e$ is a global unit field, and the 
multiplication satisfies for any two local vector fields $X,Y$
\begin{eqnarray}
\lie_{X\circ Y} (\circ ) = 
X\circ \lie_Y(\circ ) + Y\circ \lie_X(\circ )\ .
\end{eqnarray}
b) Let $(M,\circ, e)$ be an F-manifold. An {\it Euler field} $E$ {\it 
of weight} $d\in \C$ is a global vector field $E$ which satisfies
\begin{eqnarray}
\lie_E (\circ )= d\cdot \circ\ .
\end{eqnarray}
\end{definition}

\begin{remarks}
i) We do not require that the algebras $(T_pM,\circ, e(p))$ are  
Frobenius algebras (cf. chapter 1). Nevertheless, this is a distinguished 
class. Frobenius manifolds are F-manifolds
\cite{HM}\cite{Man}(I.5). This will be discussed in chapter 5.

\smallskip
ii) Definition 3.2 differs slightly from the definition of F-manifolds
in \cite{HM} by the addition of a global unit field $e$. This unit field
is important, for example, for the definition of $\specan (\tm)$. 
Also, the Euler fields were called weak Euler fields in \cite{HM} in
order to separate them from the Euler fields with stronger properties
of Frobenius manifolds. This is not necessary here.

\smallskip
iii) (3.1) is equivalent to 
\begin{eqnarray}
 [X\circ Y,Z\circ W] - [X\circ Y,Z]\circ W - [X\circ Y,W]\circ Z && 
\nonumber\\
 - X\circ [Y,Z\circ W] + X\circ [Y,Z]\circ W + X\circ [Y,W]\circ Z &&
 \\
 - Y\circ [X,Z\circ W] + Y\circ [X,Z]\circ W + Y\circ [X,W]\circ Z &=&0
\nonumber 
\end{eqnarray}
for any four (local) vector fields $X,Y,Z,W$. (3.2) is equivalent to
\begin{eqnarray}
[E,X\circ Y]-[E,X]\circ Y - X\circ[E,Y] - d\cdot X\circ Y =0 
\end{eqnarray}
for any two (local) vector fields $X,Y$. The left hand sides of 
(3.3) and (3.4) are $\OO_M$--polylinear with respect to 
$X,Y,Z,W$ resp. $X,Y$. Hence they define a $(4,1)$-- resp. $(2,1)$--tensor.
In order to check (3.1) and (3.2) for a manifold with multiplication,
it suffices to check (3.3) and (3.4) for a basis of vector fields.

\smallskip
iv) The unit field $e$ in an F-manifold $(M,\circ,e)$ plays a 
distinguished role. It is automatically nowhere vanishing. It is an
Euler field of weight 0,
\begin{eqnarray}
\lie_e (\circ )= 0\cdot \circ\ ,
\end{eqnarray}
because of (3.1) for $X=Y=e$. So, the multiplication of the F-manifold
is constant along the unit field.

\smallskip
v) An Euler field $E$ of weight $d\neq 0$ is not constant along the 
unit field. But one has for any $d\in \C$
\begin{eqnarray}
[e,E]=d\cdot e \ ,
\end{eqnarray}
because of (3.4) for $X=Y=e$. More generally, in \cite{HM}\cite{Man}(I.5) the
identity
\begin{eqnarray}
[E^{\circ n},E^{\circ m}]=d(m-n)\cdot E^{\circ (m+n-1)} 
\end{eqnarray}
is proved. Chapter 6 will show how intrinsic the notion of an Euler 
field is for an F-manifold.

\smallskip
vi) The sheaf of Euler fields of an F-manifold $(M,\circ, e)$ is a 
Lie subalgebra of $\tm$. If $E_1$ and $E_2$ are Euler fields of weight
$d_1$ and $d_2$, then $c\cdot E_1$ ($c\in \C$) is an Euler field of weight
$c\cdot d_1$, $E_1+E_2$ is an Euler field of weight $d_1+d_2$, and 
$[E_1,E_2]$ is an Euler field of weight 0. The last holds because of
$\lie_{[E_1,E_2]}=[\lie_{E_1},\lie_{E_2}]$ (cf. Remark 3.1 and 
\cite{HM}\cite{Man}(I.5)).

\smallskip
vii) The caustic $\KK$ of an F-manifold is the subvariety of points 
$p\in M$ such that the eigenspace decomposition of $T_pM$ has fewer 
components than for generic points (cf. chapter 2). The caustic is 
invariant with respect to $e$ because of (3.5).
\end{remarks}

\section{Decomposition of F-manifolds and examples}
\setcounter{equation}{0}

\begin{proposition}
The product of two F-manifolds $(M_1,\circ_1,e_1)$ and $(M_2,\circ_2,e_2)$
is an F-manifold
$(M,\circ, e)=(M_1\times M_2, \circ_1 \oplus \circ_2, e_1+e_2)$.

If $E_1$ and $E_2$ are Euler fields on $M_1$ and $M_2$ of the same weight
$d$ then the sum $E_1+E_2$ (of the lifts to $M$) is an Euler field 
of weight $d$ on $M$.
\end{proposition}

{\bf Proof:}
$\tm = \OO_M \cdot pr_1^{-1} \TT_{M_1} \oplus 
\OO_M \cdot pr_2^{-1} \TT_{M_2}$. Vector fields $X_i,Y_i\in 
pr_i^{-1} \TT_{M_i}$, $i=1,2,$ satisfy 
\begin{eqnarray*}
&& X_i\circ Y_i\in pr_i^{-1} \TT_{M_i} , 
\ [X_i,Y_i]\in pr_i^{-1} \TT_{M_i} \ ,\\
&&X_1\circ Y_2 =0,\ [X_1,Y_2]=0\ .
\end{eqnarray*}
This together with (3.3) for vector fields in $\TT_{M_1}$ and for vector 
fields in $\TT_{M_2}$ gives (3.3) for vector fields in 
$pr_1^{-1} \TT_{M_1} \cup pr_2^{-1} \TT_{M_2}$. Because of the 
$\OO_M$--polylinearity then (3.3) holds for all vector fields.
For the same reasons, $E_1+E_2$ satisfies (3.4).
\hfill $\qed$

\begin{theorem}
Let $(M,p)$ be the germ in $p\in M$ of an F-manifold $(M,\circ,e)$.

Then the eigenspace decomposition 
$T_pM = \bigoplus_{k=1}^{l} (T_pM)_k$ of the algebra $T_pM$
extends to a unique decomposition 
$$(M,p) = \prod_{k=1}^{l} (M_k,p)$$ 
of the germ $(M,p)$ into a product
of germs of F-manifolds. These germs $(M_k,p)$ are irreducible germs
of F-manifolds, as already the algebras $T_pM_k \cong (T_pM)_k$ are 
irreducible.

An Euler field $E$ on $(M,p)$ decomposes into a sum of Euler fields
of the same weights on the germs $(M_k,p)$ of F-manifolds.
\end{theorem}

{\bf Proof:}
By Corollary 2.2 the eigenspace decomposition of $T_pM$ extends in some 
neighborhood of $p$ to a decomposition of the tangent bundle into a sum
of multiplication invariant subbundles. First we have to show that these
subbundles and any sum of them are integrable.

Let $\TT_{M,p} = \bigoplus_{k=1}^{l} (\TT_{M,p})_k$ be the according
decomposition of $\TT_{M,p}$ into multiplication invariant free
$\OO_{M,p}$-submodules, and $e=\sum_k e_k$ with 
$e_k \in (\TT_{M,p})_k$. Then $e_k\circ: \TT_{M,p} \to (\TT_{M,p})_k$
is the projection; $e_j\circ e_k = \delta_{jk} e_k$.

\smallskip
{\bf Claim:}
\begin{list}{}{}
\item[i)]
$\lie_{e_k}(\circ )=0\cdot \circ \ ;$
\item[ii)]
$[e_j,e_k]=0\ ;$
\item[iii)]
$[e_j,(\TT_{M,p})_k]\subset (\TT_{M,p})_k\ ;$
\item[iv)]
$[(\TT_{M,p})_j,(\TT_{M,p})_k]\subset (\TT_{M,p})_j + (\TT_{M,p})_k\ .$
\end{list}

\smallskip
{\bf Proof of the claim:}
i) 
\begin{eqnarray*}
\delta_{jk}\cdot \lie_{e_k}(\circ )= \lie_{e_j\circ e_k}(\circ )
= e_j\circ \lie_{e_k}(\circ ) + e_k\circ \lie_{e_j}(\circ )\ .
\end{eqnarray*}
This implies for $j\neq k$ as well as for $j=k$ that
$e_j\circ\lie_{e_k}(\circ)= 0\cdot \circ\ .$ 
Thus $\lie_{e_k}(\circ ) = 0\cdot \circ\ .$

ii) $0= \lie_{e_j}(\circ )(e_k,e_k)= 
[e_j,e_k\circ e_k] - 2e_k\circ [e_j,e_k]$, hence 
$[e_j,e_k]\in (\TT_{M,p})_k$, so for $j\neq k$ we have $[e_j,e_k]=0$, for 
$j=k$ anyway.

iii) Suppose $X=e_k\circ X\in (\TT_{M,p})_k$; then 
$0=\lie_{e_j}(\circ )(e_k,X)= [e_j,X] -e_k\circ[e_j,X].$

iv) Suppose $X\in (\TT_{M,p})_j,\ Y\in (\TT_{M,p})_k,\ k\neq i \neq j$;
then $e_i\circ X=0$ and 
\begin{align*}
0 &= \lie_{e_i\circ X}(\circ ) (e_i,Y) 
   = e_i\circ \lie_X(\circ )(e_i,Y)    \\
  &= e_i\circ [X,e_i\circ Y ] - e_i \circ [X,e_i]\circ Y 
    -e_i\circ[X,Y]\circ e_i            \\ 
  &= - e_i\circ [X,Y] \ .                  \tag{$\qed$}
\end{align*}

\medskip
Claim iv) shows that for any $k$ the subbundle with germs of sections
$\bigoplus_{j\neq k}(\TT_{M,p})_j$ is integrable. According to the 
Frobenius theorem there is a (germ of a) submersion
$f_k:(M,p)\to (\C^{\dim (T_pM)_k},0)$ such that the fibers are the 
integral manifolds of this subbundle. Then 
$\bigoplus_k f_k :(M,p) \to (\C^{\dim M},0)$ is an isomorphism.

The submanifolds $(M_k,p):= ((\bigoplus_{j\neq k}f_j)^{-1}(0),p)$
yield the decomposition 
$(M,p)\cong \prod_{k=1}^{l} (M_k,p)$
with 
$$\TT_{M,p} = \bigoplus_k \OO_{M,p}\cdot  pr_k^{-1} \TT_{M_k,p}
= \bigoplus_k (\TT_{M,p})_k\ .$$

{\bf Claim:}
\begin{list}{}{}
\item[v)] 
$X\circ Y \in pr_k^{-1}\TT_{M_k,p}$
\qquad if $X,Y\in pr_k^{-1}\TT_{M_k,p}\ .$
\item[vi)]
$e_k\circ E \in pr_k^{-1}\TT_{M_k,p}$
\qquad if $E$ is an Euler field.
\end{list}

\smallskip
{\bf Proof of the claim:}
v) $X\circ Y\in \OO_{M,p}\cdot pr_k^{-1}\TT_{M_k,p} $ because this is 
multiplication invariant. Now $X\circ Y \in pr_k^{-1}\TT_{M_k,p}$ 
is true if and only if 
$[Z,X\circ Y]\in (\TT_{M,p})_j$ for any $j$ and any 
$Z\in (\TT_{M,p})_j$; but
\begin{eqnarray*}
[Z,X\circ Y] &=& \lie_Z(X\circ Y)= \lie_{e_j\circ Z}(X\circ Y)  \\
             &=& e_j\circ \lie_Z(X\circ Y) \in (\TT_{M,p})_j\ .
\end{eqnarray*}
vi) Analogously, for any $Z\in (\TT_{M,p})_j$
\begin{eqnarray*}
&& -[Z,e_k\circ E] \\
&=& \lie_{e_k\circ E}(e_j\circ Z)\\ 
&=& \lie_{e_k\circ E}(\circ )(e_j,Z) + \lie_{e_k\circ E}(e_j)\circ Z
+ \lie_{e_k\circ E} (Z)\circ e_j        \\
&=& e_k\circ e_j \circ Z + \lie_{e_k\circ E}(e_j)\circ Z
+ \lie_{e_k\circ E} (Z)\circ e_j   \in  (\TT_{M,p})_j \ .
\end{eqnarray*}
\hfill ($\qed$)

\medskip
Claim v) and vi) show that the multiplication on $(M,p)$ and an 
Euler field $E$ come from multiplication and vector fields on 
the submanifolds $(M_k,p)$ via the decomposition. 
These satisfy (3.3) and (3.4): this is just the restriction of 
(3.3) and (3.4) to $pr_k^{-1}\TT_{M_k,p}$.
\hfill $\qed$

\begin{examples} 
i) $M=\C$ with coordinate $u$ and unit field 
$e={\partial \over \partial u}$ with multiplication $e\circ e=e$ is an
F-manifold. The field $E=u\cdot e =u{\partial \over \partial u}$
is an Euler field of weight 1. The space of all Euler fields of weight
$d$ is $d\cdot E + \C\cdot e$. One has just to check (3.3) and (3.4)
for $X=Y=Z=W=e$ and compare (3.6).

Any 1-dimensional F-manifold is locally isomorphic to an open subset
of this F-manifold $(\C,\circ, e)$. It will be called $A_1$.

\smallskip
ii) From i) and Proposition 4.1 one obtains the F-manifold
$A_1^n = (\C^n,\circ, e)$ with coordinates $u_1,...,u_n$, idempotent
vector fields $e_i = {\partial \over \partial u_i}$, semisimple
multiplication $e_i\circ e_j=\delta_{ij}e_i $, unit field
$e=\sum_i e_i$ and an Euler field $E=\sum_i u_i\cdot e_i$ of weight 1.
Because of Theorem 4.2, the space of Euler fields of weight $d$ is 
$d\cdot E+\sum_i\C\cdot e_i$. 

Also because of Theorem 4.2, any F-manifold $M$ with semisimple multiplication
is locally isomorphic to an open subset of the F-manifold $A_1^n$. 
The induced local coordinates $u_1,...,u_n$ on $M$ are unique up to
renumbering and shift. They are called canonical coordinates, following 
Dubrovin. They are the eigenvalues of a locally defined Euler field 
of weight 1.

\smallskip
iii) Any complex Frobenius manifold is an F-manifold \cite{HM}\cite{Man}(I.5).

\smallskip
iv) Especially, the complex orbit space of a finite Coxeter group
carries the structure of a Frobenius manifold \cite{Du1}\cite{Du2}(Lecture 4). 
The F-manifold structure will be discussed in chapter 18, the 
Frobenius manifold structure in chapter 19.
Here we only give the multiplication and the Euler fields for the 
2-dimensional F-manifolds $I_2(m)$, $m\geq 2$, with 
$I_2(2)=A_1^2,$ $I_2(3)=A_2,$ $I_2(4)=B_2=C_2,$ $I_2(5)=:H_2,$ 
$I_2(6)=G_2$.

\smallskip
The manifold is $M=\C^2$ with coordinates $t_1,t_2$; we denote
$\delta_i := {\partial \over \partial t_i}$. 
Unit field $e$ and multiplication $\circ$ are given by 
$e=\delta_1$ and $\delta_2\circ \delta_2 = t_2^{m-2}\cdot \delta_1$.
An Euler field $E$ of weight 1 is 
$E= t_1\delta_1 + {2\over m} t_2\delta_2$. The space of global Euler fields
of weight $d$ is $d\cdot E+\C\cdot e$ for $m\geq 3$. 
The caustic is $\KK=\{t\in M\ |\ t_2=0\}$ for $m\geq 3$ and $\KK=\emptyset$
for $m=2$. The multiplication is semisimple outside of $\KK$;
the germ $(M,t)$ is an irreducible germ 
of an F-manifold iff $t\in \KK$. One can check all of this by hand.
We will come back to it in Theorem 12.1, when more general results allow
more insight.

\smallskip
v) Another 2-dimensional F-manifold is $\C^2$ with coordinates 
$t_1,t_2$, unit field $e=\delta_1$ and multiplication $\circ$ given by
$\delta_2\circ \delta_2=0$. Here all germs $(M,t)$ are irreducible and
isomorphic. $E_1:=t_1\delta_1$ is an Euler field of weight 1.
Contrary to the cases above with generically semisimple multiplication,
here the space of 
Euler fields of weight $0$ is infinite dimensional, by (3.4):
\begin{eqnarray}
&& \{E \ |\ \lie_E(\circ )=0\cdot \circ\}\nonumber \\
&=& \{E\ |\ [\delta_1, E]=0,\ \delta_2 \circ [\delta_2,E]=0 \}  \\
&=& \{ \varepsilon_1 \delta_1+ \varepsilon_2(t_2)\delta_2\ |\ 
\varepsilon_1 \in \C,\ \varepsilon_2\in \OO_{\C^2}(\C^2),\ 
\delta_1(\varepsilon_2)=0\}\ . \nonumber
\end{eqnarray}

\smallskip
vi) The base space of a semiuniversal unfolding of an isolated 
hypersurface singularity is (a germ of) an F-manifold.
The multiplication was defined first by K. Saito \cite{SaK4}(1.5)
\cite{SaK5}(1.3).
A good part of the geometry of F-manifolds that will be developed in the
next chapters is classical in the case of hypersurface singularities, 
from different points of view. We will discuss this in chapter 16.

\smallskip
vii) Also the base of a semiuniversal unfolding of a boundary 
singularity is (a germ of) an F-manifold, compare chapter 17.
There are certainly more classes of semiuniversal unfoldings which carry
the structure of F-manifolds.
\end{examples}

\section{F-manifolds and potentiality}
\setcounter{equation}{0}

The integrability condition (3.1) for the multiplication in F-manifolds
and the potentiality condition in Frobenius manifolds are closely related.
For semisimple multiplication this has been known previously
(with 6.2 i)$\iff$ii) and e.g. \cite{Hi}(Theorem 3.1)). 
Here we give a general version, requiring neither semisimple 
multiplication nor flatness of the metric.

\begin{remark}
Let $M$ be a manifold with a connection $\nabla$.
The covariant derivative $\nabla_XT$ of a $(k,l)$-tensor with respect to
a vector field is defined exactly as the Lie derivative
$\lie_XT$ in Remark 3.1., starting with the covariant derivatives of 
vector fields.  $\nabla_X$ is a $\OO_M$-derivation on the sheaf
of $(k,l)$-tensors just as $\lie_X$. But $\nabla_X$ is also 
$\OO_M$-linear in $X$, opposite to $\lie_X$. Therefore $\nabla T$ is
a $(k+1,l)$-tensor.
\end{remark}

\begin{theorem}
Let $(M,\circ,\nabla)$ be a manifold $M$ with a commutative and 
associative multiplication $\circ$ on $TM$ and with a torsion free
connection $\nabla$. By definition,
$\nabla \circ (X,Y,Z)$ is symmetric in $Y$ and $Z$.

If the $(3,1)$-tensor $\nabla \circ$ is symmetric in all three
arguments, then the multiplication satisfies for any local 
vector fields $X$ and $Y$
\begin{eqnarray}
\lie_{X\circ Y}(\circ) = X\circ \lie_Y (\circ )+Y\circ\lie_X(\circ)\ .
\end{eqnarray}
\end{theorem}

{\bf Proof:}
$\nabla \circ (X,Y,Z) = \nabla_X(Y\circ Z) - \nabla_X (Y)\circ Z 
            - Y\circ \nabla_X (Z)$
is symmetric in $Y$ and $Z$. The $(4,1)$-tensor 
\begin{eqnarray}
&& (X,Y,Z,W)\nonumber \\ 
&\mapsto &\nabla\circ (X,Y\circ Z,W) + W\circ \nabla\circ (X,Y,Z)
      \nonumber \\
&=& \nabla_X (Y\circ Z\circ W) - \nabla_X(Y) \circ Z\circ W  \\
   && - Y\circ \nabla_X(Z)\circ W - Y\circ Z \circ \nabla_X(W)
      \nonumber
\end{eqnarray}
is symmetric in $Y,Z,W$.
A simple calculation using the torsion freeness of $\nabla$ shows
\begin{eqnarray}
&& (\lie_{X\circ Y} (\circ) - X\circ \lie_Y(\circ) - Y\circ \lie_X(\circ))
(Z,W)   \\
&=& \nabla \circ (X\circ Y,Z,W) - X\circ \nabla \circ (Y,Z,W) 
   - Y\circ \nabla \circ (X,Z,W)     \nonumber \\
&-& \nabla \circ (Z\circ W,X,Y) + Z\circ \nabla \circ (W,X,Y) 
   + W\circ \nabla \circ (Z,X,Y)\ .  \nonumber
\end{eqnarray}
If $\nabla \circ$ is symmetric in all three arguments then
\begin{eqnarray}
\nabla \circ (X\circ Y,Z,W) + Z \circ \nabla \circ (W,X,Y) + 
W\circ \nabla \circ (Z,X,Y) 
\end{eqnarray}
is symmetric in $X,Y,Z,W$ because of the symmetry of the tensor in (5.2).
Then the right hand side of (5.3) vanishes.
\hfill $\qed$

\begin{theorem}
Let $(M,\circ, e, g)$ be a manifold with a commutative and associative
multiplication $\circ$ on $TM$, a unit field $e$, and a metric $g$
(a symmetric nondegenerate bilinear form) on $TM$ which is multiplication 
invariant,
i.e. the $(3,0)$-tensor $A$,
\begin{eqnarray}
A(X,Y,Z) := g(X,Y\circ Z)\ ,
\end{eqnarray}
is symmetric in all three arguments. $\nabla$ denotes the Levi-Civita
connection of the metric. The coidentity $\varepsilon$ is the 1--form
which is defined by $\varepsilon (X) = g(X,e)$. The following 
conditions are equivalent:
\begin{list}{}{}
\item[i)]
$(M,\circ ,e)$ is an F-manifold and $\varepsilon$ is closed.
\item[ii)]
The $(4,0)$-tensor $\nabla A$ is symmetric in all four arguments.
\item[iii)]
The $(3,1)$-tensor $\nabla \circ$ is symmetric in all three arguments.
\end{list}
\end{theorem}

{\bf Proof:}
The Levi-Civita connection satisfies $\nabla g=0$. Therefore
\begin{eqnarray}
&& \nabla A(X,Y,Z,W) \nonumber \\ 
&=& X g(Y,Z\circ W) - g(\nabla_XY, Z\circ W) \nonumber \\
  &&- g(Y,W\circ\nabla_XZ) - g(Y,Z\circ \nabla_XW)  \\
&=& g(Y,\nabla_X (Z\circ W) - W\circ \nabla_XZ - Z\circ \nabla_XW) \nonumber \\
&=& g(Y,\nabla \circ (X,Z,W))\ . \nonumber 
\end{eqnarray}
$g$ is nondegenerate and $\nabla A (X,Y,Z,W)$ is always symmetric in
$Y,Z,W$. Equation (5.6) shows $ii) \iff iii)$. 
Because of $\nabla g=0$ and the torsion freeness
$\nabla_X Y -\nabla_Y X=[X,Y]$,
the 1--form $\varepsilon$ satisfies
\begin{eqnarray}
d \varepsilon (X,Y) &=& X(\varepsilon (Y)) - Y (\varepsilon(X)) 
    - \varepsilon ([X,Y]) \nonumber   \\
&=&  g(Y,\nabla_X e) - g(X,\nabla_Y e) \\
&=& - \nabla A (X,Y,e,e) + \nabla A (Y,X,e,e)\ . \nonumber
\end{eqnarray}
Hence $ii)\Rightarrow d\varepsilon =0$; with Theorem 5.2 this gives
$ii)\Rightarrow i)$. It rests to show $i)\Rightarrow ii)$.

(5.8) and (5.9)  follow from the definition of $\nabla \circ$ and from (5.6),
\begin{eqnarray}
\nabla \circ (X,Y,e) &=& Y\circ \nabla \circ (X,e,e)\ , \\
\nabla A (X,U,Y,e)   &=& \nabla A (X,U\circ Y,e,e) \ .
\end{eqnarray}
One calculates with (5.3) and (5.6) 
\begin{eqnarray}
&& g(e,(\lie_{X\circ Y}(\circ )- X\circ \lie_Y(\circ) -Y\circ \lie_X(\circ))
(Z,W))  \\
&=& \nabla A (X\circ Y,e,Z,W) - \nabla A (Y,X,Z,W) - \nabla A (X,Y,Z,W) 
              \nonumber \\
&-& \nabla A (Z\circ W,e,X,Y) + \nabla A (W,Z,X,Y) + \nabla A (Z,W,X,Y)\ .
\nonumber
\end{eqnarray}
If $i)$ holds then (5.7), (5.9), and (5.10) imply
\begin{eqnarray}
&&\nabla A (Y,X,Z,W) - \nabla A (W,Z,X,Y) \nonumber \\
&=& -\nabla A (X,Y,Z,W) 
  + \nabla A (Z,W,X,Y) \ .
\end{eqnarray}
The left hand side is symmetric in $X$ and $Z$, the right hand side
is skewsymmetric in $X$ and $Z$, so both sides vanish. $\nabla A $ 
is symmetric in all four arguments.
\hfill $\qed$

\begin{lemma}
Let $(M,g,\nabla)$ be a manifold with metric $g$ and Levi-Civita connection
$\nabla$. Then a vector field $Z$ is flat, i.e. $\nabla Z=0$, iff 
$\lie_Z(g)=0$ and the 1--form $\varepsilon_Z:=g(Z,.)$ is closed.
\end{lemma}

{\bf Proof:}
$\nabla$ is torsion free and satisfies $\nabla g=0$. Therefore (cf. (5.7))
\begin{eqnarray}
d \varepsilon_Z (X,Y) 
&=&  g(Y,\nabla_X Z) - g(X,\nabla_Y Z) \ ,\\
\lie_Z (g) (X,Y) &=& g(Y,\nabla_X Z) + g(X,\nabla_Y Z) \ .
\end{eqnarray}
\hfill $\qed$

\begin{remarks} 
a) Let $(M,\circ,e,g)$ satisfy the hypotheses of 
Theorem 5.3 and let $g$ be flat. Then 
condition $ii)$ in Theorem 5.3 
is equivalent to the existence of a local potential
$\Phi\in \OO_{M,p}$ (for any $p\in M$) with 
$(XYZ)\Phi = A(X,Y,Z)$ for any flat local vector fields $X,Y,Z$.

b) In view of this the conditions $ii)$ and $iii)$ in Theorem 5.3 
are called {\it potentiality} conditions.

c) $(M,\circ,e,E,g)$ is a Frobenius manifold if it satisfies the hypotheses
and conditions in Theorem 5.3, if $g$ is flat, 
if $\lie_e(g)=0$ (resp. $e$ is flat),
and if $E$ is an Euler field (with respect to $M$ as F-manifold, 
preferably of weight 1), with $\lie_E(g)=D\cdot g$ for some $D\in \C$.
\end{remarks}

\section{Lagrange property of massive F-manifolds}
\setcounter{equation}{0}

Consider an $n$-dimensional manifold $(M,\circ,e)$ with commutative
and associative multiplication on the tangent bundle and with unit field
$e$. 
Its analytic spectrum $L:=\specan (\TT_M)$ is a subvariety of the 
cotangent bundle $T^{*}M$.
The cotangent bundle carries a natural symplectic structure, given by
the 2--form $d\alpha$. Here $\alpha$ is the canonical 
1--form, which is written as $\alpha=\sum_i y_idt_i$ in local coordinates
$t_1,...,t_n$ for the base and dual coordinates $y_1,...,y_n$ for the 
fibers
($\TT_M\to (\pi_{T^{*}M})_{*}\OO_{T^{*}M},
\ {\partial \over \partial t_i} \mapsto y_i$).

The isomorphism $\aaa: \TT_M \to \pi_{*}\OO_L$ from (2.2) can be 
expressed with $\alpha$ by
\begin{eqnarray}
\aaa(X) = \alpha(\widetilde{X})|_L\ ,
\end{eqnarray}
where $X\in \TT_M$ and $\widetilde{X}$ is any lift of $X$ to a neighborhood 
of $L$ in $T^{*}M$. 

The values of the function $\aaa(X)$ on $\pi^{-1}(p)$ are the
eigenvalues of $X\circ$ on $T_pM$; this follows from Theorem 2.1
and Lemma 1.1.

\begin{definition}
A manifold $(M,\circ, e)$ with commutative and associative multiplication 
on the tangent bundle and with unit field $e$ is {\it massive} if the 
multiplication is generically semisimple.
\end{definition}

Then the set of points where the multiplication is not semisimple is empty
or a hypersurface, the caustic $\KK$ (Proposition 2.4). In the rest
of the paper we will study the local structure of massive F-manifolds
at points where the multiplication is not semisimple.

It is known that the analytic spectrum of a massive Frobenius manifold
is Lagrange (compare \cite{Au} and the references cited there). Theorem 6.2 
together with Theorem 5.3 makes the relations between the 
different conditions transparent.

\begin{theorem}
Let $(M,\circ, e)$ be a massive $n$-dimensional manifold $M$. 

The analytic spectrum $L= \specan (\TT_M)\subset T^{*}M$ is an 
everywhere reduced subvariety. The map $\pi: L\to M$ is finite and flat.
It is a branched covering of degree $n$, branched above the caustic
$\KK$. The following conditions are equivalent.
\begin{list}{}{}
\item[i)]
$(M,\circ, e)$ is a massive F-manifold;
\item[ii)]
At any semisimple point $p\in M-\KK$, the idempotent vector fields
$e_1,...,e_n\in \TT_{M,p}$ commute.
\item[iii)]
$L\subset T^{*}M$ is a Lagrange variety, i.e. $\alpha|L_{reg}$ is closed.
\end{list}
\end{theorem}

{\bf Proof:}
$L-\pi^{-1}(\KK)$ is smooth, $\pi:L-\pi^{-1}(\KK)\to M-\KK$ is a 
covering. $\pi_{*}\OO_L\ (\cong \TT_M)$ is a free $\OO_M$--module, so 
a Cohen--Macaulay $\OO_M$--module and a Cohen--Macaulay ring. 
Then $L$ is Cohen--Macaulay and  everywhere
reduced, as it is reduced at generic points (cf. \cite{Lo2} pp 49--51).

It rests to show the equivalences i)$\iff$ii)$\iff$iii).

i) $\Rightarrow$ ii) follows from Theorem 4.2 and has been discussed in
Example 4.3 ii).

ii) $\Rightarrow$ i) is clear because (3.4) holds everywhere if it holds 
at generic points (in fact, one point suffices).

ii) $\Rightarrow$ iii) We fix canonical coordinates $u_i$ with
${\partial \over \partial u_i}=e_i$ on $(M,p)$ for a point $p\in M-\KK$
and the dual coordinates $x_i$ on the fibers of the contangent bundle
($\TT_M\to (\pi_{T^{*}M})_{*}\OO_{T^{*}M},\ e_i\mapsto x_i$). 
Then locally above
$(M,p)$ the analytic spectrum $L$ is in these coordinates
\begin{eqnarray}
L &\cong& \{(x_j,u_j)\ |\ x_1+...+x_n=1,\ x_ix_j=\delta_{ij}x_j\} \nonumber \\
&=& 
\bigcup_{i=1}^n \{(x_j,u_j)\ |\ x_j=\delta_{ij}\}\ .
\end{eqnarray}
The 1--form $\alpha=\sum x_idu_i$ is closed in $L-\pi^{-1}(\KK)$.
This is open and dense in $L_{reg}$, hence $L$ is a Lagrange variety.

iii) $\Rightarrow$ ii) Above a small neighborhood $U$ of $p\in M-\KK$ 
the analytic spectrum consists of $n$ smooth components $L_k$, $k=1,...,n$,
with $\pi:L_k\stackrel{\cong}{\longrightarrow} U$. 
An idempotent vector field $e_i$ can be lifted to vector fields 
$\widetilde{e_i}$ in neighborhoods $U_k$ of $L_k$ in $T^{*}M$ such that they
are tangent to all $L_k$. 
The commutator $[\widetilde{e_i},\widetilde{e_j}]$ is a lift of the commutator 
$[e_i,e_j]$ in these neighborhoods $U_k$.
\begin{eqnarray*}
\aaa([e_i,e_j]|_U) &=& \bigcup_k 
\bigl(\alpha([\widetilde{e_i},\widetilde{e_j}])\bigr)|_{L_k}\\
&=& \bigcup_k \bigl(\widetilde{e_i}(\alpha(\widetilde{e_j})) - \widetilde{e_j}(\alpha(\widetilde{e_i}))
-d\alpha(\widetilde{e_i},\widetilde{e_j})\bigr)|_{L_k}       \\
&=& \bigcup_k \bigl( \widetilde{e_i}|_{L_k}(\delta_{jk}) 
- \widetilde{e_j}|_{L_k}(\delta_{ik}) 
- d\alpha|_{L_k}(\widetilde{e_i}|_{L_k},\widetilde{e_j}|_{L_k})\bigr) =0 \ .
\end{eqnarray*}
But
$\aaa:\TT_M \to \pi_{*}\OO_L$ is an isomorphism, so $[e_i,e_j]=0$
\hfill $\qed$

\begin{theorem}
a) Let $(M,\circ, e)$ be a massive F-manifold. 

A vector field $E$ is an Euler field of weight $c\in \C$ iff
\begin{eqnarray}
d(\aaa(E))|L_{reg} = c\cdot \alpha|L_{reg} \ .
\end{eqnarray}
b) Let $(M,p)= \prod_{k=1}^l (M_k,p)$ be the decomposition of the germ
of a massive F-manifold into irreducible germs of F-manifolds
$(M_k,\circ, e_k)$.

i) The space of (germs of) Euler fields of weight 0 for $(M,p)$
is the abelian Lie algebra $\sum_{k=1}^l \C\cdot e_k$.

ii) There is a unique continuous function $F:(L,\pi^{-1}(p))
\to (\C,0)$ on the multigerm $(L,\pi^{-1}(p))$ 
which has value 0 on $\pi^{-1}(p)$, is holomorphic on 
$L_{reg}$ and satisfies 
$(dF)|L_{reg}= \alpha|L_{reg}$.

iii) An Euler field of weight $c\neq 0$ for $(M,p)$ exists iff this 
function $F$ is holomorphic. In that case, $c\cdot \aaa^{-1}(F)$
is an Euler field of weight $c$ and 
$\C\cdot \aaa^{-1}(F) + \sum_{k=1}^l \C\cdot e_k$
is the Lie algebra of all Euler fields on the germ 
$(M,p)$.
\end{theorem}

{\bf Proof:} a) It is sufficient to prove this locally for a germ
$(M,p)$ with $p\in M-\KK$. This germ is isomorphic to $A_1^n$. 
A vector field $E=\sum_{i=1}^n\varepsilon_i e_i$, 
$\varepsilon_i\in \OO_{M,p}$, is an Euler field of weight $c$ iff
$d\varepsilon_i = c\cdot du_i$ (Theorem 4.2 and Example 4.3 ii)).

Going into the proof of 6.2 ii) $\Rightarrow$ iii), one sees that this
is equivalent to (6.3).

b) The multigerm $(L,\pi^{-1}(p))$ has $l$ components, the space
of locally constant functions on it has dimension $l$. The function
(multigerm) $F$ exists because $\alpha|L_{reg}$ is closed. This will be 
explained in the next chapter (Lemma 7.1). All 
statements follow now with a). 
\hfill $\qed$

\section{Existence of Euler fields}
\setcounter{equation}{0}

By Theorem 6.2, the analytic spectrum $(L,\lambda)$ of an irreducible 
germ $(M,p)$ of a massive F-manifold is a germ of an (often singular)
Lagrangian variety, and 
$(L,\lambda) \hookrightarrow (T^{*}M,\lambda) \to (M,p)$
is a germ of a Lagrangian map. The paper \cite{Gi2} of Givental is 
devoted to such objects. It contains implicitly many results on 
massive F-manifolds. It will be extremely useful and often cited 
in the following.

The question when a germ of a massive F-manifold has an Euler field of 
weight 1 is reduced by Theorem 6.3 b)iii) to the question when the function
germ $F$ there is holomorphic. Partial answers are given in 
Corollary 7.5 and Lemma 7.6.
We start with a more general situation, as in \cite{Gi2}(chapter 1.1).

Let $(L,0)\subset (\C^N,0)$ be a reduced complex space germ. 
Statements on germs will often be formulated using representatives, 
but they are welldefined for the germs, e.g. "$\alpha|_{L_{reg}}$ 
is closed" for $\alpha\in \Omega_{\C^N,0}^k$.

\begin{lemma}
Let $\alpha\in \Omega_{\C^N,0}^1$ be closed on $L_{reg}$. Then there exists
a unique function germ $F:(L,0)\to (\C,0)$ which is holomorphic on $L_{reg}$,
continuous on $L$ and satisfies $dF|_{L_{reg}}= \alpha|_{L_{reg}}$.
\end{lemma}

{\bf Proof:} 
$(L,0)$ is homeomorphic to a cone as it admits a Whitney stratification.
One can integrate $\alpha$ along paths corresponding to such a cone structure,
starting from 0. One obtains a continuous function $F$ on $L$, which is 
holomorphic on $L_{reg}$ because of $d\alpha|_{L_{reg}}=0$ 
and which satisfies $dF|_{L_{reg}}= \alpha|_{L_{reg}}$.
The unicity of $F$ with value $F(0)=0$ is clear. \hfill $\qed$

\bigskip
Which germs $(L,0)$ have the property that all such function germs are
holomorphic on $(L,0)$? This property has not been studied much. It
can be seen in a line with normality and maximality of complex structures
and is weaker than maximality.

It can be rephrased as $H^1_{Giv}((L,0))=0$.
Here $H^{*}_{Giv}((L,0))$ is the cohomology of the de Rham complex
\begin{eqnarray}
&& (\Omega^{*}_{\C^N,0}/\{\omega \in \Omega_{\C^N,0}^{*}\ |\ 
\omega|_{L_{reg}}=0\},d) \ ,
\end{eqnarray}
which is considered in \cite{Gi2}(chapter 1.1).
We state some known results on this cohomology.

\begin{theorem}
a) {\rm (Poincare-Lemma, \cite{Gi2}(chapter 1.1))} 
If $(L,0)$ is weighted homogeneous with positive weights then
$H^{*}_{Giv}((L,0))=0$. 

\smallskip
b) {\rm (\cite{Va})} Suppose that $(L,0)$ is a germ of a hypersurface
with an isolated singularity, $(L,0)=(f^{-1}(0),0)\subset (\C^{n+1},0)$ 
and $f: (\C^{n+1},0)\to (\C,0)$ is a holomorphic function with an
isolated singularity. Then
\begin{eqnarray}
&& \dim H^n_{Giv}((L,0)) = \mu-\tau \\
&=& 
  \dim\OO_{\C^{n+1},0}/({\partial f\over \partial x_i})
- \dim\OO_{\C^{n+1},0}/(f,{\partial f\over \partial x_i}) \ . \nonumber
\end{eqnarray}

c) {\rm (essentially Varchenko and Givental, \cite{Gi2}(chapter 1.2))}
Let $(L,0)$ be as in b) with $\mu-\tau\neq 0$.
The class $[\eta ]\in H^n_{Giv}((L,0))$ of $\eta\in \Omega^n_{\C^{n+1},0}$
is not vanishing if $d\eta$ is a volume form, i.e. 
$d\eta=hdx_0...dx_n$ with $h(0)\neq 0$.
\end{theorem}

\begin{remarks}
i) The proofs of b) and c) use the Gau{\ss}-Manin connection for isolated
hypersurface singularities.

\smallskip
ii) c) was formulated in \cite{Gi2}(chapter 1.2) only for $n=1$. The missing piece
for the proof for all $n$ was the following fact, which was at that time
only known for $n=1$:
{\it The exponent of a form $hdx_0...dx_n$ is the minimal exponent
iff $h(0)\neq 0$.}

This fact has been established by M. Saito \cite{SaM}(3.11) for all $n$.

\smallskip
iii) By a result of K. Saito \cite{SaK1}, an isolated hypersurface 
singularity $(L,0)=(f^{-1}(0),0)\subset (\C^{n+1},0)$ is weighted 
homogeneous (with positive weights) iff 
$\mu - \tau = 0$.

\smallskip
iv) For us only the case $n=1$ in Theorem 7.2 b)+c) is relevant. 
Proposition 7.4, which is also due to Givental, implies the following:

Of all isolated hypersurface singularities
$(L,0)=(f^{-1}(0),0)\subset (\C^{n+1},0)$
only the curve singularities ($n=1$) turn up as germs of Lagrange 
varieties. These are, of course, germs of Lagrange varieties with
respect to any volume form on $(\C^2,0)$.

\smallskip
v) If $(L,0)\subset((S,0),\omega)$ is the germ of a Lagrange variety 
in a symplectic space $S$ with symplectic form $\omega$, then the class
$[\alpha ]\in H^1_{Giv}((L,0))$ of some $\alpha$ with $d\alpha=\omega$
is independent of the choice of $\alpha$. 
It is called the characteristic class of 
$(L,0)\subset((S,0),\omega)$.

\smallskip
vi) Givental made the conjecture \cite{Gi2}(chapter 1.2): 
{\it Let $(L,0)$ be an $n$-dimensional Lagrange germ. If 
$H^n_{Giv}((L,0))\neq 0$ then $H^1_{Giv}((L,0))\neq 0$ and the 
characteristic class $[\alpha ]\in H^1_{Giv}((L,0))$ is nonzero.}

It is true for $n=1$ because of Theorem 7.2 and Remark 7.3 iii).
Givental sees the conjecture in analogy with a conjecture of Arnold 
which was proved by Gromov 1985 (cf. \cite{Gi2}(chapter 1.2)):
{\it any real closed Lagrange manifold $L\subset T^{*}\R^n$ has 
nonvanishing characteristic class $[\alpha ]\in H^1(L,\R)$.}
\end{remarks}

\begin{proposition}
{\rm (\cite{Gi2}(chapter 1.1))} 
An $n$-dimensional germ $(L,0)$ of a Lagrange variety with embedding 
dimension $\embdim\, (L,0) = n+k$ is a product of a 
$k$-dimensional Lagrange germ $(L',0)$ with 
$\embdim\, (L',0)=2k$ and a smooth $(n-k)$-dimensional 
Lagrange germ $(L'',0)$; here the decomposition of $(L,0)$ 
corresponds to a decomposition
\begin{eqnarray}
((S,0),\omega) \cong ((S',0),\omega') \times ((S'',0),\omega'')
\end{eqnarray}
of the symplectic space germ $(S,0)\supset (L,0)$.
\end{proposition}

{\bf Proof:} 
If $k<n$ then a holomorphic function $f$ on $S$ exists with smooth 
fiber $f^{-1}(0)\supset L$. The Hamilton flow of this function $f$
respects $L$ and the fibers of $f$. The spaces of orbits in 
$f^ {-1}(0)$ and $L$ give a symplectic space germ of dimension
$2n-2$ and in it a Lagrange germ (e.g. \cite{AGV} 18.2).

To obtain a decomposition as in (7.3) one chooses a germ 
$(\Sigma,0)\subset (S,0)$ of a $2n-1$-dimensional submanifold $\Sigma$
in $S$ which is transversal to the Hamilton field $H_f$ of $f$.
There is a unique section $v$ in $(TS)|_\Sigma$ with $\omega (H_f,v)=1$
and $\omega (T_p\Sigma,v)=0$ for $p\in \Sigma$. 
The shift $\widetilde{v}$ 
of $v$ with the Hamilton flow of $f$ forms together with 
$H_f$ a 2-dimensional integrable distribution on $S$, 
because of $0=\lie_{H_f} (\widetilde{v})=[H_f,\widetilde{v}]$.

This distribution is everywhere complementary and orthogonal to the 
integrable distribution whose integral manifolds are the intersections
of the fibers of $f$ with the shifts of $\Sigma$ by the Hamilton flow
of $f$. This  yields a decomposition
$(S,0)\cong (\C^2,0)\times (\Sigma\cap f^{-1}(0),0)$. One can check
that the symplectic form decomposes as required.
If $k<n-1$ one repeats this process.\hfill $\qed$

\begin{corollary}
a) Let $(L,0)$ be an $n$-dimensional Lagrange germ isomorphic to 
$(L',0)\times (\C^{n-1},0)$ as complex space germ. Then 
$(L',0)$ is a plane curve singularity. The characteristic class
$[\alpha ]\in H^1_{Giv}((L,0))$ is vanishing iff $(L',0)$ 
is weighted homogeneous.

b) Let $(L,\lambda)\subset (T^{*}M,\lambda)$ be the analytic spectrum
of an irreducible germ $(M,p)$ of a massive F-manifold.
Suppose $(L,\lambda)\cong (L',0)\times (\C^{n-1},0)$. Then there exists an
Euler field of weight 1 on $(M,p)$ iff $(L',0)$ is weighted homogeneous.
\end{corollary}

{\bf Proof:}
a) Proposition 7.4, Theorem 7.2, and Remark 7.3 iii). 

b) Part a) and Theorem 6.3 b)iii). \hfill $\qed$

\bigskip
In Proposition 20.1 for any plane curve singularity $(L',0)$ irreducible germs
of F-manifolds with analytic spectrum 
$(L,\lambda)\cong (L',0)\times (\C^{n-1},0)$ for some $n$ will be
constructed. So, often there exists no Euler field of weight 1 on 
a germ of a massive F-manifold. On the other hand, the Poincare-Lemma 7.2 a)
and Proposition 7.4 say that an Euler field of weight 1 exists on a germ
of a massive F-manifold $(M,p)$ if the multigerm 
$(L,\pi^{-1}(p))$ of the analytic spectrum
is at all points of $\pi^{-1}(p)$ a product of a
smooth germ and a germ which is weighted homogeneous with positive
weights. Also, we have the following.

\begin{lemma}
Let $M$ be a massive F-manifold and $F:L\to \C$ a continuous function 
with $dF|_{L_{reg}} =\alpha|_{L_{reg}}$. 
Then $\aaa^{-1}(F|(L-\pi^{-1}(\KK)))$ is an Euler field of weight 1 on 
$M-\KK$. It extends to an Euler field on $M$ if $(L,\lambda)$ is at
all points $\lambda \in L$ outside of a subset of codimension 
$\geq 2$ a product of a smooth germ and a germ which is weighted
homogeneous with positive weights.
\end{lemma}

{\bf Proof:}
Suppose $K\subset L$ is a subset of codimension $\geq 2$ with this
property. Then $F$ is holomorphic in $L-K$ because of the 
Poincare-Lemma 7.2 a) and Proposition 7.4. The Euler field extends to
$M-\pi^{-1}(K)$. But $\pi (K)$ has also codimension $\geq 2$.
So the Euler field extends to $M$ (and $F$ is holomorphic on $L$).
\hfill $\qed$

\section{Lyashko-Looijenga maps and graphs of Lagrange maps}
\setcounter{equation}{0}

In this chapter classical facts on Lagrange maps are presented, close
to \cite{Gi2}(chapter 1.3), but slightly more general.
They will be used in chapters 9 -- 11.

Let $L\subset T^{*}M$ be a Lagrange variety (not necessarily smooth)
in the cotangent bundle of an $m$-dimensional connected manifold $M$.
We assume:
\begin{list}{}{}
\item[a)] The projection $\pi:L\to M$ is a branched covering of degree $n$, 
that is, there exists a subvariety $K{\subset}M$ such that 
$\pi:L-\pi^{-1}(K)\to M-K$ is a covering of degree $n$
($\pi:L\to M$ is not necessarily flat).
\item[b)] There exists a {\it generating function} $F:L\to \C$, that is, 
a continuous function which is holomorphic on $L_{reg}$ with 
$dF|_{L_{reg}}=\alpha|_{L_{reg}}$
(locally such a function exists by Lemma 7.1).
\end{list}

Such a function $F$ will be fixed. It can be considered as a 
multivalued function on $M-K$; the 1--graph of this multivalued function
is $L-\pi^{-1}(K)$.

The {\it Lyashko-Looijenga map} 
$\Lambda = (\Lambda_1,...,\Lambda_n):M\to \C^n$ of $L\subset T^{*}M$ and 
$F$ is defined as follows:
for $q\in M-K$, the roots of the unitary polynomial 
$z^n+\sum_{i=1}^n \Lambda_i(q)z^{n-i}$ are the values of $F$ on 
$\pi^{-1}(q)$. It extends to a holomorphic map on $M$ because 
$F$ is holomorphic on $L_{reg}$ and continuous on $L$.

The {\it reduced Lyashko-Looijenga} map 
$\Lambda^{(red)} = (\Lambda_2^{(red)},...,\Lambda_n^{(red)}):
M\to \C^{n-1}$ of $L\subset T^{*}M$ and 
$F$ is defined as follows:
for $q\in M-K$, the roots of the unitary polynomial 
$z^n+\sum_{i=2}^n \Lambda_i^{(red)}(q)z^{n-i}$ are the values of $F$ on 
$\pi^{-1}(q)$, shifted by their center 
$-{1\over n}\Lambda_1(q) = 
{1\over n}\sum_{\lambda \in \pi^{-1}(q)}F(\lambda)$.
It also extends to a holomorphic map on $M$.
Its significance will be discussed after the Remarks 8.2.

The {\it front} $\Phi_L$ of $L\subset T^{*}M$ and $F$ is the image 
$\Im (F,pr)\subset \C\times M$ of $(F,pr):L\to \C\times M$.
It is the zero set of the polynomial 
$z^n+\sum_{i=1}^n\Lambda_i\cdot z^{n-i}$. 
So, it is an analytic hypersurface even if $F$ is not holomorphic on all
of $L$.

Following Teissier (\cite{Tei}(2.4, 5.5), \cite{Lo2}(4.C)), 
the {\it development} 
$\widetilde{\Phi_L} \subset \P T^{*}(\C\times M)$ of this hypersurface
$\Phi_L$ in $\C\times M$ is defined as the closure in 
$\P T^{*}(\C\times M)$ of the set of tangent hyperplanes at the smooth 
points of $\Phi_L$. It is an analytic subvariety and a Legendre
variety with respect to the canonical contact structure on 
$\P T^{*}(\C\times M)$.

The map
\begin{eqnarray}
\C\times T^{*}M\to \P T^{*}(\C\times M),\ 
(c,\lambda)\mapsto ((dz-\lambda)^{-1}(0),(c,p)) &&
\end{eqnarray}
($\lambda \in T_p^{*}M$ and 
$dz-\lambda \in T_{(c,p)}^{*}(\C\times M) 
\cong T_c^{*}\C\times T_p^{*}M$)
identifies $\C\times T^{*}M$ with the open subset in $\P T^{*}(\C\times M)$
of hypersurfaces in the tangent spaces which do not contain
$\C\cdot {\partial \over \partial z}$. 
The induced contact structure on $\C\times T^{*}M$ is given by the 
1--form $dz-\alpha$.

The following fact is wellknown. It is one form of appearance of the 
relation between Lagrange and Legendre maps
(e.g. \cite{AGV} 18.--20.). To check it, one has to consider $F$ as a 
multivalued function on $M-K$ and $\Phi_L$ as its graph.

\begin{proposition}
The embedding $\C\times T^{*}M\hookrightarrow \P T^{*}(\C\times M)$
identifies the graph $\Im (F,id)\subset \C\times T^{*}M$ of $F:L\to \C$
with the development $\widetilde{\Phi_L}$ of the front $\Phi_L$.
\end{proposition}

\begin{remarks}
It has some nontrivial consequences.

i) The polynomial $z^n+\sum_{i=1}^n\Lambda_iz^{n-i}$ has no multiple
factors and the branched covering $\Phi_L\to M$ has degree $n$: over any 
point $p\in M-K$ the varieties $L$ and $\widetilde{\Phi_L}$ have $n$
points and the points of $\Phi_L$ have $n$ tangent planes; so, over a 
generic point $p\in M- K$ also $\Phi_L$ has $n$ points.

\smallskip
ii) The graph $\Im(F,id)\cong \widetilde{\Phi_L}$ is an analytic variety
even if $F$ is not holomorphic on all of $L$.

\smallskip
iii) The composition of maps
$\widetilde{\Phi_L}\stackrel{\cong}{\longrightarrow}
\Im(F,id)\stackrel{pr}{\longrightarrow} L$ is a bijective morphism.
It is an isomorphism iff $F$ is holomorphic on $L$. Also, the 
continuous map $L\to \Phi_L$ is a morphism iff $F$ is holomorphic on $L$.

\smallskip
iv) The Lagrange variety $L\subset T^{*}M$ together with the values of $F$
at one point of each connectivity component of $L$ and any of the following 
data determine each other uniquely: the front $\Phi_L$, the development
$\widetilde{\Phi_L}$, the Lyashko-Looijenga map $\Lambda$, the generating
function $F$ as a multivalued function on the base $M$.
\end{remarks}

To motivate the reduced Lyashko-Looijenga map, we have to talk about 
Lagrange maps and their isomorphisms (\cite{AGV}18., \cite{Gi2}3.1).

A Lagrange map is a diagram $L\hookrightarrow (S,\omega)\to M$ 
where $L$ is a Lagrange variety in a symplectic manifold $(S,\omega)$ and
$S\to M$ is a Lagrange fibration. 
An isomorphism between two Lagrange maps is given by an isomorphism of the
Lagrange fibrations which maps one Lagrange variety to the other.

An automorphism of $T^{*}M\to M$ as Lagrange fibration which fixes the base
is given by a shift in the fibers,
\begin{eqnarray*}
T^{*}M\to T^{*}M,\ \lambda \mapsto \lambda + dS,
\end{eqnarray*}
where $S:M\to \C$ is holomorphic (\cite{AGV}18.5). So, regarding
$T^{*}M\to M$ as a Lagrange fibration means to forget the 0-section and the
1--form $\alpha$, but to keep the Lagrange fibration and the class
$\alpha+\{dS\ |\ S:M\to \C\mbox{ holomorphic}\}$ of 1--forms.

\begin{corollary}
Let $L\hookrightarrow T^{*}M\to M $ be as above (satisfying the assumptions
a) and b)) with $l$ connectivity
components and points $\lambda_1,...,\lambda_l$, one in each connectivity
component.
The data in i) -- iii) are equivalent.
\begin{list}{}{}
\item[i)] $L\hookrightarrow T^{*}M\to M $ as a Lagrange map and the
differences $F(\lambda_i)-F(\lambda_j)\in \C$ of values of $F$,
\item[ii)] the generating function $F$ modulo addition of a function on 
the base,
\item[iii)] the reduced Lyashko-Looijenga map 
$\Lambda^{(red)}:M\to \C^{n-1}$.
\end{list}
\end{corollary}

{\bf Proof:}
i) $\Rightarrow$ ii): Integrating the 1--forms in 
$\alpha+\{dS\ |\ S:M\to \C\mbox{ holomorphic}\}$ gives ii).

ii) $\Rightarrow$ iii) Definition of $\Lambda^{(red)}$.

iii) $\Rightarrow$ i): $(0,\Lambda^{(red)})=(0,\Lambda^{(red)}_2,...,
\Lambda^{(red)}_n):M\to \C^n$ is the Lyashko-Looijenga map of a 
Lagrange variety in $T^{*}M$ which differs from the original Lagrange 
variety only by the shift of $d({-1\over n}\Lambda_1)$ in the fibers.
$(0,\Lambda^{(red)})$ determines this Lagrange variety and a generating 
function for it because of Remark 8.2 iv). One recovers i).
\hfill $\qed$

\section{Miniversal Lagrange maps and F-manifolds}
\setcounter{equation}{0}

The notion of a miniversal germ of a Lagrange map is central in Givental's
paper \cite{Gi2}. 
We need a slight generalization to multigerms, taking a semilocal viewpoint.

Let $L\subset T^{*}M$ be a Lagrange variety with finite projection 
$\pi :L\to M$. 
The {\it germ at the base point }$p\in M$ of 
$L\hookrightarrow T^{*}M\to M$ is the diagram 
\begin{eqnarray*}
(L,\pi^{-1}(p))\hookrightarrow (T^{*}M,T^{*}_p M)\to (M,p).
\end{eqnarray*}
Here $(L,\pi^{-1}(p))$ is a multigerm. $(T^{*}M,T^{*}_p M)\to (M,p)$
is the cotangent bundle of the germ $(M,p)$; it is a germ in the base, but 
not in the fiber.

For this diagram the morphisms
\begin{eqnarray}
\OO_{M,p}\oplus \TT_{M,p} \to \pi_{*}(\OO_L)_p , \ 
(c,X)\mapsto (c+\alpha(\widetilde{X}))|_{L}        
\end{eqnarray}
and
\begin{eqnarray}
\C\oplus T_pM &\to& \pi_{*}(\OO_L)_p/\mmm_p\cdot \pi_{*}(\OO_L)_p\ , \\ 
(c,X) &\mapsto& (c+\alpha(\widetilde{X}))
|\pi_{*}(\OO_L)_p/\mmm_p\cdot \pi_{*}(\OO_L)_p       \nonumber
\end{eqnarray}
are welldefined. Here $\widetilde{X}$ is in both cases a lift of $X$
to $T^{*}M$. These morphisms are not invariants of the diagram as a 
{\it germ at the base point $p\in M$ of a Lagrange map} 
because the identification of the Lagrange fibration with the 
cotangent bundle of $(M,p)$ is unique only up to shifts in the fibers and
only the class of 1--forms 
$\alpha+\{dS\ |\ S:M\to \C\mbox{ holomorphic}\}$
is uniquely determined (cf. chapter 8).

But being an isomorphism or epimorphism in (9.1) and (9.2) are clearly 
properties of the germ at $p$ of the Lagrange map.

\begin{definition}
The germ at $p\in M$ of $L\hookrightarrow T^{*}M\to M$
as a Lagrange map is called {\it miniversal} ({\it versal}) if the morphism
in (9.2) is an isomorphism (epimorphism) (cf. \cite{Gi2}(chapter 1.3)).
\end{definition}

We are only interested in the case of a flat projection $\pi :L\to M$.
Wellknown criteria of flatness for finite maps give the next lemma.

\begin{lemma}
The following conditions are equivalent.
\begin{list}{}{}
\item[i)] The germ at $p\in M$ of $L\hookrightarrow T^{*}M\to M$
as a Lagrange map is miniversal with flat projection 
$\pi:(L,\pi^{-1}(p))\to (M,p)$,
\item[ii)] it is miniversal with $\deg \pi =1+\dim M$,
\item[iii)] the morphism in (9.1) is an isomorphism,
\item[iv)] the Lagrange map is miniversal at all points in a 
neighborhood of $p\in M$.
\end{list}
\end{lemma}

\begin{example}
A miniversal germ at a base point of a Lagrange map with not flat projection
$\pi:L\to M$ is given by the germ at $0\in \C^2$ of the Lagrange fibration
$\C^4\to \C^2$,  $(y_2,y_3,t_2,t_3)\mapsto (t_2,t_3)$ with 
$\omega= dy_2dt_2+dy_3dt_3$ and by the Lagrange variety $L$ which is the union
of two appropriate planes and which is defined by the ideal
\begin{eqnarray}
(y_2,y_3)\cap (y_2-t_2,y_3-t_3)=(y_2,y_3)\cdot (y_2-t_2,y_3-t_3)\ .
\end{eqnarray}
\end{example}

\smallskip
Now let $M$ be a massive $n$-dimensional F-manifold with analytic spectrum
$L\subset T^{*}M$. Then
$(L,\pi^{-1}(p))\hookrightarrow (T^{*}M,T^{*}_p M)\to (M,p)$ is for any 
$p\in M$ a versal, but not a miniversal germ at the base point 
$p\in M$ of a Lagrange map. But there is a miniversal one.

$pr_e:(M,p)\to (M^{(r)},p^{(r)})$ denotes the germ of the fibration at
$p$ whose fibers are the orbits of the unit field $e$. The fiberwise linear
function on $T^{*}M$ which corresponds to $e$ is called $y_1$. Its 
Hamilton field $\widetilde{e}:=H_{y_1}$ is a lift of $e$ to $T^{*}M$.
It leaves the hypersurface $y_1^{-1}(1)\subset T^{*}M$ and the Lagrange variety
$L\subset y_1^{-1}(1)$ invariant. The orbits of $\widetilde{e}$ in 
$y_1^{-1}(1)$ form a germ of 
a $2n-2$-dimensional symplectic manifold with a Lagrange fibration,
which can be identified with the cotangent bundle
\begin{eqnarray*}
(T^{*}M^{(r)},T^{*}_{p^{(r)}}M^{(r)})\to (M^{(r)},p^{(r)})\ .
\end{eqnarray*}
But this identification is only unique up to shifts in the fibers.
The orbits of $\widetilde{e}$ in $L$ form a Lagrange variety
$L^{(r)}\subset T^{*}M^{(r)}$. 

The germ at $p^{(r)}\in M^{(r)}$ of the diagram
$L^{(r)}\hookrightarrow T^{*}M^{(r)}\to M^{(r)}$ 
is unique up to isomorphism of {\it germs in the base of Lagrange maps}.
It will be called the {\it restricted Lagrange map} of the germ $(M,p)$ of 
the F-manifold $M$. An explicit description will be given in the proof 
of the next result.

\begin{theorem}
a) The restricted Lagrange map of the germ $(M,p)$ of a massive F-manifold
is miniversal with flat projection $\pi^{(r)}:L^{(r)}\to M^{(r)}$.

b) It determines the germ $(M,p)$ of the F-manifold uniquely.

c) Any miniversal germ at a base point of a Lagrange map 
$L' \hookrightarrow T^{*}M' \to M' $ with flat projection $L' \to M' $ 
is the restricted Lagrange map of a germ of a massive F-manifold.
\end{theorem}

{\bf Proof:}
a) In order to be as explicit as possible we choose coordinates
$t=(t_1,t')=(t_1,...,t_n):(M,p)\to (\C^n,0)$ with $e(t_1)=1$.
The dual coordinates on $T^{*}M$ are $(y_1,...,y_n)=(y_1,y')=y$.
The multiplication is given by
${\partial \over \partial t_i}\circ {\partial\over \partial t_j}
= \sum_k a_{ij}^k(t') {\partial \over \partial t_k}$
and the analytic spectrum $L$ is
\begin{eqnarray}
L\cong \{(y,t)\in \C^n\times (\C^n,0)\ |\ 
y_1=1,\ y_iy_j=\sum_k a_{ij}^k(t')y_k\}\ .
\end{eqnarray}
The restricted Lagrange map is represented by the Lagrange fibration
\begin{eqnarray}
\C^{n-1}\times (\C^{n-1},0) \to (\C^{n-1},0),\ (y',t')\mapsto t'
\end{eqnarray}
with canonical 1--form $\alpha':= \sum_{i\geq 2}y_idt_i$ and by the 
Lagrange variety
\begin{eqnarray}
&&  \{(y',t')\in \C^{n-1}\times (\C^{n-1},0) \ | \nonumber   \\
&&  \quad 
y_iy_j=a_{ij}^1(t')+\sum_{k\geq 2}a_{ij}^k(t')y_k \mbox{ for } i,j\geq 2\}
\cong L^{(r)}\ .
\end{eqnarray}
The equations for the Lagrange variety in (9.6) show that the morphism in 
(9.1) for this Lagrange map with fixed canonical 1--form 
$\alpha'$ is an isomorphism. This implies that the restricted Lagrange map
for $(M,p)$ is miniversal with flat projection.

\smallskip
b)+c) Any miniversal germ at a base point $p'\in M' $ of a Lagrange map
$L'\hookrightarrow T^{*}M' \to M' $ with flat projection 
$\pi':L'\to M'$ can be represented by a Lagrange fibration as in (9.5)
and a Lagrange variety as in (9.6).

Defining $L$ by (9.4) and $M:= \C\times (M',p')$ with 
$e:={\partial \over \partial t_1}$, one obtains an F-manifold with unit field
$e$ and analytic spectrum $L$.

It rests to show that this does not depend on the way in which the 
Lagrange fibration is identified with the cotangent bundle of 
$(\C^{n-1},0)$ in (9.5).
But one sees easily that a shift in the fibers of (9.5) of the type 
$y_i\mapsto y_i+{\partial S\over \partial t_i}$ for some holomorphic 
function $S:(\C^{n-1},0)\to \C$ on the base 
corresponds only to a change of the
coordinate fields ${\partial\over \partial t_2},...,
{\partial\over \partial t_n}$ and the coordinate $t_1$ in $M$
and thus to a shift of the section $\{0\}\times(M',p') $ in 
$M\to (M',p')$. It does not affect $L$ and the multiplication on $(M,p)$.
\hfill $\qed$

\bigskip

Let $(M,p)$ be a germ of a massive F-manifold. The germ
\begin{eqnarray}
\GG_{p^{(r)}} := \{ X\in \TT_{M,p}\ |\ [e,X]=0\}
\end{eqnarray}
is a free $\OO_{M^{(r)},p^{(r)}}$-module of rank $n$. 
It is an $\OO_{M^{(r)},p^{(r)}}$-algebra because of 
$\lie_e(\circ)=0\cdot \circ$.
The functions $\aaa(X)$ for $X\in \GG_{p^{(r)}}$ are invariant with
respect to $\widetilde{e}$ and induce holomorphic functions on $L^{(r)}$.
One obtains a map
\begin{eqnarray}
\aaa^{(r)}:\GG_{p^{(r)}} \to (\pi^{(r)}_{*}\OO_{L^{(r)}})_{p^{(r)}}\ .
\end{eqnarray}

\begin{lemma}
$\aaa^{(r)}$ is an isomorphism of $\OO_{M^{(r)},p^{(r)}}$-algebras.
\end{lemma}

{\bf Proof:}
The isomorphism $\aaa:\TT_{M,p}\to (\pi_{*}\OO_L)_p$ maps the $e$-invariant
vector fields in $(M,p)$ to the $\widetilde{e}$-invariant functions in 
$(\pi_{*}\OO_L)_p$.\hfill $\qed$

\bigskip
This isomorphism $\aaa^{(r)}$ is closely related to (9.1) for the
restricted Lagrange map of $(M,p)$:

An isomorphism as in (9.1) requires the choice of a 1--form for its Lagrange
fibration. 
The choice of a function $t_1:(M,p)\to (\C,0)$ with $e(t_1)=1$
yields such a 1--form: the 1--form which is induced by $\alpha-dt_1$
($\alpha-dt_1$ on $T^{*}M$ is $\widetilde{e}$-invariant and vanishes
on $\widetilde{e}$ and induces a 1--form on the space of 
$\widetilde{e}$-orbits of $y_1^{-1}(1)$).

The choice of $t_1$ also yields an isomorphism
\begin{eqnarray}
&&\OO_{M^{(r)},p^{(r)}} \oplus \TT_{M^{(r)},p^{(r)}} \nonumber \\
&\to& 
\OO_{M^{(r)},p^{(r)}} \cdot e \oplus \{X\in \GG_{p^{(r)}}\ |\ X(t_1)=0\}
=\GG_{p^{(r)}}\ .
\end{eqnarray}
One sees with the proof of Theorem 9.4 a) that the composition of this 
isomorphism with $\aaa^{(r)}$ gives the isomorphism in (9.1)
for the restricted Lagrange map of $(M,p)$ 
(the germ $(M,p)$ in (9.1) is here $(M^{(r)},p^{(r)})$).

\section{Lyashko-Looijenga map of an F-manifold}
\setcounter{equation}{0}

\begin{definition}
Let $(M,\circ,e)$ be a massive $n$-dimensional F-manifold with analytic
spectrum $L\subset T^{*}M$.

a) A {\it generating function} $F$ for $(M,\circ,e)$ is a generating 
function for $L$, that is, a continuous function $F:L\to \C$ which is
holomorphic on $L_{reg}$ with $dF|_{L_{reg}} = \alpha|_{L_{reg}}$.

b) Let $F$ be a generating function for $(M,\circ,e)$.
\begin{list}{}{}
\item[i)] 
The {\it bifurcation diagram} $\BB \subset M$ of $(M,\circ,e,F)$ is the set of 
points $p\in M$ such that $F$ has less than $n$ different values on 
$\pi^{-1}(p)$.
\item[ii)]
The {\it Lyashko-Looijenga map} 
$\Lambda=(\Lambda_1,...,\Lambda_n):M\to \C^n$ of $(M,\circ,e,F)$ is the
Lyashko-Looijenga map of $F$ as generating function for $L\subset T^{*}M$
(cf. chapter 8).
\item[iii)]
The {\it discriminant} $\DD\subset M$ of $(M,\circ,e,F)$ is 
$\DD:=\Lambda_n^{-1}(0)$.
\end{list}
\end{definition}

The discriminant will be discussed in chapter 11. A generating function for 
an F-manifold exists locally (Lemma 7.1), but not necessarily globally.

A holomorphic generating function $F$ corresponds to an Euler field
$E:=\aaa^{-1}(F)$ of weight 1 (Theorem 6.3); then the values of $F$ on
$\pi^{-1}(p)$, $p\in M$, are the eigenvalues of $E\circ:T_pM\to T_pM$.
The objects $\BB,\ \Lambda,\ \DD$ of Definition 10.1 b) are welldefined
for $(M,\circ,e,E)$ if $E$ is such an Euler field.

The restriction of $\pi:L\to M$ to the complement of the caustic $\KK$ is a 
covering $\pi:\pi^{-1}(M-\KK)\to M-\KK$ of degree $n$, and 
$\pi^{-1}(M-\KK)$ is smooth. Hence there a generating function $F$ is 
holomorphic and corresponds to an Euler field $E$ on $M-\KK$.
Results and examples about the extendability of $E$ to $M$ are given in
Lemma 7.6 and Theorem 20.4.

The bifurcation diagram $\BB$ of $(M,\circ,e,F)$ contains the caustic
$\KK$. The caustic is a hypersurface or empty (Proposition 2.4) and 
invariant with respect to the unit field $e$ (Remark 3.3 vii)).
The bifurcation diagram has the same properties: the restriction of $F$ to 
an open set $U\subset M-\KK$ with canonical coordinates $(u_1,...,u_n)$
corresponds to an Euler field $E=\sum (u_i+r_i)e_i$ for some $r_i\in \C$,
and the bifurcation diagram is the hypersurface
\begin{eqnarray*}
\BB\cap U = U\cap \{u\ |\ u_i+r_i=u_j+r_j \mbox{ for some }i\neq j\}\ .
\end{eqnarray*}
It is invariant with respect to $e$ because of $e(u_i-u_j)=0$.

\medskip
The Lyashko-Looijenga map for the F-manifold $A_1^n= (\C^n,\circ,e)$
(Example 4.3 ii)) with Euler field $E=\sum u_ie_i$ and 
Euler field-function $F:=\aaa(E)$ is 
\begin{eqnarray}
\Lambda^{(n)}:\C^n\to \C^n,\ u \mapsto ((-1)^i\sigma_i(u))
\end{eqnarray}
($\sigma_1(u),...,\sigma_n(u)$ are the symmetric polynomials).
The group of automorphisms of the F-manifold $A_1^n$ which respect the 
Euler field $E$ is the symmetric group $S_n$ which permutes the 
coordinates $u_1,...,u_n$. The map $\Lambda^{(n)}$ is the quotient map
for this group. It is branched along the bifurcation diagram
$\BB=\{u\ |\ u_i=u_j \mbox{ for some }i\neq j\}$. 
The image $\Lambda^{(n)}(\BB)$ is the hypersurface
\begin{eqnarray}
\DD^{(n)} := \{a\in \C^n\ |\ z^n+\sum a_iz^{n-i}
\mbox{ has multiple roots } \} \subset \C^n . &&
\end{eqnarray}
The restriction 
$\Lambda^{(n)}:\C^n - \BB \to \C^n- \DD^{(n)}$ 
induces an F-manifold structure on $\C^n - \DD^{(n)}$, 
with unit field
\begin{eqnarray}
e^{(n)}:= d \Lambda^{(n)}(e)= -n{\partial\over \partial a_1}
- \sum_{i\geq 2} (n-i+1)a_{i-1}{\partial\over \partial a_i}
\end{eqnarray}
and Euler field 
\begin{eqnarray}
E^{(n)} := d \Lambda^{(n)}(E) = \sum_i i a_i {\partial\over \partial a_i}\ .
\end{eqnarray}
This F-manifold $(\C^n-\DD^{(n)},\circ, e^{(n)})$ will be denoted by
$A_1^n/S_n$.

\begin{theorem}
Let $(M,\circ,e)$ be a massive F-manifold with generating function 
$F:L\to \C$ and Lyashko-Looijenga map $\Lambda:M\to \C^n$.

Then $\Lambda^{-1}(\DD^{(n)})=\BB$ and $d \Lambda(e)=e^{(n)}$.
The restriction $\Lambda :M-\BB\to \C^n-\DD^{(n)}$ is an immersion
and locally an isomorphism of F-manifolds. It maps the Euler field
$\aaa^{-1}(F|_{M-\BB})$ on $M-\BB$ to the Euler field
$E^{(n)}$.
\end{theorem}

{\bf Proof:} In $M-\KK$ the multiplication is semisimple and locally
the values of the generating function are canonical coordinates.
The map $\Lambda$ factors on $M-\BB$ locally into an isomorphism to $A_1^n$
and into the map $\Lambda^{(n)}$.
\hfill $\qed$

\bigskip

The most important part of Theorem 10.2 is that 
$\Lambda:M-\BB\to \C^n-\DD^{(n)}$ is locally biholomorphic. The following
statements for germs will also be useful.

\begin{lemma}
Let $(M,p)= \prod_{k=1}^l (M_k,p)$ be a germ of a massive F-manifold
with analytic spectrum $L$ and with 
decomposition into irreducible germs $(M_k,p)$ of dimension
$n_k$, $\sum n_k=n$.

\smallskip
a) There exists precisely one generating function on the multigerm  
$(L,\pi^{-1}(p))$ for any choice of its values on 
$\pi^{-1}(p)=\{ \lambda_1,...,\lambda_l\} $.

\smallskip
b) Choose a function $t_1:(M,p)\to \C$ with $e(t_1)=1$. 

The values of a generating function for the points in $L$ above an orbit
of $e$ are of the form $t_1+\ constant$.

The entry $\Lambda_i$ of a Lyashko-Looijenga map
$\Lambda : (M,p)\to \C^n$ is a polynomial of degree $i$ in $t_1$ 
with coefficients in $\{g\in \OO_{M,p}\ |\ e(g)=0\}$ and leading 
coefficient $(-1)^i\left( {m\atop i}\right)$. 

\smallskip
c) Choose representatives $M_k$ for the germs $(M_k,p)$ and 
Lyashko-Looijenga maps $\Lambda^{[k]}:M_k\to \C^{n_k}$. 
Then the function 
$\Lambda = (\Lambda_1,...,\Lambda_n):\prod_kM_k \to \C^n$ 
which is defined by
\begin{eqnarray}
z^n+ \sum_{i=1}^n \Lambda_i z^{n-i} = 
\prod_{k=1}^l (z^{n_k}+ \sum_{i=1}^{n_k}\Lambda_i^{[k]}z^{n-i})
\end{eqnarray}
is a Lyashko-Looijenga map for the representative $\prod M_k$ of the 
germ $(M,p)$. Any Lyashko-Looijenga map for $(M,p)$ is of this type.
\end{lemma}

{\bf Proof:}
a) Lemma 7.1. 

b) It suffices to prove the first part for an orbit of $e$ in $M-\KK$.
There the generating function comes from an Euler field.
(3.5) and (3.6) imply $\lie_e(E\circ)=id$. The values of $F$ are the
eigenvalues of $E\circ$.

c) $\Lambda^{[k]}$ corresponds to an Euler field 
$E^{[k]}$ (at least) on $M_k-\BB_{M_k}$. The sum $\sum_k E^{[k]}$ is an 
Euler field on $\prod (M_k-\BB_{M_k})$ by Proposition 4.1.  
The corresponding generating function extends to $\prod M_k$ 
and has the given $\Lambda$ 
as Lyashko-Looijenga map. The last statement follows with a).
\hfill $\qed$.

\bigskip

Consider the projection $pr_e:(M,p)\to (M^{(r)},p^{(r)})$ whose fibers are 
the orbits of $e$ (chapter 9). The $e$-invariant hypersurfaces $\BB$ and 
$\KK$ project to hypersurfaces in $M^{(r)}$, which are called 
the restricted bifurcation diagram $\BB^{(r)}$ and the restricted
caustic $\KK^{(r)}$. In chapter 8 the restricted Lagrange map was defined
as the germ at $p^{(r)}\in M^{(r)}$ of a Lagrange map 
$L^{(r)}\hookrightarrow T^{*}M^{(r)}\to M^{(r)}$. Because of Corollary 8.3 
the notion of a reduced Lyashko-Looijenga map is welldefined for the 
restricted Lagrange map (independently of the identification of the
Lagrange fibration with the cotangent bundle of $M^{(r)}$).

The space of orbits of the field $e^{(n)}$ (formula (10.3)) in $\C^n$ can be 
identified with $\{a\in \C^n\ |\ a_1=0\} = \{0\}\times \C^{n-1}\cong \C^{n-1}$
and is equipped with the coordinate system $(a_2,...,a_n)=a'$.
The projection to this orbit space is denoted by $pr^{(n)}:\C^n\to \C^{n-1}$,
the image of $\DD^{(n)}$ is 
\begin{eqnarray}
\DD (A_{n-1}):= \{ a'\in \C^{n-1} \ |
\ z^n+\sum_{i=2}^n a_iz^{n-i} \mbox{ has multiple
roots}\}
\end{eqnarray}
(it is isomorphic to the discriminant of the singularity or F-manifold
$A_{n-1}$, cf. chapter 16).

\begin{corollary}
Let $(M,p)$ be the germ of a massive F-manifold, $F$ a generating function,
$\Lambda$ ($\Lambda^{(red)}$) the (reduced) Lyashko-Looijenga map of
$(M,\circ,e,F)$.

$\Lambda^{(red)}:(M,p)\to\C^{n-1}$ 
is constant along the orbits of $e$. The induced map 
$\Lambda^{(red)(r)}: (M^{(r)},p^{(r)})\to \C^{n-1}$ is a reduced 
Lyashko-Looijenga map for the restricted Lagrange map. The following 
diagram commutes, the diagonal morphism is $\Lambda^{(red)}=
pr^{(n)}\circ \Lambda = \Lambda^{(red)(r)}\circ pr_e$,
\begin{eqnarray}
\begin{CD}
(M,p)      @>{\Lambda}>> \C^n                          \\
@VV{pr_e}V              @VV{pr^{(n)}}V                 \\
(M^{(r)},p^{(r)}) @>{\Lambda^{(red)(r)}}>>    \C^{n-1} 
\end{CD}
\end{eqnarray}
The restriction
\begin{eqnarray}
\Lambda^{(red)(r)}: M^{(r)}-\BB^{(r)} \to \C^{n-1}-\DD (A_{n-1})
\end{eqnarray}
is locally biholomorphic.
\end{corollary}

{\bf Proof:}
$\Lambda^{(red)}$ is constant along the orbits of $e$ because of Lemma 
10.3 b).
The formulas (9.4), (9.5), (9.6) show that $\Lambda^{(red)(r)}$ is a 
reduced Lyashko-Looijenga map for the restricted Lagrange map. 
The rest follows from Theorem 10.2.
\hfill $\qed$

\section{Discriminant of an F-manifold}
\setcounter{equation}{0}

Let $(M,\circ, e,F)$ be a massive $n$-dimensional F-manifold with a 
generating function $F:L\to \C$ and Lyashko-Looijenga map 
$\Lambda=(\Lambda_1,...,\Lambda_n):M\to \C^n$; the discriminant of 
$(M,\circ,e,F)$ is the hypersurface $\DD=\Lambda_n^{-1}(0)\subset M$
(Definition 10.1).

If $F$ is holomorphic and $E={\bf a}^{-1}(F)$ is its Euler field 
then $\Lambda_n=(-1)^n\cdot \det (E\circ)$, 
and the discriminant is the set of points
where the multiplication with $E$ is not invertible.

By definition of $\Lambda_n$, the discriminant is $\DD=\pi(F^{-1}(0))$.
Theorem 11.1 will give an isomorphism between $F^{-1}(0)$ and the development
$\widetilde{\DD}\subset \P T^{*}M$ of $\DD$.

We need an identification of subsets of $T^{*}M$ and $\P T^{*}M$. The 
fiberwise linear function on $T^{*}M$ which corresponds to $e$ is called
$y_1$. The canonical map
\begin{eqnarray}
y_1^{-1}(1) \longrightarrow \P T^{*}M
\end{eqnarray}
identifies $y_1^{-1}(1)\subset T^{*}M$ with the open subset in $\P T^{*}M$
of hypersurfaces in the tangent spaces of $M$ which do not contain 
$\C\cdot e$.

The restriction to $y_1^{-1}(1)$ of the canonical 1--form $\alpha$ on $T^{*}M$
gives the contact structure on $y_1^{-1}(1)$ which is induced by the 
canonical contact structure on $\P T^{*}M$.

\begin{theorem}
Let $(M,\circ,e,F)$ and $F^{-1}(0)\subset L \subset y_1^{-1}(1) \subset 
T^{*}M$ be as above. 

The canonical map $y_1^{-1}(1)\hookrightarrow \P T^{*}M$ identifies 
$F^{-1}(0)$ with the development $\widetilde{\DD}\subset \P T^{*}M$ of
the discriminant $\DD$.
\end{theorem}

{\bf Proof:}
We want to make use of the discussion of fronts and graphs (chapter 8)
and of the restricted Lagrange map (chapter 9).
It is sufficient to consider the germ $(M,p)$ for some $p\in \DD$.
W choose coordinates $(t_1,...,t_n)$ and $(y_1,...,y_n)$ as in the proof of 
Theorem 9.4. The generating function $F$ for 
$L\hookrightarrow T^{*}M \to M$ takes the form
\begin{eqnarray}
F(y',t)=t_1+F^{(r)}(y',t')\ ,
\end{eqnarray}
$F^{(r)}$ is a generating function of $L^{(r)}$ with respect to (9.6), (9.5), 
and $\alpha'= \sum_{i\geq 2}y_idt_i$. The isomorphism with a sign
\begin{eqnarray}
(-t_1,pr_e):(M,p)\to (\C\times M^{(r)},0\times p^{(r)})
\end{eqnarray}
maps the discriminant to the front $\Im (F^{(r)},\pi^{(r)})$ of 
$L^{(r)}$ and $F^{(r)}$.

The development of this front is identified with the graph
$\Im (F^{(r)},id) \subset \C\times T^{*}M^{(r)}$ of $F^{(r)}:L^{(r)}\to \C$,
by Proposition 8.1 and the embedding (cf. formula (8.1))
\begin{eqnarray}
\C\times T^{*}M^{(r)} &\hookrightarrow& \P T^{*}(\C\times M^{(r)}),\\ 
(-t_1,\lambda)&\mapsto &
\bigl( (dt_1+\lambda)^{-1}(0),(-t_1,\pi^{(r)}(\lambda))\bigr)\ .\nonumber
\end{eqnarray}
But (11.1), (11.3), and (11.4) together yield also an isomorphism 
$y_1^{-1}(1)\to \C\times T^{*}M^{(r)}$, which maps $F^{-1}(0)$ to this graph.
\hfill $\qed$

\begin{remarks}
i) In the proof only the choice of $t_1$ is essential. It is equivalent to
various other choices: the choice of a section of the projection
$(M,p)\to (M^{(r)},p^{(r)})$, the choice of a 1--form $\alpha'$ for the
Lagrange fibration in the restricted Lagrange map.

\smallskip
ii) Fixing such a choice of $t_1$, one obtains together with $\alpha'$ 
and $F^{(r)}$ a Lyashko-Looijenga map 
$\Lambda^{(r)}:(M^{(r)},p^{(r)})\to \C^{n-1}$ for $F^{(r)}$.
Then (11.3) identifies the entry $\Lambda_n$ of the Lyashko-Looijenga
map $\Lambda:(M,p)\to \C^n$ with the polynomial 
$(-t_1)^n+\sum_{i=2}^n \Lambda_i^{(r)}(-t_1)^{n-i}$.

\smallskip
iii) $F^{-1}(0)\subset L$ is not an analytic hypersurface of $L$ at points of 
$L_{sing}$ where $F$ is not holomorphic. But Theorem 11.1 shows that it is 
everywhere a subvariety of $L$ of pure codimension 1.
Examples with not holomorphic $F$ will be given in chapter 20.

\smallskip
iv) Let $(M,p)=\prod_{k=1}^l(M_k,p)$ be the decomposition of a germ $(M,p)$
into irreducible germs of F-manifolds. A Lyashko-Looijenga map
$\Lambda:(M,p)\to \C^n$ corresponds to Lyashko-Looijenga maps
$\Lambda^{[k]}:(M_k,p)\to \C^{n_k}$ for the irreducible germs, in a way which
was described in Lemma 10.3 c). 
Especially $\Lambda_n=\prod_k \Lambda_{n_k}^{[k]}$,
and the germ $(\DD,p)$ of the discriminant is
\begin{eqnarray}
(\DD,p)= \bigcup_k \left( \prod_{j<k}M_j\times 
(\Lambda_{n_k}^{[k]})^{-1} (0) \times \prod_{j>k}M_j, p\right)\ ,
\end{eqnarray}
the union of  products of smooth germs with the discriminants for the 
irreducible germs
(of course, $((\Lambda_{n_k}^{[k]})^{-1} (0),p)=\emptyset$ for some or all $k$
is possible).

\smallskip
v) The development $\widetilde{\DD}$ of the discriminant $\DD$ gives the
tangent hyperplanes to $\DD$. 
Let $(M,p)=\prod_{k=1}^l(M_k,p)$ be as in iv), and 
$\pi^{-1}(p)=\{\lambda_1,...,\lambda_l\}$. 

Theorem 11.1 says that the tangent hyperplanes to $(\DD,p)$ are those 
hyperplanes $\lambda_k^{-1}(0)\subset T_pM$ for which $F(\lambda_k)=0$.

Especially, if $l=1$ and $F(\lambda_1)=0$, then 
$\lambda_1^{-1}(0)\subset T_pM$ is the nilpotent subalgebra of $T_pM$ and 
the unique tangent hyperplane to $(\DD,p)$.
The general case fits with Lemma 1.1 iv) and (11.5).
 
\smallskip
vi) $\alpha|_{L_{reg}} = dF|_{L_{reg}}$ shows immediately that 
$F^{-1}(0) \subset y_1^{-1}(1)$ is a Legendre subvariety.

\smallskip
vii) If $\lambda \in F^{-1}(0)\cap T^{*}_pM$, there is a canonical projection
from the Legendre germ $(F^{-1}(0),\lambda)$ to one component of the Lagrange
multigerm $(L^{(r)},{\pi^{(r)}}^{-1}(p^{(r)}))$. It is a bijective morphism.
It is an isomorphism iff $F$ is holomorphic at $(L,\lambda)$
(see Remark 8.2 iii)).
\end{remarks}

One can recover the multiplication on a massive $n$-dimensional F-manifold
$M$ from the unit field $e$ and a discriminant $\DD$ if the orbits of $e$
are sufficiently large. To make this precise, we introduce the 
following notion.

\begin{definition}
A massive F-manifold $(M,\circ,e,F)$ with generating function $F$ is in 
{\it standard form} if there exists globally a projection $pr_e:M\to 
M^{(r)}$ to a manifold $M^{(r)}$ such that 
\begin{list}{}{}
\item[$\alpha)$] 
the fibers are the orbits of $e$ (and thus connected),
\item[$\beta)$]  
they are with their affine linear structure isomorphic to open 
(connected) subset of $\C$,
\item[$\gamma)$]
the projection $pr_e:\DD \to M^{(r)}$ is a branched covering of degree $n$.
\end{list}
\end{definition}

\begin{remarks}
i) If $(M,\circ,e,F)$ is a massive F-manifold with generating function $F$
and properties $\alpha)$ and $\beta)$ then
\begin{eqnarray}
( -{1\over n}\Lambda_1,pr_e):M\to \C\times M^{(r)}
\end{eqnarray}
is an embedding because of $e(-{1\over n}\Lambda_1)=1$ (Lemma 10.3 b)).

The F-manifold $M$ can be extended uniquely to an F-manifold isomorphic to
$\C\times M^{(r)}$. Also the generating function $F$ can be extended.
The discriminant of this extended F-manifold satisfies $\gamma)$ because
of $\DD=\Lambda_n^{-1}(0)$ and Lemma 10.3 b).

ii) For $(M,\circ,e,F)$ as in i) the coordinate $t_1:=-{1\over n}\Lambda_1$
is distinguished and, up to addition of a constant, even independent of 
the choice of $F$. Nevertheless it does not seem to have beautiful properties:
In the case of the simplest 3-dimensional irreducible germs of 
F-manifolds, $A_3,\ B_3,\ H_3,$ it is not part of the coordinate system
of a nice normal form (chapter 18). Using the data in \cite{Du1} one can 
also check that $-{1\over n}\Lambda_1$ is not a flat coordinate of the
Frobenius manifolds $A_3,\ B_3,\ H_3$.
\end{remarks}

\begin{corollary}
Let $(M,\circ,e,F)$ be a massive F-manifold with generating function $F$ 
and in standard form.

a) The branch locus of the branched covering $pr_e:\DD \to M^{(r)}$
is $\DD_{sing}$, the set $pr_e(\DD_{sing})$ of critical values is the 
restricted bifurcation diagram $\BB^{(r)}=pr_e(\BB)$.

b) The union of the shifts of $F^{-1}(0)$ with the Hamilton field 
$\widetilde{e}=H_{y_1}$ is the analytic spectrum $L\subset T^{*}M$.

c) The data $(M,\circ,e,F)$ and $(M,e,\DD )$ are equivalent.
\end{corollary}

{\bf Proof:} a)
Theorem 11.1 implies that all tangent hyperplanes to $\DD$ are 
transversal to the unit field. Therefore the branch locus is 
only $\DD_{sing}$.

b)+c) $\DD$ and $e$ determine $F^{-1}(0)\subset y_1^{-1}(1)$ because of
Theorem 11.1. The union of shifts of $F^{-1}(0)$ with the Hamilton
field $\widetilde{e}=H_{y_1}$ is finite of degree $n$ over $M$ because 
of $\gamma)$ and it is contained in $L$, so it is $L$.
The function $F:L\to \C$ is determined by
$F^{-1}(0)\subset L$ and by the linearity of $F$ along the orbits of 
$\widetilde{e}$.
\hfill $\qed$

\bigskip
Theorem 11.1 and Corollary 11.5 b) give one possibility to recover the
multiplication of a massive F-manifold $(M,\circ,e,F)$ in standard form
from the discriminant $\DD$ and the unit field $e$. A more elementary way
is the following.

\begin{corollary}
Let $(M,\circ,e,F)$ be a massive $n$-dimensional F-manifold 
with generating function $F$ and in standard form. The multiplication can be
recovered from the discriminant $\DD$ and the unit field $e$ in the 
following way.
The multiplication is semisimple outside of the bifurcation diagram
$\BB=pr_e^{-1}(pr_e(\DD_{sing}))$. For a point $p\in M-\BB$, 
the idempotent vectors $e_i(p)\in T_pM$ with 
$e_i(p)\circ e_j(p)=\delta_{ij}e_i(p)$ are uniquely determined by i) and ii):
\begin{list}{}{}
\item[i)]
$e(p)=\sum_{i=1}^n e_i(p)$,
\item[ii)] the multigerm $(\DD,\DD\cap pr_e^{-1}(pr_e(p)))$ has exactly $n$
tangent hyperplanes; their shifts to $T_pM$ with $e$ are the 
hyperplanes
$\bigoplus_{i\neq k}\C\cdot e_i \subset T_pM$,
$k=1,...,n$.
\end{list}
\end{corollary}

{\bf Proof:}
Remark 11.2 v).  \hfill $\qed$

\section{2-dimensional F-manifolds}
\setcounter{equation}{0}

The only 1-dimensional germ of an F-manifold is $A_1$ (Example 4.3 i)).
The class of 3-dimensional germs of massive F-manifolds is already
vast. Examples and a partial classification will be given in chapter 20.
But the classification of 2-dimensional germs of F-manifolds is nice.

\begin{theorem}
a) The only germs of 2-dimensional massive F-manifolds are, up to isomorphism, 
the germs $I_2(m)$, $m\in \N_{\geq 2}$, with $I_2(2)=A_1^2$, $I_2(3)=A_2$, 
$I_2(4)=B_2$, $I_2(5)=:H_2$, $I_2(6)=G_2$, from Example 4.3 iv): 

The multiplication on $(M,p)=(\C^2,0)$ with coordinates $t_1,t_2$ and
$\delta_i:={\partial \over \partial t_i}$ is given by
$e:=\delta_1$ and $\delta_2\circ \delta_2=t_2^{m-2}\cdot \delta_1$.
An Euler field of weight 1 is $E=t_1\delta_1 + {2\over m}t_2\delta_2$.
Its discriminant is $\DD =\{t\ |\ t_2^2-{4\over m^2}t_2^m =0\}$.
$I_2(m)$ is irreducible for $m\geq 3$ with caustic and bifurcation 
diagram $\KK=\BB= \{t\ |\ t_2=0\}$. The space of Euler fields of weight
$d$
is $d\cdot E + \C\cdot e$ for $m\geq 3$.

b) The only germ of a 2-dimensional not massive F-manifold is the germ
$(\C^2,0)$ from Example 4.3 v) with multiplication given by $e:=\delta_1$
and $\delta_2\circ \delta_2=0$. The caustic is empty. An Euler field
of weight 1 is $E=t_1\delta_1$. The space of all Euler fields of weight 0 is 
$\{\varepsilon_1\delta_1 + \varepsilon_2(t_2)\delta_2\ |\ \varepsilon_1\in \C,
\ \varepsilon_2(t_2)\in \C\{t_2\}\}$.
\end{theorem}

{\bf Proof:}
a) Givental \cite{Gi2}(1.3, p 3253) classified the 1-dimensional miniversal
germs of Lagrange maps with flat projection. Together with Theorem 9.4
this yields implicitly the classification of the 2-dimensional irreducible
germs of massive F-manifolds. But we can recover this in a simple way
and we need to be more explicit.

Let $(M,p)$ be a 2-dimensional germ of a massive F-manifold with projection
$pr_e:(M,p)\to (M^{(r)},p^{(r)})$ to the space of orbits of $e$.
There is a unique generating function $F:(L,\pi^{-1}(p))\to (\C,0)$.
Its bifurcation diagram $\BB\subset M$ and restricted bifurcation diagram
$\BB^{(r)}\subset M^{(r)}$ are the hypersurfaces $\BB=pr_e^{-1}(p^{(r)})$
and $\BB^{(r)}=\{p^{(r)}\}$. 

By Corollary 10.4, the reduced Lyashko-Looijenga map 
$\Lambda^{(red)(r)}:(M^{(r)},p^{(r)})\to (\C,0)$ of the 
restricted Lagrange map is a cyclic branched covering of some order
$\widetilde{m}$, and the Lyashko-Looijenga map $\Lambda:M\to \C^2$ is also a 
cyclic branched covering of order $\widetilde{m}$, branched along $\BB$.

Because of Corollary 8.3 and Theorem 9.4, this branching order $\widetilde{m}$ 
of $\Lambda^{(red)(r)}$ determines the germ $(M,p)$ of the F-manifold
up to isomorphism. It rests to determine the allowed $\widetilde{m}$ and
explicit formulas for the F-manifolds.

Now consider the manifolds $\C^2$ with multiplication on $T\C^2$ given by
$e=\delta_1$ and $\delta_2\circ \delta_2=t_2^{m-2}\cdot \delta_1$. 
The analytic spectrum 
$L\subset T^{*}\C^2$ is
\begin{eqnarray}
L=   \{(y_1,y_2,t_1,t_2)\ |\ y_1=1,\ y_2\cdot y_2=t_2^{m-2}\}\ .
\end{eqnarray}
It is an exercise to see that $L$ is a Lagrange variety with generating 
function $F=t_1+{2\over m}y_2t_2$ with respect to 
$\alpha=y_1dt_1+y_2dt_2$, i.e. one has $(\alpha-dF)|_{L_{reg}}=0$.
Then this gives an F-manifold and
$E=\aaa^{-1}(F)=t_1\delta_1+{2\over m}t_2\delta_2$ is an Euler field
of weight 1. The Lyashko-Looijenga map $\Lambda$ is
\begin{eqnarray}
\Lambda:\C^2\to \C^2,\ (t_1,t_2)\mapsto (-2t_1,t_1^2-{4\over m^2}t_2^m)\ .
\end{eqnarray}
It is branched along $\BB=\{t\ |\ t_2=0\}$ of degree $m$.

So $I_2(m)$ is the desired F-manifold for any branching order 
$m=\widetilde{m}\geq 2$.
The same calculation yields for $m=\widetilde{m}=1$ the F-manifold 
$A_1^2/S_2$ on $\C^2-\{t\ |\ t_2=0\}$ (chapter 10) with meromorphic
multiplication along $\{t\ |\ t_2=0\}$.

\smallskip
b) Let $(M,p)=(\C^2,0)$ be the germ of a 2-dimensional not massive F-manifold
with $e=\delta_1$. Then $\lie_e(\circ)=0$ and $[e,\delta_2]=0$ imply
$\lie_e(\delta_2\circ\delta_2)=0$. 
Hence $\delta_2\circ\delta_2=\beta(t_2)\delta_1+\gamma(t_2)\delta_2$
for some $\beta(t_2),\gamma(t_2)\in \C\{t_2\}$.

The field $\widetilde{\delta_2}:=\delta_2-{1\over 2}\gamma(t_2)\delta_1$
satisfies
$\widetilde{\delta_2}\circ\widetilde{\delta_2}=
(\beta(t_2)+{1\over 4}\gamma(t_2)^2)\delta_1$
and $[\delta_1,\widetilde{\delta_2}]=0$. Changing coordinates we may 
suppose $\delta_2=\widetilde{\delta_2}$, $\gamma(t_2)=0$.

The analytic spectrum is 
$L=\{(y_1,y_2,t_1,t_2)\ |\ y_1=1,\ y_2\cdot y_2=\beta(t_2)\}$.
The F-manifold is not massive, hence $\beta(t_2)=0$.
One checks with (3.3) and (3.4) easily that this multiplication gives an 
F-manifold and that the space of Euler fields is as claimed.
\hfill $\qed$

\bigskip
Let us discuss the role which 2-dimensional germs of F-manifolds can 
play for higher dimensional massive F-manifolds. The set
\begin{eqnarray}
\DD^{(n,3)} &:=& \{ a\in\C^n\ |\ z^n+\sum_{i=1}^na_iz^{n-i}
\mbox{ has a root}\\
&&\mbox{\qquad of multiplicity } \geq 3\}
\subset \DD^{(n)} \subset \C^n \nonumber
\end{eqnarray}
is an algebraic subvariety of $\C^n$ 
of codimension 2 (see the proof of Proposition 2.3).
Given a massive F-manifold $M$, the space
\begin{eqnarray}
\KK^{(3)}:= \{ p\in M\ |\ \PP(T_pM)\succ (3,1,...,1)\} \subset \KK\subset M
\end{eqnarray}
of points $p$ such that $(M,p)$ does not decompose into 1- and 2-dimensional
germs of F-manifolds is empty or an analytic subvariety (Propostion 2.3).

\begin{theorem}
Let $(M,\circ,e)$ be a massive F-manifold with generating function $F$.

a) $F$ is holomorphic on $\pi^{-1}(M-\KK^{(3)})$ and gives rise to an
Euler field of weight 1 on $M-\KK^{(3)}$.

b) If $\codim \KK^{(3)} \geq 2$ then $F$ is holomorphic on $L$ and 
$E=\aaa^{-1}(F)$ is an Euler field of weight 1 on M.

c) $\KK^{(3)}\subset \Lambda^{-1}(\DD^{(n,3)})$, and 
$\overline{\Lambda^{-1}(\DD^{(n,3)})-\KK^{(3)}}$ is analytic of pure 
codimension 2. Thus 
$\codim \KK^{(3)}\geq 2 \iff \codim \Lambda^{-1}(\DD^{(n,3)})\geq 2$.

d) The restriction of the Lyashko-Looijenga map
\begin{eqnarray*}
\Lambda:M-\Lambda^{-1}(\DD^{(n,3)}) \to \C^n-\DD^{(n,3)}
\end{eqnarray*}
is locally a branched covering, branched along 
$\BB-\Lambda^{-1}(\DD^{(n,3)})$.

If $p\in \BB-\Lambda^{-1}(\DD^{(n,3)})$ and 
$\Lambda (p)\in \DD^{(n)} -\DD^{(n,3)}$ are smooth points of the 
hypersurfaces $\BB$ and $\DD^{(n)}$ and if there the branching order is $m$,
then $(M,p)$ is the germ of an F-manifold of type $I_2(m)\times A_1^{n-2}$.
\end{theorem}

{\bf Proof:}
a)+b) Each germ $(L,\lambda_k)$ of the analytic spectrum $(L,\pi^{-1}(p))
=(L,\{\lambda_1,...,\lambda_l\})$ of a reducible germ 
$(M,p)=\prod_k (M_k,p)$ is the product of a smooth germ 
with the analytic spectrum of $(M_k,p)$. The analytic spectrum of 
$I_2(m)$ ($m\geq 2$) is isomorphic to 
$(\C,0)\times (\{y_2,t_2)\ |\ y_2^2=t_2^{m-2}\},0)$. One applies
Lemma 7.6.

c) A Lyashko-Looijenga map of $I_2(m)$ is a cyclic branched covering of 
order $m$, branched along the bifurcation diagram. This together with
Lemma 10.3 c) implies that locally around a point 
$p\in M-\KK^{(3)}$ the fibers of the Lyashko-Looijenga map 
$\Lambda:M\to \C^n$ are finite. Therefore 
$\codim_M (\Lambda^{-1}(\DD^{(n,3)}),p)=\codim_{\C^n}(\DD^{(n,3)}) =2$.

d) $\Lambda$ determines the multiplication of the F-manifold $M$
(Theorem 10.2). One uses this, Lemma 10.3 c) and properties of $I_2(m)$.
\hfill $\qed$

\bigskip

Many interesting F-manifolds, e.g. those for hypersurface singularities,
boundary singularities, finite Coxeter groups (chapters 16, 17, 18),
satisfy the property
$\codim \KK^{(3)}\geq 2$ and have an Euler field of weight 1.

\section{Logarithmic vector fields}
\setcounter{equation}{0}

K. Saito \cite{SaK3} introduced the notions of logarithmic vector fields
and free divisors. Let $H\subset M$ be a reduced hypersurface in an
$n$-dimensional manifold $M$. The sheaf $\Der_M (\log H)\subset \tm$
of logarithmic vector fields consists of those holomorphic vector fields
which are tangent to $H_{reg}$. The hypersurface $H$ is a free divisor
if $\Der_M(\log H)$ is a free $\OO_M$-module of rank $n$.

The results in this chapter are not really new. They had been established
in various generality by Bruce \cite{Bru}, Givental \cite{Gi2}(chapter 1.4),
Lyashko \cite{Ly1}\cite{Ly3}, K. Saito \cite{SaK4}\cite{SaK5}, Terao
\cite{Ter}, and Zakalyukin \cite{Za} as results for hypersurface singularities,
boundary singularities or miniversal Lagrange maps.
But the formulation using the multiplication of F-manifolds is 
especially nice.

\begin{theorem}
Let $(M,\circ,e)$ be a massive F-manifold with Euler field $E$ of weight 1,
generating function $F=\aaa(E)$ and discriminant 
$\DD = (\det (E\circ ))^{-1} (0) = \pi (F^{-1}(0))$.

a) The discriminant is a free divisor with $\Der_M (\log \DD )=E\circ \tm$.

b) The kernel of the map
\begin{eqnarray}
\aaa_\DD :\tm \to \pi_{*} \OO_{F^{-1}(0)},\ 
X\mapsto \aaa(X)|_{F^{-1}(0)}
\end{eqnarray}
is 
\begin{eqnarray}
\ker \aaa_\DD = E\circ \tm = \Der_M (\log \DD )\ .
\end{eqnarray}
\end{theorem}

{\bf Proof:}
a) $E\circ \tm$ is a free $\OO_M$-module of rank $n$. Therefore a) follows 
from (13.2).

b) The $\OO_M$-module $\pi_{*}\OO_{F^{-1}(0)}$ has support $\DD$.
(13.2) holds in $M-\DD$. The set $\DD\cap \BB$ has codimension 2 in $M$.
Hence (13.2) holds in $M$ if it holds in $\DD-\BB$.

Let $p\in \DD-\BB$. We choose a small neighborhood $U$ of $p$ with 
canonical coordinates $u_1,...,u_n$ centered at $p$, with 
$\DD\cap U = \{u\ |\ u_1=0\}$ and with Euler field 
$E=u_1e_1+\sum_{i\geq 2}(u_i+r_i)e_i$ for some $r_i\in\C-\{0\}$.
With the notation of the proof of Theorem 6.2 ii)$\Rightarrow$iii) we have 
$\alpha=\sum x_idu_i$,
\begin{eqnarray*}
F^{-1}(0)\cap \pi^{-1}(U) = \{(x,u)\ |\ x_j=\delta_{1j},u_1=0\}\ ,
\end{eqnarray*}
and for any vector field $X=\sum\xi_ie_i\in \tm(U)$
\begin{eqnarray*}
\aaa(X)|F^{-1}(0)\cap \pi^{-1}(U)=\xi_1(0,u_2,...,u_n)\ .
\end{eqnarray*}
Therefore
\begin{eqnarray}
(\ker \aaa_\DD)_p &=& \OO_{M,p} \cdot u_1e_1
\oplus \bigoplus_{i=2}^n \OO_{M,p}\cdot e_i \nonumber \\
&=& E\circ \TT_{M,p} = \Der_{M,p}(\log \DD) \ .
\end{eqnarray}
\hfill $\qed$

\begin{remark}
One can see 13.1 a) in a different way: there is a criterion of K. Saito
\cite{SaK3}(Lemma (1.9)). To apply it, one has to show 
\begin{eqnarray}
[E\circ \tm,E\circ \tm] \subset E\circ \tm\ .
\end{eqnarray}
With (3.1) and (3.2) one calculates for any two (local) vector fields $X,Y$
\begin{eqnarray}
[E\circ X,E\circ Y]=E\circ ([X,E\circ Y]-[Y,E\circ X]-E\circ[X,Y])\ .
\end{eqnarray}
\end{remark}

\medskip
In the rest of this chapter $(M,\circ, e,E)$ will be a massive F-manifold
which is equipped with an Euler field $E$ of weight 1 and which is in 
standard form (Definition 11.3). $pr_e:M\to M^{(r)}$ is the projection
to the space of orbits of $e$.
The sheaf of $e$-invariant vector fields
\begin{eqnarray}
\GG:=\{X\in (pr_e)_{*}\tm\ |\ [e,X]=0\}
\end{eqnarray}
is a free $\OO_{M^{(r)}}$-module of rank $n$. Because of $\lie_e(\circ)=0$ 
it is also an $\OO_{M^{(r)}}$-algebra.

\begin{theorem}
Let $(M,\circ,e,E)$ be a massive F-manifold with Euler field $E$ of weight 1
and in standard form.

a) 
\begin{eqnarray}
(pr_e)_{*}\tm = \GG \oplus (pr_e)_{*} (E\circ \tm)\ .
\end{eqnarray}

b) The kernel of the map 
\begin{eqnarray}
(pr_e)_{*} \aaa_\DD : (pr_e)_{*}\tm &\to& (pr_e\circ \pi)_{*} 
\OO_{F^{-1}(0)},\\
 X&\mapsto& \aaa(X)|_{F^{-1}(0)}   \nonumber
\end{eqnarray}
is $(pr_e)_{*} (E\circ \tm)$. The restriction
\begin{eqnarray}
(pr_e)_{*} \aaa_\DD : \GG \to (pr_e\circ \pi)_{*} 
\OO_{F^{-1}(0)}
\end{eqnarray}
is an isomorphism of $\OO_{M^{(r)}}$-algebras.
\end{theorem}

{\bf Proof:}
a) follows from b).

b) The kernel of $(pr_e)_{*} \aaa_\DD $ is $(pr_e)_{*} (E\circ \tm)$
because of Theorem 13.1 b).
The F-manifold in standard form has a global restricted Lagrange map 
$L^{(r)}\hookrightarrow T^{*}M^{(r)} \to M^{(r)}$ (the identification of its
Lagrange fibration with $T^{*}M^{(r)} \to M^{(r)}$ is unique only up to shifts
in the fibers). 

The canonical projection $F^{-1}(0)\to L^{(r)}$ is 
bijective (Corollary 11.5 b)), and then an isomorphism because $F$ is
holomorphic.
It induces an isomorphism $(pr_e\circ \pi)_{*}\OO_{F^{-1}(0)}\cong 
(\pi^{(r)})_{*}\OO_{L^{(r)}}$. The composition with (13.9) is the 
isomorphism
\begin{eqnarray}
\aaa^{(r)}:\GG \to (\pi^{(r)})_{*} \OO_{L^{(r)}}
\end{eqnarray}
from Lemma 9.5.
\hfill $\qed$

\bigskip
Theorem 13.1 and Theorem 13.3 are translations to F-manifolds of statements
in \cite{SaK4}(1.6)\cite{SaK5}(1.7) for hypersurface singularities. 
In fact, K. Saito
used essentially (13.9) to define the multiplication on $\GG$ for 
hypersurface singularities. 

The arguments in Lemma 13.4 and Theorem 13.5 are 
due to Lyashko \cite {Ly1}\cite{Ly3} and Terao \cite {Ter}, see also 
Bruce \cite{Bru}. 

Again $(M,\circ, e,E)$ is a massive F-manifold with Euler field of weight 1
and in standard form. We choose a function $t_1:M\to \C$ with $e(t_1)=1$
(e.g. $t_1=-{1\over n}\Lambda_1$, cf. Lemma 10.3 b)).
This choice simplifies the formulation of the results in Lemma 13.4.
The vector fields in $\TT_{M^{(r)}}$ will be identified with their 
(unique) lifts in
$\{X\in \GG\ |\ X(t_1)=0\} \subset \GG \subset (pr_e)_{*}\tm$. 
The projection to $\TT_{M^{(r)}}$ of all possible lifts to $M$ 
of vector fields in $M^{(r)}$ is 
\begin{eqnarray}
d(pr_e): \TT_{M^{(r)}} \oplus (pr_e)_{*}\OO_M \cdot e \to \TT_{M^{(r)}}\ .
\end{eqnarray}

\begin{lemma}
Let $(M,\circ,e,E,t_1)$ be as above.

a)
\begin{eqnarray}
&& (pr_e)_{*}(E\circ \tm) \cap 
\left(\TT_{M^{(r)}}\oplus
 \bigoplus_{k=0}^{n-1}\OO_{M^{(r)}} \cdot t_1^k\cdot e 
\right) \nonumber \\
&=& \bigoplus_{k=1}^{n-1}
 \OO_{M^{(r)}}\cdot \left( t_1^k\cdot e - (t_1e-E)^{\circ k}\right)\ .
\end{eqnarray}

b) Each vector field in $M^{(r)}$ which lifts to a vector field in 
$(pr_e)_{*}\Der_M (\log \DD)$ lifts to a unique vector field in (13.12).

c) The vector fields in $M^{(r)}$ which lift to vector fields in 
$(pr_e)_{*}\Der_M (\log \DD)$ are tangent to the restricted bifurcation diagram
$\BB^{(r)}=pr_e (\BB)\subset M^{(r)}$ and form the free $\OO_{M^{(r)}}$-module
of rank $n-1$
\begin{eqnarray}
\bigoplus_{k=1}^{n-1} \OO_{M^{(r)}} \cdot d(pr_e) ((t_1e-E)^{\circ k})
\subset \Der_{M^{(r)}} (\log \BB^{(r)})\ .
\end{eqnarray}
\end{lemma}

{\bf Proof:}
a) One sees inductively by multiplication with $t_1e-E$ that for any $k\geq 1$
\begin{eqnarray}
t_1^ke-(t_1e-E)^{\circ k} \in (pr_e)_{*}(E\circ \tm) 
\end{eqnarray}
holds. $t_1e-E\in \GG$ and $\lie_e(\circ)=0$ imply $(t_1e-E)^{\circ k}\in \GG$,
therefore
\begin{eqnarray}
\TT_{M^{(r)}}\oplus\bigoplus_{k=0}^{n-1}\OO_{M^{(r)}}t_1^k e = \GG\oplus
\bigoplus_{k=1}^{n-1}\OO_{M^{(r)}}\left(t_1^ke-(t_1e-E)^{\circ k}\right)\ .
\end{eqnarray}
Now the decomposition (13.7) yields (13.12).

b) $pr_e:\DD\to M^{(r)}$ is a branched covering of degree $n$ 
(Definition 11.3), so
\begin{eqnarray}
\bigoplus_{k=0}^{n-1}\OO_{M^{(r)}}\cdot t_1^k \to (pr_e)_{*}(\OO_\DD)
\end{eqnarray}
is an isomorphism. Therefore any lift $h\cdot e +X$, $h\in (pr_e)_{*}\OO_M$,
of $X\in \TT_{M^{(r)}}$ can be replaced by a unique lift 
$\widetilde{h}\cdot e+X$ with 
$\widetilde{h}\in \bigoplus_{k=0}^{n-1}\OO_{M^{(r)}}t_1^k$ and 
$(h-\widetilde{h})|_\DD=0$. 
If $h\cdot e+X$ is tangent to $\DD$, then $(h-\widetilde{h})|_\DD=0$ is 
necessary and sufficient for $\widetilde{h}\cdot e + X$ to be tangent
to $\DD$.

c) A generic point $p^{(r)}\in (\BB^{(r)})_{reg}$ has a preimage 
$p\in (\DD_{sing})_{reg}$ such that the projection of germs 
$pr_e:(\DD_{sing},p)\to (\BB^{(r)},p^{(r)})$ is an isomorphism 
(Corollary 11.5 a)). A vector field $h\cdot e+X$, $X\in \TT_{M^{(r)}}$, which
is tangent to $\DD_{reg}$ is also tangent to $(\DD_{sing})_{reg}$.
Then $X$ is tangent to $(\BB^{(r)})_{reg}$. One obtains the generators
in (13.13) by projection to $\TT_{M^{(r)}}$ of the generators in (13.12).
\hfill $\qed$

\bigskip
$\KK^{(3)}\subset \KK \subset M$ is the set of points $p\in M$ such that
$(M,p)$ does not decompose into 1- and 2-dimensional germs of F-manifolds
(chapter 12).

\begin{theorem}
Let $(M,\circ,e,E)$ be a massive F-manifold with Euler field $E$ of weight
1 and in standard form. Suppose that $\codim \KK^{(3)}\geq 2$.

Then the restricted bifurcation diagram $\BB^{(r)}$ is a free divisor and 
(13.13) is an equality.
\end{theorem}

{\bf Proof:}
In viev of Lemma 13.4 c) it is sufficient to show that any vector field
tangent to $\BB^{(r)}$ lifts to a vector field tangent to $\DD$.

$pr_e:\DD-\BB\to M^{(r)}-\BB^{(r)}$ is a covering of degree $n$. For any vector
field $X\in \TT_{M^{(r)}}$ there exists a unique function 
$h_X\in (pr_e)_{*}\OO_{\DD-\BB}$ such that $\widetilde{h}\cdot e+X$ is tangent
to $\DD-\BB$ iff $\widetilde{h}|_{\DD-\BB}=h_X$.

One has to show that $h_X$ extends to a function in $(pr_e)_{*}\OO_\DD$ if
$X\in \Der_{M^{(r)}}(\log \BB^{(r)})$. Then the unique lift 
$\widetilde{h}\cdot e +X$ with 
$\widetilde{h}\in \bigoplus_{k=0}^{n-1}\OO_{M^{(r)}} t_1^k$
and $\widetilde{h}|_\DD =h_X$ is tangent to $\DD$.

Let $p$ be a point in the set
\begin{eqnarray}
&&\{p\in (\DD_{sing})_{reg}\ |\ p^{(r)}\in (\BB^{(r)})_{reg},\nonumber \\ 
&& \qquad \qquad pr_e:(\DD,p)\to (M^{(r)},p^{(r)}) 
\mbox{ has degree }2\}\ .
\end{eqnarray}
Then the germ $(\DD,p)$ is the product of $(\C^{n-2},0)$ and the 
discriminant of the germ of an F-manifold of type $I_2(m)$ ($m\geq 2$)
(Remark 11.2 iv)). 
One can find coordinates $(t_1,t')=(t_1,...,t_n)$ around $p\in M$ such that
$(\DD,p)\subset (M,p)\to (M^{(r)},p^{(r)})$
corresponds to
\begin{eqnarray}
(\{(t_1,t')\ |\ t_1^2-t_2^m=0\},0) \subset (\C^n,0)\to (\C^{n-1},0),\ 
t\mapsto t'\ .
\end{eqnarray}
Then $(\BB^{(r)},p^{(r)})\cong (\{t'\ |\ t_2=0\},0)$. Obviously the vector
fields tangent to $(\BB^{(r)},p^{(r)})$ have locally lifts to vector fields
tangent to $(\DD,p)$. The function $h_X$ of a field
$X\in \Der_{M^{(r)}}(\log \BB^{(r)})$ extends holomorphically to the set in 
(13.17). The complement in $\DD$ of $\DD_{reg}=\DD-\BB$ and of the 
set in (13.17) has codimension
$\geq 2$ because of $\codim \KK^{(3)}\geq 2$. 
Therefore $h_X\in (pr_e)_{*}\OO_\DD$.
\hfill $\qed$

\section{Isomorphisms and modality of germs of F-manifolds}
\setcounter{equation}{0}

The following three results are applications of Theorem 10.2 for the
Lyashko-Looijenga map. They will be proved
together. $((M,p),\circ,e,\Lambda)$ denotes the germ of an F-manifold
with the function germ $\Lambda:(M,p)\to \C^n$ as additional structure.
A map germ $\varphi:(M,p)\to (M,p)$ respects $\Lambda$ if 
$\Lambda \circ \varphi = \Lambda$.

\begin{theorem}
The automorphism group of a germ $(M,p)$ of a massive F-manifold is finite.
\end{theorem}

\begin{theorem}
Let $(M,\circ, e,F)$ be a massive F-manifold with generating function 
$F$ and Lyashko-Looijenga map $\Lambda:M\to \C^n$. For any 
$p_1\in M$ the set 
\begin{eqnarray*}
\{q\in M\ |\ ((M,p_1),\circ,e,\Lambda)\cong ((M,q),\circ,e,\Lambda)\}
\end{eqnarray*}
is discret and closed in $M$.
\end{theorem}

\begin{corollary}
Let $(M,\circ, e,E)$ be a massive F-manifold with Euler field $E$ of weight
1. For any $p_1\in M$ the set 
\begin{eqnarray*}
\{q\in M\ |\ ((M,p_1),\circ,e,E)\cong ((M,q),\circ,e,E)\}
\end{eqnarray*}
is discret and closed in $M$.
\end{corollary}

{\bf Proof:}
14.3 follows from 14.2. For Theorem 14.1, it suffices to regard an 
irreducible germ of a massive F-manifold. The automorphisms of an 
irreducible germ respect a given Lyashko-Looijenga map $\Lambda$ 
because of Lemma 10.3 a).
So we may fix for 14.1 and 14.2 a massive F-manifold $(M,\circ,e)$ 
and a Lyashko-Looijenga map $\Lambda :M\to \C^n$.

The set 
\begin{eqnarray*}
\Delta := \{ (p,p')\in M\times M\ |\ \Lambda(p)=\Lambda(p')\}
\end{eqnarray*}
has a reduced complex structure. It is a subset of 
$(M-\BB)\times (M-\BB)\cup \BB\times \BB$ and the intersection
$\Delta \cap (M-\BB)\times (M-\BB)$ is smooth of dimension $n$.
This follows from Theorem 10.2

Now consider an isomorphism 
$\varphi: ((M,p),\circ,e,\Lambda) \to ((M,p'),\circ,e,\Lambda)$.
The graph germ 
\begin{eqnarray*}
(\GG (\varphi),(p,p')) := (\{(q,\varphi (q))\in M\times M\ |\ 
q \mbox{ near }p\},
(p,p'))
\end{eqnarray*}
is a smooth analytic germ of dimension $n$ and is contained in the germ
$(\Delta,(p,p'))$. 
It meets
$\Delta \cap (M-\BB)\times (M-\BB)$.
Because of the purity of the dimension of an irreducible analytic germ,
it is an irreducible component of the analytic germ $(\Delta,(p,p'))$.
One can recover the map germ $\varphi$ from the graph germ
$(\GG (\varphi),(p,p'))$. The germ $(\Delta,(p,p'))$ consists
of finitely many irreducible components. The case $p=p'$ together with the
remarks at the beginning of the proof gives Theorem 14.1.

For Theorem 14.2, we assume that there is an infinite sequence
$(p_i,\varphi_i)_{i\in \N}$ of different points $p_i\in M$ and 
map germs
\begin{eqnarray*}
\varphi_i: ((M,p_1),\circ,e,\Lambda) 
\stackrel{\cong}{\to} ((M,p_i),\circ,e,\Lambda)
\end{eqnarray*}
and one accumulation point $p_\infty\in M$. 
The set $\overline{\Delta - \BB \times \BB}$ is analytic of pure 
dimension $n$. It contains the germs 
$(\GG (\varphi_i),(p_1,p_i))$ and the point $(p_1,p_\infty)$.

We can choose a suitable open neighborhood $U$ of $(p_1,p_\infty)$ in 
$M\times M$ and a stratification 
\begin{eqnarray*}
\bigcup_\alpha S_\alpha = U \cap \overline{\Delta - \BB \times \BB}
\end{eqnarray*}
of $U\cap \overline{\Delta - \BB \times \BB}$ which consists of 
finitely many disjoint smooth connected constructible sets $S_\alpha$
and satisfies the boundary condition:
{\it The boundary $\overline{S_\alpha}-S_\alpha$  of a stratum $S_\alpha$
is a union of other strata.}

$(\GG(\varphi_i),(p_1,p_i))$ is an $n$-dimensional irreducible 
component of the $n$-dimensional germ 
$(\overline{\Delta - \BB \times \BB},(p_1,p_i))$. 
There is a unique $n$-dimensional stratum whose closure contains 
$(\GG(\varphi_i),(p_1,p_i))$.
If $(p_1,p_i)\in S_\alpha$ then this together with the boundary 
condition implies 
$(S_\alpha,(p_1,p_i)) \subset (\GG(\varphi_i),(p_1,p_i))$.
The germ $(\GG(\varphi_i),(p_1,p_i))$  is the graph of the isomorphism
$\varphi_i$. Therefore it intersects the germ 
$(\{p_1\}\times M,(p_1,p_i))$ only in $(p_1,p_i)$; the same holds for
$(S_\alpha,(p_1,p_i))$.

Now there exists at least one stratum $S_{\alpha_0}$ which contains 
infinitely many of the points $(p_1,p_i)$. The intersection of the 
analytic sets $\overline{S_{\alpha_0}}$ and $U\cap (\{p_1\}\times M)$
contains these points as isolated points. This is impossible.
The assumption above was wrong.
\hfill $\qed$

\bigskip
In singularity theory there are the notions of $\mu$-constant stratum and
(proper) modality of an isolated hypersurface singularity. One can define
versions of them for the germ $(M,p)$ of an F-manifold $(M,\circ,e)$
(massive or not massive): The {\it $\mu$-constant stratum} $(S_\mu,p)$ is the 
analytic germ of points $q\in M$ such that the eigenspace decompositions
of $T_qM$ and $T_pM$ have the same partition (cf. Proposition 2.3).

The idempotent fields $e_1,...,e_l$ of the decomposition 
$(M,p)=\prod_{k=1}^l (M_k,p)$ into irreducible germs of F-manifolds 
commute and satisfy $\lie_{e_i}(\circ)=0\cdot \circ$.
So the germs $(M,q)$ of points $q$ in one integral manifold of 
$e_1,...,e_l$ are isomorphic. This motivates the definition of the 
{\it modality}:
\begin{eqnarray}
\modmu (M,p):= \dim (S_\mu,p) - l\ .
\end{eqnarray}
Let $(S_\mu^{[k]},p)$ denote the $\mu$-constant stratum of $(M_k,p)$;
Then Theorem 4.2 implies
\begin{eqnarray}
(S_\mu,p) = \prod_k (S_\mu^{[k]},p) \mbox{ \ \ and}\\
\mod (M,p) = \sum_k \modmu (M_k,p)\ .
\end{eqnarray}
For massive F-manifolds, Theorem 14.2 and Lemma 10.3 give more information:

\begin{corollary}
Let $(M,p)=\prod_{k=1}^l (M_k,p)$ be the germ of a massive F-manifold and 
$\Lambda:(M,p)\to \C^n$ a Lyashko-Looijenga map. 

a) There exist a 
representative $S_\mu$ of the $\mu$-constant stratum $(S_\mu,p)$, a
neighborhood $U\subset \C^l$ of $0$ and an isomorphism
\begin{eqnarray}
\psi:S_\mu \to (S_\mu \cap \Lambda^{-1}(\Lambda(p))) \times U
\end{eqnarray}
such that $\psi^{-1}(\{q\}\times U)$ is the integral manifold of 
$e_1,...,e_l$ which contains $q$. 
Any subset of points in $S_\mu\cap \Lambda^{-1}(\Lambda(p))$ 
with isomorphic germs of F-manifolds is discret and closed.

b) 
\begin{eqnarray}
&&\modmu (M,p)=\dim (S_\mu\cap \Lambda^{-1}(\Lambda (p)),p) \ ,\\
&&\sup (\modmu (M,q)\ |\ q \mbox{ near }p)=\dim (\Lambda^{-1}(\Lambda(p)),p)\ .
\end{eqnarray}
\end{corollary}

{\bf Proof:}
a) For $l=1$, the existence of $\psi$ follows from the $e$-invariance of 
$S_\mu$ and from $e(-{1\over n}\Lambda_1)=1$ (Lemma 10.3 b)).

For arbitrary $l$, one uses (14.2) and Lemma 10.3 c): the maps $\Lambda$ and 
$(\Lambda^{[1]},...,\Lambda^{[l]})$ have the same germs of fibers, especially
\begin{eqnarray}
S_\mu\cap \Lambda^{-1}(\Lambda(p)) = 
\prod_k S_\mu^{[k]}\cap {\Lambda^{[k]}}^{-1}(\Lambda^{[k]}(p))\ .
\end{eqnarray}
A germ $(M,q)$ has only a finite number of Lyashko-Looijenga maps with
fixed value at $q$ (Lemma 10.3 a)). 
The finiteness statement in 14.4 a) follows from this and Theorem 14.2.

b) (14.5) follows from a). A representative of the germ 
$(\Lambda^{-1}(\Lambda (p)),p)$ is stratified into constructible subsets which
consist of the points $q$ with the same partition for the eigenspace
decomposition of $T_qM$ (Proposition 2.3).
A point $q\in \Lambda^{-1}(\Lambda (p))$ in an open stratum with maximal 
dimension satisfies
\begin{eqnarray}
\modmu (M,q)&=& \dim (S_\mu(q)\cap \Lambda^{-1}(\Lambda(p)),q) \nonumber \\
&=& \dim (\Lambda^{-1}(\Lambda(p)),q)
=\dim (\Lambda^{-1}(\Lambda(p)),p)\ .
\end{eqnarray}
This shows
\begin{eqnarray}
\sup (\modmu (M,q)\ |\ q\in \Lambda^{-1}(\Lambda(p)) \mbox{ near }p)
=\dim (\Lambda^{-1}(\Lambda (p)),p)\ .
\end{eqnarray}
The upper semicontinuity of the fiber dimension of $\Lambda$ gives (14.6).
\hfill $\qed$

\begin{remark}
Gabrielov \cite{Ga}
proved in the case of isolated hypersurface singularities the upper
semicontinuity of the modality,
\begin{eqnarray}
\modmu (M,q)\leq \modmu (M,p) \mbox{ \ for }q\mbox{ near }p
\end{eqnarray}
(and the equality with another version of modality which was defined by
Arnold). 
He used (14.5), (14.6), and a result of himself, Lazzeri, and L\^e, which,
translated to the F-manifold of a singularity (chapter 16), says:
\begin{eqnarray}
(S_\mu\cap \Lambda^{-1}(\Lambda(p)),p) = 
(\Lambda^{-1}(\Lambda(p)),p)\ .
\end{eqnarray}
(14.10) is an immediate consequence of (14.5), (14.6) and (14.11).
But for other F-manifolds (14.11) and (14.10) are not clear.
\end{remark}

\medskip
In the case of the simple hypersurface singularities, the base of the 
semiuniversal unfolding is an F-manifold $M\cong \C^n$ and the map 
$\Lambda: M- \BB \to \C^n- \DD^{(n)}$ is a finite covering.
Therefore the complement $M-\BB$ is a $K(\pi,1)$ space and the 
fundamental group is a subgroup of finite index of the braid group
$Br(n)$. This is the application of Looijenga \cite{Lo1} and Lyashko
\cite{Ar1} of the map $\Lambda$, which led to the name 
Lyashko-Looijenga map.

It can be generalized to F-manifolds. We call a massive F-manifold $M$
{\it simple} if $\modmu (M,p)=0$ for all $p\in M$.
This fits with the notions of simple hypersurface singularities, simple
boundary singularities, and simple Lagrange maps 
(\cite{Gi2} 1.3, p 3251).

A distinguished class of simple F-manifolds are the F-manifolds of the 
finite Coxeter groups (chapter 18 and 
\cite{Lo1}\cite{Ar1}\cite{Ly1}\cite{Ly3}\cite{Gi2}).
There are other examples (Proposition 20.6 and Remark 20.7). 

A Lyashko-Looijenga map of a massive F-manifold is locally a branched
covering iff $M$ is simple ((14.6) and Theorem 10.2). A detailed proof
of the following result had been given by Looijenga \cite{Lo1} (Theorem 2.1)
(cf. also \cite{Gi2} 1.4, Theorem 5).

\begin{theorem}
Let $(M,p)=(\C^n,p)$ be the germ of a simple F-manifold with fixed 
coordinates.
Then, if $\varepsilon < \varepsilon_0$ for some $\varepsilon_0$, the space
$\{z\in \C^n\ |\ |z|<\varepsilon\}- \BB$ is a $K(\pi,1)$ space. 
Its fundamental group is a subgroup of finite index of the braid group 
$Br(n)$.
\end{theorem}

\section{Analytic spectrum embedded differently}
\setcounter{equation}{0}

The analytic spectrum $L\subset T^{*}M$ of an F-manifold determines the
multiplication on $\tm$ via the isomorphism ((6.1) and (2.2))
\begin{eqnarray}
\aaa:\tm \to \pi_{*}\OO_L,\ X\mapsto \alpha(\widetilde{X})|_L\ .
\end{eqnarray}
One can generalize this and replace $L$, $T^{*}M$, and $\alpha$ by other
spaces and other 1--forms. This allows to find F-manifolds in natural
geometric situations and to encode them in an economic way.
Corollary 15.2 and Definition 15.4 are the two most interesting
special cases of Theorem 15.1.

\begin{theorem}
Let the following data be given:\\
manifolds $Z$ and $M$, where $M$ is connected and $n$-dimensional;\\
a surjective map $\pi_Z:Z\to M$ which is everywhere a submersion;\\
an everywhere $n$-dimensional reduced subvariety $C\subset Z$ 
such that the restriction $\pi_C:C\to M$ is finite;\\
a 1--form $\alpha_Z$ on $Z$ with the property: 
\begin{eqnarray}
&& \mbox{any local lift }
\widetilde{X}\in \TT_Z \mbox{  of the zero vector field }0\in \tm\nonumber\\
&& \mbox{ satisfies }\alpha_Z(\widetilde{X})|_C=0\ .
\end{eqnarray}
Then

a) The map 
\begin{eqnarray}
\aaa_C:\tm\to (\pi_C)_{*}\OO_C, \ X\mapsto \alpha_Z(\widetilde{X})|_C
\end{eqnarray}
is welldefined; here $\widetilde{X}\in \TT_Z$ is any lift of $X$ to a 
neighborhood of $C$ in $Z$.

\smallskip
b) The image $L\subset T^{*}M$ of the map
\begin{eqnarray}
\qqq:C &\to & T^{*}M,\\ 
z &\mapsto & \qqq(z)=\bigl(X\mapsto\aaa_C(X)(z)\bigr)
\in T^{*}_{\pi_C(z)}M    \nonumber
\end{eqnarray}
is a (reduced) variety. $\qqq:C\to L$ is a finite map, the projections 
$\pi:L\to M$ and $\pi_C=\pi\circ \qqq$ are branched coverings. 
The composition of the maps 
${\bf \hat{q}}:\pi_{*}\OO_L \to (\pi_C)_{*}\OO_C$ and 
\begin{eqnarray}
\aaa :\tm\to \pi_{*}\OO_L, \ X\mapsto \alpha(\widetilde{X})|_L
\end{eqnarray}
is $\aaa_C={\bf \hat{q}}\circ \aaa$. All three are $\OO_M$-module homomorphisms.

\smallskip
c) The 1--forms $\alpha$ and $\alpha_Z$ satisfy 
$(\qqq^{*}\alpha)|_{C_{reg}} = \alpha_Z|_{C_{reg}}$. 
Therefore $L$ is a Lagrange
variety iff $\alpha_Z|_{C_{reg}}$ is exact.

\smallskip
d) $\aaa:\tm\to\pi_{*}\OO_L$ is an isomorphism iff
\begin{list}{}{}
\item[i)] 
$\aaa_C$ is injective,
\item[ii)]
its image $\aaa_C(\tm )\subset (\pi_C)_{*}\OO_C$ is multiplication
invariant,
\item[iii)]
$\aaa_C(\tm)$ contains the unit $1_C\in (\pi_C)_{*}\OO_C$.
\end{list}
In this case $\aaa_C:\tm \to (\pi_C)_{*}\OO_C$ induces a (commutative
and associative and) generically semisimple multiplication on $\tm$
with global unit field and with analytic spectrum $L$.

\smallskip
e) $\aaa_C:\tm \to (\pi_C)_{*}\OO_C$ provides $M$ with the structure of a 
massive F-manifold iff $\alpha_Z|_{C_{reg}}$ is exact and the conditions
i)--iii) in d) are satisfied.
\end{theorem}

{\bf Proof:}
a) This follows from (15.2).

b) $\dim C=n=\dim M$ and $\pi_C$ finite imply that $\pi_C$ is open.
$M$ is connected, thus $\pi_C$ is a branched covering. Using local 
coordinates for $M$ and $T^{*}M$ one sees that $\qqq:C \to T^{*}M$ is an
analytic map. $\pi_C=\pi\circ \qqq$ is clear and shows that $\qqq$ is finite.
Then $L=\qqq(C)$ is a variety and $\pi$ is a branched covering.
$\aaa_C={\bf \hat{q}}\circ \aaa$ follows from the definition of $\qqq$.

c) There is an open subset $M^{(0)}\subset M$ with analytic complement
$M-M^{(0)}$ such that $\pi_C^{-1}(M^{(0)})\subset C$ and 
$\pi^{-1}(M^{(0)}) \subset L$ are smooth, 
$\pi_C:\pi_C^{-1}(M^{(0)})\to M^{(0)}$ and
$\pi:\pi^{-1}(M^{(0)})\to M^{(0)}$ are coverings and 
$\qqq:\pi_C^{-1}(M^{(0)})\to \pi^{-1}(M^{(0)})$ is a covering on each component
of $\pi^{-1}(M^{(0)})$.
Now $\aaa_C={\bf \hat{q}}\circ \aaa$ implies 
$\qqq^{*}\alpha|\pi_C^{-1}(M^{(0)})
=\alpha_Z|\pi_C^{-1}(M^{(0)})$.

d) ${\bf \hat{q}}:\pi_{*}\OO_L\to (\pi_C)_{*}\OO_C$ is an injective homomorphism
of $\OO_M$-algebras.
If $\aaa:\tm \to \pi_{*}\OO_L$ is an isomorphism then i)--iii) are obviously
satisfied.

Suppose that i)--iii) are satisfied. Then $\aaa:\tm \to \pi_{*}\OO_L$
is injective with multiplication invariant image 
$\aaa(\tm)\subset \pi_{*}\OO_L$ and with
$1_L\in \aaa(\tm)$. The maps $\aaa$ and $\aaa_C$ induce the same (commutative
and associative) multiplication with global unit field on $TM$.

We have to show that this multiplication is generically semisimple with
analytic spectrum $L$. Then $\aaa:\tm \to \pi_{*}\OO_L$ is an isomorphism and 
the proof of d) is complete.

If for each $p\in M$ the linear forms in $\pi^{-1}(p)\subset T^{*}_pM$ would
generate a subspace of $T^{*}_pM$ of dimension $<n$ then $\aaa$ would not be
injective. So, for a generic point $p\in M$ there exist $n$ elements in 
$\pi^{-1}(p)\subset T^{*}_pM$ which form a basis of $T^{*}_pM$. 
We claim that $\pi^{-1}(p)$ contains no other than these elements:
$\pi^{-1}(p)$ does not contain $0\in T^{*}_pM$ because of $1_L\in \aaa(\tm)$.
From the multiplication invariance of $\aaa(\tm)$ one derives easily that 
$\pi^{-1}(p)$ does not contain any further elements.

This extends to a small neighborhood $U$ of the generic point $p\in M$:
$\pi^{-1}(U)$ consists of $n$ sheets which form a basis of sections of
$T^{*}M$; the map $\aaa|_U:\TT_U \to \pi_{*}(\pi^{-1}(U))$ is an isomorphism
and induces a semisimple multiplication on $TM$ with analytic spectrum 
$\pi^{-1}(U)$.

Then $L$ is the analytic spectrum of the multiplication on $\tm$ because
$M$ is connected.

e) c)+d) and Theorem 6.2.
\hfill $\qed$

\bigskip

In Theorem 15.1 the map $\pi:L\to M$ has degree $n$, but $\pi_C:C\to M$
can have degree $>n$; and even if $\pi_C:C\to M$ has degree $n$ the  map
$\qqq:C\to L$ does not need to be an isomorphism. Examples will be 
discussed below (Examples 15.5, Lemma 17.9). But the most important
special case is the following.

\begin{corollary}
Let $Z, M, \pi_Z, C\subset Z, \alpha_Z,\aaa_C, L$, and $\qqq$ be as in 
Theorem 15.1. Suppose that $\alpha_Z|_{C_{reg}}$ is exact and 
$\aaa_C:\tm \to (\pi_C)_{*}\OO_C$ is an isomorphism. 

Then $\qqq:C\to L$ is an isomorphism and $\aaa_C={\bf \hat{q}}\circ \aaa$ provides 
$M$ with the structure of a massive F-manifold with analytic spectrum $L$.
\end{corollary}

{\bf Proof:}
Theorem 15.1 e) gives all except for the isomorphism $\qqq:C\to L$.
This follows from the isomorphism 
${\bf \hat{q}}:\pi_{*}\OO_L\to (\pi_C)_{*}\OO_C$ and a universal property of the 
analytic spectrum.
\hfill $\qed$

\bigskip

One can encode an irreducible germ of a massive F-manifold with data
as in Corollary 15.2 such that the dimension of $Z$ is minimal.

\begin{lemma}
Let $(M,p)$ be an irreducible germ of a massive $n$-dimensional F-manifold.
$\mmm \subset T_pM$ denotes the maximal ideal in $T_pM$.

a) $\dim Z \geq n+\dim \mmm/\mmm^2$ for any data as in Corollary 15.2 for
$(M,p)$.

b) There exist data as in Corollary 15.2 for $(M,p)$ with
$\dim Z =n+\dim \mmm/\mmm^2$ (the construction will be given in the proof).
\end{lemma}

{\bf Proof:}
a) $\pi^{-1}_C(p)=\pi_Z^{-1}(p)\cap C$ consists of one fat point with 
structure ring $T_pM$. Its embedding dimension $\dim \mmm/\mmm^2$ is bounded 
by the dimension $\dim \pi_Z^{-1}(p)=\dim Z-n$ of the smooth fiber 
$\pi_Z^{-1}(p)$.

b) One can choose coordinates $(t_1,...,t_n)=(t_1,t')=t$ for $(M,p)$ with
$e={\partial \over \partial t_1}$ as usual and with 
\begin{eqnarray}
&& \bigoplus_{i=2}^n\C\cdot [{\partial\over \partial t_i}] =\mmm\subset T_pM
\mbox{ \ and }\\
&& \bigoplus_{i=m+1}^n\C\cdot [{\partial\over \partial t_i}] =\mmm^2\subset T_pM
\end{eqnarray}
for $m=1+\dim \mmm/\mmm^2$. The dual coordinates on $(T^{*}M,T^{*}_pM)$ are
$y_1,...,y_n$, the analytic spectrum is (cf. (2.1))
\begin{eqnarray}
L=\{(y,t)\ |\ y_1=1,\ y_iy_j=\sum a_{ij}^k(t')y_k\}\ .
\end{eqnarray}
Because of (15.7) there exist functions $b_i\in \C\{t'\}[y_2,...,y_m]$ with
\begin{eqnarray}
y_i|_L = b_i(y_2,...,y_m,t')|_L \mbox{ \ \ \ for } i=m+1,...,n\ .
\end{eqnarray}
We identify $(M,p)$ and $(\C^n,0)$ using $(t_1,...,t_n)$ and define
$(Z,0)=(\C^{m-1}\times \C^n,0)$. The embedding
\begin{eqnarray}
\iota:(Z,0)=(\C^{m-1}\times\C^n,0) \hookrightarrow  T^{*}M\ ,\qquad \qquad &&\\
(x_1,...,x_{m-1},t) \mapsto  (y,t)=(1,x_1,...,x_{m-1},b_{m+1}(x,t'),...,
b_n(x,t'),t)&&\nonumber
\end{eqnarray}
provides canonical choices for the other data,
\begin{eqnarray}
&& \pi_Z:(Z,0)\to (M,p), (x,t)\mapsto t\ ,\\
&& C=\iota^{-1}(L)\ ,\\
&&\alpha_Z=\iota^{*}\alpha = dt_1+\sum_{i=2}^mx_{i-1}dt_i
+\sum_{i=m+1}^nb_i(x,t')dt_i\ .
\end{eqnarray}
The conditions in Corollary 15.2 are obviously satisfied.
\hfill $\qed$

\bigskip

The notion of a generating family for a Lagrange map (\cite{AGV} 19., 
\cite{Gi2} 1.4) motivates to single out another special case of Theorem 15.1.

\begin{definition}
Let $Z, M, \pi_Z, C, \alpha_Z$, and $\aaa_C$ be as in Theorem 15.1 with 
$\alpha_Z|_{C_{reg}}$ exact and $\aaa_C:\tm \to (\pi_C)_{*}\OO_C$ injective
with multiplication invariant image $\aaa_C(\tm)\supset \{ 1_C\}$.
These data yield a massive F-manifold $(M,\circ,e)$.

A function $F:Z\to \C$ is a {\it generating family} for this F-manifold 
if $\alpha_Z=dF$ and if $C$ is the critical set of the map 
$(F,\pi_Z):Z\to \C\times M$.
\end{definition}

\smallskip
The name {\it generating family} has two reasons: \\
1) $F$ is considered as a family of functions 
on the fibers $\pi^{-1}_Z(p)$, $p\in M$.\\
2) The restriction of $F$ to $C$ is the lift of a generating function
$\widetilde{F}:L\to \C$, i.e. $F=\widetilde{F}\circ \qqq$;
so the 1-graph of $F$ as a multivalued function on $M$ is $L$.

In the case of a generating family the conditions (15.2) and $\alpha_Z$ exact
are obvious. The most difficult condition is the multiplication invariance
of $\aaa (\tm)$. 

It is not clear whether for any massive F-manifold $M$ data $(Z,\pi_Z,F)$ as 
in Definition 15.4 exist. But often even many nonisomorphic data exist.
We illustrate this for the 2-dimensional germs $I_2(m)$ of F-manifolds
(chapter 11).

\begin{examples}
Always $(Z,0)=(\C\times \C^2,0)$ and $(M,p)=(\C^2,0)$ with projection 
$\pi_Z:(Z,0)\to (M,0),\ (x,t_1,t_2)\mapsto(t_1,t_2)$ and 
$e:=\delta_1:= {\partial\over \partial t_1}, 
\delta_2:={\partial \over \partial t_2}$.

\smallskip
a) $C=\{(x,t)\ |\ x^{m-2}-t_2^2=0\},\ \alpha_Z=dt_1+xdt_2\ .$\\
These are data as in 
Corollary 15.2 for $I_2(m)$.

\smallskip
b) Generating family $F=t_1+\int_0^x(t_2-u^2)^kdu$ ($k\geq 1$), \\
$C=\{(x,t)\ |\ t_2-x^2=0\}$, \\
$\alpha_Z|_C=dF|_C=(dt_1+c\cdot x^{2k-1}dt_2)|_C$ for some $c\in \C-\{0\}$,\\
$\aaa_C(\delta_2)\cdot \aaa_C(\delta_2)=c^2\cdot t_2^{2k-1}\cdot 1_C\ .$\\
These are 
data as in Definition 15.4 for $I_2(2k+1)$, $\pi_C:C\to M$ has degree 2,
$\qqq:C\to L$ is the normalisation and the maximalisation of $L$.

\smallskip
c) Generating family $F=t_1+x^{k+1}t_2-{k+1\over k+2}x^{k+2}$ ($k\geq 1$),\\
$C=\{(x,t)\ |\ (t_2-x)x=0\}$, \\
$\alpha_Z|_{C_{reg}}=dF|_{C_{reg}}=(dt_1+x^{k+1}dt_2)|_{C_{reg}}$,\\
$\aaa_C(\delta_2-{1\over 2}t_2^{k+1}\cdot \delta_1)^2 
= {1\over 4}t_2^{2k+2}\cdot 1_C\ .$\\
These are
data as in Definition 15.4 for $I_2(2k+4)$, $\pi_C:C\to M$ has degree 2, 
$\qqq:C\to L$ is the maximalisation of $L$
(for the missing case $I_2(4)$ compare Lemma 17.9).

\smallskip
d) Generating family $F=t_1+\int_0^x (u^2-t_2)^kudu $ ($k\geq 1$),\\
$C=\{(x,t)\ |\ (t_2-x^2)x=0\}$, \\
$\alpha_Z|_{C_{reg}} = dF|_{C_{reg}} = (dt_1+c\cdot x^{2k}dt_2)|_{C_{reg}}$
for some $c\in \C-\{0\}$, \\
$\aaa_C(\delta_2-{1\over 2}ct_2^k\cdot \delta_1)^2 = 
{1\over 4}c^2t_2^{2k}\cdot 1_C\ .$ \\
These are
data as in Definition 15.4 for $I_2(2k+2)$, 
$\pi_C:C\to M$ has degree 3, $\qqq:C\to L$ covers one component with degree
1, the other with degree 2.
\end{examples}

\section{Hypersurface singularities}
\setcounter{equation}{0}

A distinguished class of germs of massive F-manifolds is related to isolated
hypersurface singularities: the base space of a semiuniversal unfolding
of an isolated hypersurface singularity is an irreducible germ of a 
massive F-manifold with smooth analytic spectrum (Theorem 16.3). In fact,
there is a 1-1 correspondence between such germs of F-manifolds and
singularities up to stable right equivalence (Theorem 16.6).

The structure of an F-manifold on the base space has beautiful geometric
implications and interpretations (Theorem 16.4, Remarks 16.5). 
Many of them have been known since long time from different points of view.
The concept of an F-manifold unifies them. On the other hand,
for much of the general treatment of F-manifolds in this paper the
singularity case had been the model case.

\bigskip
An isolated hypersurface singularity is a holomorphic function germ
$f:(\C^m,0)\to (\C,0)$ with an isolated singularity at 0. 
Its Milnor number $\mu\in \N$ is the dimension of the Jacobi algebra
$\OO_{\C^m,0}/({\partial f\over \partial x_1},...,
{\partial f\over \partial x_m})
=\OO_{\C^m,0}/J_f$.

The notion of an unfolding of an isolated hypersurface singularity goes
back to Thom and Mather.
An {\it unfolding} of $f$ is a holomorphic function germ 
$F:(\C^m\times \C^n,0)\to (\C,0)$ such that $F|_{\C^m\times \{0\}}=f$.
We will write the parameter space as $(M,0)=(\C^n,0)$.

The {\it critical space} $(C,0)\subset (\C^m\times M,0)$ of the unfolding 
$F=F(x_1,...,x_m,t_1,...,t_n)$ is the critical space of the map
$(F,pr):(\C^m\times M,0)\to (\C\times M,0)$. 
It is the zero set of the ideal
\begin{eqnarray}
J_F := ({\partial F\over \partial x_1},...,
{\partial F\over \partial x_m})
\end{eqnarray}
with the complex structure $\OO_{C,0}=\OO_{\C^m\times M,0}/J_F|_{(C,0)}$.

The intersection $C\cap (\C^m\times \{0\})=\{0\}$ is a point and 
$(C,0)$ is a complete intersection of dimension $n$. Therefore the projection
$pr:(C,0)\to (M,0)$ is finite and flat with degree $\mu$
and $\OO_{C,0}$ is a free $\OO_{M,0}$-module of rank $\mu$.

A kind of {\it Kodaira-Spencer map} is the $\OO_{M,0}$-linear map
\begin{eqnarray}
\aaa_C:\TT_{M,0}\to \OO_{C,0},\ X\mapsto \widetilde{X}(F)|_{(C,0)}
\end{eqnarray}
where $\widetilde{X}$ is any lift of $X\in \TT_{M,0}$ to $(\C^m\times M,0)$.
Dividing out the submodules $\mmm_{M,0}\cdot \TT_{M,0}$ and 
$\mmm_{M,0}\cdot \OO_{C,0}$ one obtains the {\it reduced Kodaira-Spencer map}
\begin{eqnarray}
\aaa_C|_0:T_0M \to \OO_{\C^m,0}/J_f\ .
\end{eqnarray}
All these objects are independent of the choice of coordinates. In fact, 
they even behave well with respect to morphisms of unfoldings.

\medskip
There are several possibilities to define morphisms of unfoldings
(cf. Remark 16.2 iv)). We need the following.

Let $F_i:(\C^m\times M_i,0)\to (\C,0),\ i=1,2,$ be two unfoldings of $f$
with projections $pr_i:(\C^m\times M_i,0)\to (M_i,0)$, critical spaces
$C_i$, and Kodaira-Spencer maps $\aaa_{C_i}$. A {\it morphism}
from $F_1$ to $F_2$ is a pair $(\phi, \phi_{base})$ of map germs
such that the following diagram commutes,
\begin{eqnarray}
\begin{CD} 
(\C^m\times M_1,0) @>\phi>> (\C^m\times M_2,0)\\
@VV{pr_1}V                 @VV{pr_2}V       \\  
(M_1,0)      @>{\phi_{base}}>> (M_2,0)\ ,
\end{CD}
\end{eqnarray}
and  
\begin{eqnarray}
\phi|_{\C^m\times \{0\}}=id\ ,\\
F_1=F_2\circ \phi
\end{eqnarray}
hold. One says that $F_1$ is {\it induced} by $(\phi,\phi_{base})$ from $F_2$.

The definition of critical spaces is compatible with the morphism
$(\phi,\phi_{base})$, that is, $\phi^{*}J_{F_2}=J_{F_1}$ and 
$(C_1,0)=\phi^{-1}((C_2,0))$. 
Also the Kodaira-Spencer maps behave well: the $\OO_{M_1,0}$-linear maps
\begin{eqnarray}
d\phi_{base} &:&  \TT_{M_1,0}\to 
\OO_{M_1,0}\otimes _{\OO_{M_2,0}} \TT_{M_2,0}\ ,\\
\aaa_{C_2}    &:&  \OO_{M_1,0}\otimes_{\OO_{M_2,0}} \TT_{M_2,0} \to
\OO_{M_1,0}\otimes_{\OO_{M_2,0}} \OO_{C_2,0}\ ,\\
\phi^{*}|_{(C_2,0)} &:& \OO_{M_1,0}\otimes_{\OO_{M_2,0}} \OO_{C_2,0} \to 
\OO_{C_1,0}
\end{eqnarray}
are defined in the obvious way; their composition is 
\begin{eqnarray}
\aaa_{C_1}=\phi^{*}|_{(C_2,0)}\circ \aaa_{C_2} \circ d\phi_{base}\ .
\end{eqnarray}
(16.9) restricts to the identity on the Jacobi algebra of $f$ because of
(16.5). Therefore the reduced Kodaira-Spencer maps satisfy
\begin{eqnarray}
\aaa_{C_1}|_0 = \aaa_{C_2}|_0 \circ d\phi_{base}|_0\ .
\end{eqnarray}
An unfolding of $f$ is {\it versal} if any unfolding is induced from it
by a suitable morphism. A versal unfolding 
$F:(\C^m\times M,0)\to (\C,0)$ is {\it semiuniversal} if the dimension of
the parameter space $(M,0)$ is minimal. Semiuniversal unfoldings 
of an isolated hypersurface singularity exist by work of Thom and Mather.
Detailed proofs can nowadays be found at many places, e.g. 
\cite{Was}\cite{AGV}(8.).

\begin{theorem}
An unfolding $F:(\C^m\times M,0)\to (\C,0)$ of an isolated hypersurface
singularity $f:(\C^m,0)\to (\C,0)$ is versal iff the reduced Kodaira-Spencer
map $\aaa_C|_0:T_0M\to \OO_{\C^m,0}/J_f$ is surjective. 
It is semiuniversal iff $\aaa_C|_0$ is an isomorphism.
\end{theorem}

\begin{remarks}
i) Because of the lemma of Nakayama $\aaa_C|_0$ is surjective (an isomorphism)
iff $\aaa_C$ is surjective (an isomorphism).

\smallskip
ii) A convenient choice of a semiuniversal unfolding 
$F:(\C^m\times \C^\mu,0)\to (\C,0)$ is 
$F(x_1,...,x_m,t_1,...,t_\mu)=f+\sum_{i=1}^\mu m_it_i$, where 
$m_1,...,m_\mu \in \OO_{\C^m,0}$ represent a basis of the Jacobi algebra
of $f$, preferably with $m_1=1$.

\smallskip
iii) The critical space of an unfolding $F:(\C^m\times \C^n,0)\to(\C,0)$ is 
reduced and smooth iff the matrix 
$\bigl( {\partial^2 F\over \partial x_i \partial x_j},
{\partial^2 F\over \partial x_i \partial t_k}\bigr) (0)$ has maximal rank $m$.
This is satisfied for versal unfoldings.

\smallskip
iv) In the literature (e.g. \cite{Was}) one often finds a slightly different
notion of morphisms of unfoldings: An {\it $(r)$-morphism} between
unfoldings $F_1$ and $F_2$ as above is a triple $(\phi,\phi_{base},\tau)$
of map germs $\phi$ and $\phi_{base}$ with (16.4) and (16.5)  and 
$\tau:(M_1,0)\to (\C,0)$ with (16.6) replaced by
\begin{eqnarray}
F_1=F_2\circ \phi +\tau\ .
\end{eqnarray}
{\it $(r)$-versal} and {\it $(r)$-semiuniversal} unfoldings are defined 
analogously. They exist (\cite{Was}): an unfolding
$F:(\C^m\times M,0)\to (\C,0)$ is $(r)$-versal ($(r)$-semiuniversal)
iff the map 
\begin{eqnarray}
\C\oplus T_0M \to \OO_{\C^m,0}/J_f,\ 
(c,X)\mapsto c+\aaa_C|_0(X)
\end{eqnarray}
is surjective (an isomorphism).

So one gains a bit: the base space of an $(r)$-semiuniversal unfolding 
$F^{(r)}$ has dimension $\mu-1$; if 
$F^{(r)}=F^{(r)}(x_1,...,x_m,t_2,...,t_\mu)$ is ${(r)}$-semiuniversal
then $t_1+F^{(r)}$ is semiuniversal; between two semiunversal unfoldings
$t_1+F^{(r)}_1$ and $t_1+F^{(r)}_2$ of this form there exist isomorphisms which
come from $(r)$-isomorphisms of $F^{(r)}_1$ and $F^{(r)}_2$.
(The relation between $F^{(r)}$ and $t_1+F^{(r)}$ motivates the ``$(r)$'',
which stands for ``restricted'').

On the other hand, one looses (16.10). Anyway, one should keep 
$(r)$-semiuniversal unfoldings in mind. They are closely related to 
miniversal Lagrange maps (see the proof of Theorem 16.6 and \cite{AGV} 19.).

\smallskip
v) One can generalize the notion of a morphism between unfoldings if one
weakens the condition (16.5): 
Let $F_i:(\C^m\times M_i),0)\to (\C,0),\ i=1,2,$ be unfoldings of two 
isolated hypersurface singularities $f_1$ and $f_2$. 
A {\it generalized morphism} from $F_1$ to $F_2$ is a pair $(\phi,\phi_{base})$
of map germs with a commutative diagram as in (16.4) such that (16.6) holds
and $\phi|_{\C^m\times\{0\}}$ is a coordinate change (between $f_1$ and
$f_2$). 

Then $f_1$ and $f_2$ are {\it right equivalent}. If the generalized
morphism is invertible then also $F_1$ and $F_2$ are called 
{\it right equivalent}.

Critical spaces and Kodaira-Spencer maps behave well also for 
generalized morphisms; (16.10) holds, in (16.11) one has to take into
account the isomorphism of the Jacobi algebras of $f_1$ and $f_2$ which
is induced by $\phi|_{\C^m\times\{0\}}$.
\end{remarks}

The multiplication on the base space of a semiuniversal unfolding was first
defined by K. Saito \cite{SaK4}(1.5)\cite{SaK5}(1.3).

\begin{theorem}
Let $f:(\C^m,0)\to (\C,0)$ be an isolated hypersurface singularity and
$F:(\C^m\times M,0)\to (\C,0)$ a semiuniversal unfolding.

The Kodaira-Spencer map $\aaa_C:\TT_{M,0}\to \OO_{\C,0}$ is an isomorphism
and induces a multiplication on $\TT_{M,0}$. Then $(M,0)$ is an irreducible
germ of a massive F-manifold with smooth analytic spectrum.
$E:=\aaa_C^{-1}(F|_{C})$ is an Euler field of weight 1.
\end{theorem}

{\bf Proof:}
$\aaa_C:\TT_{M,0}\to \OO_{C,0}$ is an isomorphism because of Theorem 16.1 and
Remark 16.2 i). The critical space $(C,0)$ is reduced and smooth.
One applies Corollary 15.2 to $(Z,0)=(\C^m\times M,0)$ and 
$\alpha_Z=dF$. 
The map $\qqq:(C,0)\to (L,\pi^{-1}(0))$ is an isomorphism. 
$\pi^{-1}(0)$ is a point. Theorem 15.1 c) shows that 
$F|_{C}\circ \qqq^{-1}$ is a holomorphic generating function.
Therefore $E$ is an Euler field of weight 1.
\hfill $\qed$

\begin{theorem}
Let $f:(\C^m,0)\to (\C,0)$ be an isolated hypersurface singularity and 
$F_i:(\C^m\times M_i,0)\to (\C,0)$, $i=1,2,$ be two semiuniversal
unfoldings. 

There exists a unique isomorphism $\varphi:(M_1,0)\to (M_2,0)$ of 
F-manifolds such that $\phi_{base}=\varphi$  for any isomorphism
$(\phi,\phi_{base})$ of the unfoldings $F_1$ and $F_2$.
\end{theorem}

{\bf Proof:}
$\phi_{base}:(M_1,0)\to (M_2,0)$ is an isomorphism of F-manifolds 
because of (16.10).

Suppose that $F_1=F_2$ and $(M_1,0)=(M_2,0)$. The tangent map of 
$\phi_{base}$ on $T_0M_1$ is $d\phi_{base}|_0=id$ because of (16.11).
The group of all automophisms of $(M_1,0)$ as F-manifold is finite
(Theorem 14.1). Therefore $\phi_{base}=id$.
\hfill $\qed$

\begin{remarks}
i) The rigidity of the base morphism $\phi_{base}$ in Theorem 16.4 is 
in sharp contrast to the general situation for deformations of 
geometric objects. Usually only a part of the base space of a 
miniversal deformation is rigid with respect to automorphisms
of the deformation.

\smallskip
ii) The reason for the rigidity are, via Theorem 14.1 and Theorem 10.2, 
the canonical coordinates at generic parameters. The corresponding result
for singularities is that the critical values of $F$ form coordinates
on the base at generic parameters. It had been proved by 
Looijenga \cite{Lo1}.

\smallskip
iii) Because of this rigidity also the openness of versality 
(e.g. \cite{Tei}) takes a special form:
For any point $t\in M$ in a representative of the base space 
$(M,0)=(\C^\mu,0)$ of a semiuniversal unfolding $F$, Theorem 4.2 yields
a unique decomposition $(M,t)=\prod_{k=1}^l (M_k,t)$ into a product of 
irreducible germs of F-manifolds. These germs $(M_k,t)$ are the base spaces
of semiuniversal unfoldings of the singularities of $F|_{\C^m\times \{t\}}$.
The multigerm of $F$ at $\C^m\times \{t\}\cap C$ itself is isomorphic
-- in a way which can be made precise easily-- to a transversal union of 
versal unfoldings of these singularities.

\smallskip
iv) The tangent space $T_tM\cong \bigoplus_{k=1}^l T_tM_k$ is canonically 
isomorphic to the direct sum of the Jacobi algebras of singularities of
$F|_{\C^m\times\{t\}}$. The vector in $T_tM$ of the Euler field $E$
is mapped to the direct sum of the classes of the function
$F|_{\C^m\times\{t\}}$ in these Jacobi algebras. 
A result of Scherk \cite{Sch} says:
\begin{list}{}{}
\item[] {\it The Jacobi algebra $\OO_{\C^m,0}/J_f$ of an isolated hypersurface
singularity 
$f:(\C^m,0)\to (\C,0)$ together with the class $[f]\in \OO_{\C^m,0}/J_f$
determines $f$ up to right equivalence.}
\end{list}
This result shows that the base space $M$ as an F-manifold with Euler field
$E$ determines for each parameter $t\in M$ the singularities
of $F|_{\C^m\times\{t\}}$ up to right equivalence and also
the critical values.
Theorem 16.6 will give an even stronger result.

\smallskip
v) The eigenvalues of $E\circ:T_tM\to T_tM$ are the critical values of 
$F|_{\C^m\times\{t\}}$. Therefore the discriminant of the Euler field
$E$ is
\begin{eqnarray}
\DD = \{t\in M\ |\ (\det (E\circ ))(t)=0\} = \pi_C(C\cap F^{-1}(0))
\end{eqnarray}
and coincides with the classical discriminant of the unfolding $F$.

All the results of chapter 11 apply to this discriminant.
Of course, many of them are classic in the singularity case.

For example, Theorem 11.1 and the isomorphism $\qqq:C\to L$ from
Corollary 15.2 yield an isomorphism between the development
$\widetilde{\DD}\subset \P T^{*}M$ of the discriminant and the smooth
variety $C\cap F^{-1}(0)$ which had been established by Teissier
\cite{Tei}. Implicitly it is also in \cite{AGV}(19.).

The elementary way in Corollary 11.6 how discriminant and unit field
determine the Jacobi algebras seems to be new. But the consequence from
this and Scherk's result that discriminant and unit field determine
the singularity (up to right equivalence) is known (compare below
Theorem 16.6 and Remark 16.7 iv)).
\end{remarks}

\smallskip
Arnold studied the relation between singularities and Lagrange maps 
\cite{AGV}(19.). His results (cf. also \cite{Ph1}(4.7.4.1, pp 299--301), 
\cite{Ph2},
\cite{Wi}(Corollary 10)) together with those of chapter 9 yield the following
correspondence between unfoldings and certain germs of F-manifolds.

\begin{theorem}
a) Each irreducible germ of a massive F-manifold with smooth analytic spectrum
is the base space of a semiuniversal unfolding of an isolated 
hypersurface singularity.

b) Suppose, $F_i:(\C^{m_i}\times M_i,0)\to (\C,0)$, $i=1,2,$ are semiuniversal
unfoldings of singularities $f_i:(\C^{m_i},0)\to (\C,0)$ and 
$\varphi:(M_1,0)\to (M_2,0)$ is an isomorphism of the base spaces as
F-manifolds. Suppose that $m_1\leq m_2$.

Then a coordinate change $\psi:(\C^{m_2},0)\to (\C^{m_2},0)$ exists such that
\begin{eqnarray}
f_1(x_1,...,x_{m_1})+x_{m_1+1}^2+...+x_{m_2}^2 
= f_2(x_1,...,x_{m_2})\circ \psi
\end{eqnarray}
and an isomorphism $(\phi,\phi_{base})$ of the unfoldings 
$F_1+x_{m_1+1}^2+...+x_{m_2}^2$ and $F_2\circ \psi$ exists with 
\begin{eqnarray}
F_1+x_{m_1+1}^2+...+x_{m_2}^2 = F_2\circ \psi\circ \phi
\mbox{ \ and \ }
\phi_{base}=\varphi\ .
\end{eqnarray}
\end{theorem}

{\bf Proof:}
a) The restricted Lagrange map of the germ of a massive F-manifold with
smooth analytic spectrum is a miniversal germ of a Lagrange map with 
smooth Lagrange variety (chapter 9). Arnold \cite{AGV}(19.3)
constructs a generating family $F^{(r)}=F^{(r)}(x,t_2,...,t_\mu)$ for it.
Looking at the notions of stable maps and generating families in 
\cite{AGV}(19.), one sees:
$F^{(r)}$ is an $(r)$-semiuniversal unfolding of $F^{(r)}(x,0)$
(cf. Remark 16.2 iv)).
$t_1+F^{(r)}$ is a semiuniversal unfolding of $F^{(r)}(x,0)$. 
Its base space is the given germ of a massive F-manifold.

b) The unfolding $F_i$ is isomorphic to an unfolding 
$t_1+F_i^{(r)}(x_1,...,x_{m_i},t_2,...,t_\mu)$ as in Remark 16.2 iv) over the
same base. Then $F^{(r)}_i$ is an $(r)$-semiuniversal unfolding and a 
generating family for the restricted Lagrange map of the F-manifold
$(M_i,0)$.

The isomorphism $\varphi:(M_1,0)\to (M_2,0)$ induces an isomorphism of the
restricted Lagrange maps. Then the main result in \cite{AGV}(19.4)
establishes a notion of equivalence for $F^{(r)}_1$ and $F^{(r)}_2$,
stable $\RR^+$-equivalence, which yields the desired equivalence in 
Theorem 16.6 b) for $F_1$ and $F_2$.
\hfill $\qed$

\begin{remarks}
i) Two isolated hypersurface singularities $f_i:(\C^{m_i},0)\to (\C,0)$ 
with $m_1\leq m_2$ are {\it stably right equivalent} if a 
coordinate change $\psi:(\C^{m_2},0)\to (\C^{m_2},0)$ with (16.15) exists.
They are right equivalent if furthermore $m_1=m_2$.
The splitting lemma says:
\begin{list}{}{}{\it
\item[]
An isolated hypersurface singularity $f:(\C^m,0)\to (\C,0)$ with 
$r:=m-\rank \bigl({\partial^2 f\over \partial x_i\partial x_j}\bigr)(0)$
is stably right equivalent to a singularity $g:(\C^{r},0)\to (\C,0)$ with
$\rank \bigl({\partial^2 g\over \partial x_i\partial x_j}\bigr)(0)=0$;
this singularity $g$ is unique up to right equivalence.
}\end{list}
(For the existence of $g$ see e.g. \cite{Sl}(4.2 Satz), the uniqueness of 
$g$ up to right equivalence follows from Theorem 16.6 or from Scherk's
result (Remark 16.5 iv)).)

\smallskip
ii) Theorem 16.6 gives a 1-1 correspondence between isolated hypersurface
singularities up to stable right equivalence and irreducible germs of
massive F-manifolds with smooth analytic spectrum.

But the liftability of an isomorphism $\varphi:(M_1,0)\to (M_2,0)$ to 
unfoldings which is formulated in Theorem 16.6 b) is stronger. 
The 1-1 correspondence itself follows already from Theorem 16.6 a) and 
Scherk's result (Remark 16.5 iv)).

\smallskip
iii) The proof of Theorem 16.6 a) is not so difficult. 
If $(M,0)$ is an irreducible germ of a massive F-manifold with analytic 
spectrum $(L,\lambda)\subset T^{*}M$, then a sufficiently generic extension
of a generating function on $(L,\lambda)$ to a function on 
$(T^{*}M,\lambda)$ is already a semiuniversal unfolding over $(M,0)$.
A version different from \cite{AGV}(19.3)
of the precise construction is given by Pham
\cite{Ph1}(4.7.4.1, pp 291--301), following H\"ormander.

\smallskip
iv) Theorem 16.6 b) follows also from \cite{Ph2} (again following 
H\"ormander) and from \cite{Wi}(Corollary 10). To apply Wirthm\"ullers
arguments one has to start with the discriminant $\DD$ and the unit field.
Pham \cite{Ph1}\cite{Ph2} starts with the 
{\it characteristic variety}. That is the cone in $T^{*}M$ of the development
$\widetilde{\DD}\subset \P T^{*}M$ of the discriminant.
\end{remarks}

\smallskip

A semiuniversal unfolding $F:(\C^m\times M,0)\to (\C,0)$ yields data as in 
Corollary 15.2 for the germ $(M,0)$ of an F-manifold:
\begin{eqnarray}
(Z,0)=(\C^m\times M,0),\ \ \alpha_Z=dF\ .
\end{eqnarray}
The semiuniversal unfolding $F$ is also a generating family of $(M,0)$
as germ of an F-manifold in the sense of Definition 15.4.

The following observation says that these two special cases
Corollary 15.2 and Definition 15.4 of the general construction of
F-manifolds in Theorem 15.1 meet only in the case of unfoldings of 
isolated hypersurface singularities.

\begin{lemma}
Let $Z, M,\pi_Z, C, \alpha_Z, \aaa_C$, and $F:Z\to C$ satisfy all the properties
in Corollary 15.2 and Definition 15.4.

Then $C$ is smooth. For any point $p\in M$ the multigerm 
$F:(Z,C\cap\pi_Z^{-1}(p))\to \C$ is isomorphic to a transversal product of versal
unfoldings of the singularities of 
$F|\pi_Z^{-1}(p)$ (cf. Remark 16.5 iii)). The irreducible germs $(M_k,p)$
of F-manifolds in the decomposition $(M,p)=\prod_{k=1}^l(M_k,p)$ are
base spaces of semiuniversal unfoldings of the singularities of 
$F|\pi_Z^{-1}(p)$.
\end{lemma}

{\bf Proof:}
The isomorphism $\aaa_C:\tm \to (\pi_C)_{*}\OO_C$ of Corollary 15.2 restricts
at $p\in M$ to a componentwise isomorphism of algebras
\begin{eqnarray*}
T_pM = \bigoplus_{k=1}^l T_pM_k \to 
\bigoplus_{z\in C\cap \pi_Z^{-1}(p)} \bigl(\mbox{Jacobi algebra of }
F|(\pi_Z^{-1}(p),z)\bigr).
\end{eqnarray*}
One applies Theorem 16.1.
\hfill $\qed$

\section{Boundary singularities}
\setcounter{equation}{0}

The last chapter showed that germs of F-manifolds with smooth analytic
spectrum correspond to isolated hypersurface singularities.
The simplest nonsmooth germ of an analytic spectrum of dimension $n$
is isomorphic to 
\begin{eqnarray*}
(\{(x,y)\in \C^2\ |\ xy=0\},0)\times (\C^{n-1},0)\ .
\end{eqnarray*}
We will see that irreducible germs of massive F-manifolds with such an
analytic spectrum correspond to boundary singularities (Theorem 17.6).
Boundary singularities had been introduced by Arnold \cite{Ar2}.
Because of the similarities to hypersurface singularities we will develop
things parallel to chapter 16.

\bigskip
We consider a germ $(\C^{m+1},0)$ with coordinates $x_0,...,x_m$ always together
with the hyperplane $H:=\{x\in \C^{m+1}\ |\ x_0=0\}$ of the first
coordinate.
A {\it boundary singularity} $(f,H)$ is a holomorphic function germ 
$f:(\C^{m+1},0)\to (\C,0)$ such that $f$ and $f|_H$ have isolated 
singularities at 0. It can be considered as an extension of the
hypersurface singularities $f$ and $f|_H$.

Its Jacobi algebra is
\begin{eqnarray}
\OO_{\C^{m+1},0}/J_{f,H}:=
\OO_{\C^{m+1},0}/\bigl( x_0{\partial f\over \partial x_0}, 
{\partial f\over \partial x_1},....,{\partial f\over \partial x_m}\bigr)\ ,
\end{eqnarray}
its Milnor number $\mu=\mu(f,H):=\dim \OO_{\C^{m+1},0}/J_{f,H}$
satisfies (\cite{Ar2}(\S 3), \cite{Sz}(\S 2))
\begin{eqnarray}
\mu = \mu (f)+\mu (f|_H)\ .
\end{eqnarray}
An unfolding of $(f,H)$ is simply a holomorphic function germ
$F:(\C^{m+1}\times \C^n,0)\to (\C,0)$ such that $F|\C^{m+1}\times \{0\}=f$,
that is , an unfolding of $f$.
Again we write the parameter space as $(M,0)=(\C^n,0)$.

But the critical space $(C,0)\subset (\C^{m+1}\times M,0)$ of $F$ as unfolding
of the boundary singularity $(f,H)$ is the zero set of the ideal
\begin{eqnarray}
J_{F,H}:=\bigl( x_0{\partial F\over \partial x_0}, 
{\partial F\over \partial x_1},....,{\partial F\over \partial x_m}\bigr)
\end{eqnarray}
with the complex structure 
$\OO_{C,0}=\OO_{\C^{m+1}\times M,0}/J_{F,H}|_{(C,0)}$ (cf. \cite{Sz}).
Forgetting the complex structure, $(C,0)$ is the union of the critical
sets $(C^{(1)},0)$ of $F$ and $(C^{(2)},0)$ of $F|_{H\times M}$ as
unfoldings of hypersurface singularities.

For the same reasons as in the hypersurface case the projection
$pr:(C,0)\to (M,0)$ is finite and flat with degree $\mu$ and
$\OO_{C,0}$ is a free $\OO_{M,0}$-module of rank $\mu$. The 1--form
\begin{eqnarray}
\alpha_Z:= -{\partial F\over \partial x_0} dx_0 + dF 
=\sum_{i=1}^m {\partial F\over \partial x_i} dx_i + 
\sum_{j=1}^n {\partial F\over \partial t_j} dt_j
\end{eqnarray}
on $(Z,0):=(\C^{m+1}\times M,0)$ gives rise to a kind of 
Kodaira-Spencer map
\begin{eqnarray}
\aaa_C : \TT_{M,0}\to \OO_{C,0},\ 
X\mapsto \alpha_Z(\widetilde{X})|_{(C,0)}\ ,
\end{eqnarray}
where $\widetilde{X}$ is any lift of $X\in \TT_{M,0}$ to $(Z,0)$. 
It induces a reduced Kodaira-Spencer map
\begin{eqnarray}
\aaa_{C}|_0:T_0M\to \OO_{\C^{m+1},0}/J_{f,H}\ .
\end{eqnarray}
The ideal $J_{F,H}$ and the maps $\aaa_C$ and $\aaa_C|_0$ behave well with 
respect to morphisms of unfoldings, as we will see.

A morphism between two unfoldings $F_1$ and $F_2$ as in chapter 16 of a 
boundary singularity $(f,H)$ is a pair 
$(\phi,\phi_{base})$ of a map germ with (16.4)--(16.6) and 
additionally
\begin{eqnarray}
\phi(H\times M_1)\subset H\times M_2\ .
\end{eqnarray}
Then the first entry of $\phi$ takes the form $x_0\cdot unit\in \OO_{Z,0}$.
Using this one can see with a bit more work than in the 
hypersurface case that the critical spaces behave well with respect to
morphisms:
\begin{eqnarray}
\phi^{*}J_{F_2}=J_{F_1}\mbox{ \ and \ }
(C_1,0)=\phi^{-1}((C_2,0))\ .
\end{eqnarray}
Also the Kodaira-Spencer maps behave as well as in the hypersurface
case. The $\OO_{M_1,0}$-linear maps $d\phi_{base},$ $\aaa_{C_2}$, and
$\phi^{*}|_{(C_2,0)}$ are defined as in (16.7)--(16.9); one finds again
\begin{eqnarray}
\aaa_{C_1} = \phi^{*}|_{(C_2,0)} \circ \aaa_{C_2} \circ d\phi_{base}
\end{eqnarray}
and 
\begin{eqnarray}
\aaa_{C_1}|_0 = \aaa_{C_2}|_0 \circ d\phi_{base}|_0\ .
\end{eqnarray}
Versal and semiuniversal unfoldings of boundary singularities are defined
analogously to the hypersurface case and exist.

\begin{theorem}[\cite{Ar2}]
An unfolding $F:(\C^{m+1}\times M,0)\to (\C,0)$ of a boundary
singularity $(f,H)$, $f:(\C^{m+1},0)\to (\C,0)$, is versal iff the 
reduced Kodaira-Spencer map $\aaa_C|_0:T_0M\to \OO_{\C^{m+1},0}/J_{f,H}$
is surjective. It is semiuniversal iff $\aaa_C|_0$ is an isomorphism.
\end{theorem}

\begin{remarks}
i) $\aaa_C|_0$ is surjective (an isomorphism) iff $\aaa_C$ is surjective
(an isomorphism).

\smallskip
ii) $F(x_0,...,x_m,t_1,...,t_\mu)=f+\sum_{i=1}^\mu m_it_i$ is a 
semiuniversal unfolding of the boundary singularity $(f,H)$ if 
$m_1,...,m_\mu\in \OO_{\C^{m+1},0}$ represent a basis of 
$\OO_{\C^{m+1},0}/J_{f,H}$.

\smallskip
iii) The critical space of an unfolding 
$F:(Z,0)=(\C^{m+1}\times M,0)\to (\C,0)$ of a boundary singularity 
$(f,H)$ is reduced and isomorphic to
$(\{(x,y)\in \C^2\ |\ xy=0\},0)\times (\C^{n-1},0)$ iff
${\partial F \over \partial x_0},...,
{\partial F\over \partial x_m}$ 
represent a generating
system of the vector space
$\mmm_{Z,0}/((x_0)+\mmm_{Z,0}^2)$. This is equivalent to the 
nondegeneracy condition
\begin{eqnarray}
\rank  \left
( \begin{array}{*{2}{c}}
{\partial^2 F \over \partial x_0 \partial x_j } &
{\partial^2 F \over \partial x_0 \partial t_k } \\
{\partial^2 F \over \partial x_i \partial x_j } &
{\partial^2 F \over \partial x_i \partial t_k } 
\end{array} 
\right)_{i,j\geq 1} 
(0)=m+1
\end{eqnarray}
(cf. \cite{NN}). It is satisfied for versal unfoldings.

\smallskip
iv) As in Remark 16.2 v) for hypersurface singularities, 
one can define generalized morphisms between
unfoldings of right equivalent boundary singularities. Again the critical
spaces and Kodaira-Spencer maps behave well.
\end{remarks}

\begin{theorem}
Let $F:(\C^{m+1}\times M,0)\to (\C,0)$ be a semiuniversal unfolding of a 
boundary singularity $(f,H)$. 

The Kodaira-Spencer map 
$\aaa_C:\TT_{M,0}\to \OO_{C,0}$ is an isomorphism and induces a multiplication
on $\TT_{M,0}$. Then $(M,0)$ is an irreducible germ of a massive 
F-manifold with analytic spectrum isomorphic to 
$(\{(x,y)\in \C^2\ |\ xy=0\},0)\times (\C^{\mu-1},0)$. 
The field $E:=\aaa_{C}(F|_{C})$ is an Euler field of weight 1.
\end{theorem}

{\bf Proof:}
Similar to the proof of Theorem 16.3. 
One wants to apply Corollary 15.2 and has to show that $\alpha_Z|_{C_{reg}}$
is exact. The critical space $(C,0)$ as a set is the union of the 
smooth zero sets $(C^{(1)},0)$ of $J_F$ and $(C^{(2)},0)$ of the ideal
$(x_0, {\partial F \over \partial x_1} , ....,
{ \partial F\over \partial x_m })$.
The definition $(17.4)$ of $\alpha_Z$ shows 
\begin{eqnarray}
\alpha_Z|_{(C^{(i)},0)} = dF|_{(C^{(i)},0)} \mbox{ \ \ for }i=1,2\ .
\end{eqnarray}
Therefore $\alpha_Z|_{C_{reg}}$ is exact and $F|_{C}\circ \qqq^{-1}$
is a holomorphic generating function of the analytic spectrum.
\hfill $\qed$

\begin{theorem}
Let $F_i:(\C^{m+1}\times M_i,0)\to (\C,0)$, $i=1,2$, be two semiuniversal
unfoldings of a boundary singularity $(f,H)$.

There exists a unique isomorphism $\varphi:(M_1,0)\to (M_2,0)$ of 
F-manifolds such that $\phi_{base}=\varphi$ for any isomorphism 
$(\phi,\phi_{base})$ of the unfoldings $F_1$ and $F_2$.
\end{theorem}

{\bf Proof:}
Similar to the proof of Theorem 16.4.
\hfill $\qed$

\begin{remarks}
i) Let $F:(\C^{m+1}\times M,0)\to (\C,0)$ be a semiuniversal unfolding
of a boundary singularity $(f,H)$ with critical space 
$(C,0)=(C^{(1)},0)\cup (C^{(2)},0)$. For any $t\in M$ the points in 
$\C^{m+1}\times \{ t\} \cap (C,0)$ split into three classes:
The hypersurface singularities of $F|\C^{m+1}\times \{ t\}$ in 
$C^{(1)}-C^{(2)}$, the hypersurface singularities of 
$F|H\times \{ t\}$ in $C^{(2)}-C^{(1)}$, and the boundary singularities
of $F|\C^{m+1}\times \{ t\}$ in $C^{(1)}\cap C^{(2)}$.

The algebra $\OO_C |\C^{m+1}\times \{ 0\}$
is the direct sum of  their Jacobi algebras.
The reduced Kodaira-Spencer map at $t\in M$ is an isomorphism from 
$T_tM$ to this algebra.

Hence the multigerms of $F$ at $\C^{m+1}\times \{ t\}\cap C^{(1)}$ and of 
$F|H\times M$ at $H\times \{t\} \cap (C^{(2)}-C^{(1)})$ 
together form a transversal union of versal unfoldings of these 
hypersurface and boundary singularities.

The components $(M_k,t)$ of the decomposition $(M,t)=\prod_{k=1}^l(M_k,t)$ 
into irreducible germs of F-manifolds are bases of semiuniversal unfoldings
of the hypersurface and boundary singularities.

\smallskip
ii) The eigenvalues of $E\circ :T_tM\to T_tM$ are by definition of $E$ 
the values of $F$ on $\C^{m+1}\times \{ t\}\cap C$. 
The discriminant of the Euler field is 
\begin{eqnarray}
\DD = \{t\in M\ |\ (\det (E\circ ))(t)=0\} = \pi_C (C\cap F^{-1}(0))\ .
\end{eqnarray}
It is the union of the discriminants of $F$ and $F|H\times M$ as 
unfoldings of hypersurface singularities and coincides with the classical
discriminant of $F$ as unfolding of a boundary singularity 
\cite{Ar2}\cite{Sz}. All the results of chapter 11 apply to this discriminant.
\end{remarks}

\smallskip
Nguyen huu Duc and Nguyen tien Dai studied the relation 
between boundary singularities and Lagrange maps \cite{NN}.
Their results together with chapter 9 yield the following correspondence
between unfoldings of boundary singularities and certain germs of
F-manifolds.

\begin{theorem}
Let $(M,0)$ be an irreducible germ of a massive F-manifold with analytic 
spectrum $(L,\lambda )$ isomorphic to 
\begin{eqnarray*}
(\{(x,y)\in \C^2\ |\ xy=0\},0)\times (\C^{n-1},0)
\end{eqnarray*}
and ordered components $(L^{(1)},\lambda)\cup (L^{(2)},\lambda )=(L,\lambda)$.

a) There exists a semiuniversal unfolding $F$ of a boundary singularity
such that the base space is isomorphic to $(M,0)$ as F-manifold and
the isomorphism $\qqq:(C,0)\to (L,\lambda )$ maps $C^{(i)}$ to 
$L^{(i)}$.

b) Suppose, $F_i:(\C^{m_i+1}\times M_i,0)\to (\C,0)$, $i=1,2$, are 
semiuniversal unfoldings of boundary singularities $(f_i,H_i)$ and
$\varphi:(M_1,0)\to (M_2,0)$ is an isomorphism of the base spaces as
F-manifolds. Suppose that $m_1\leq m_2$.

Then a coordinate change $\psi:(\C^{m_2+1},0)\to (\C^{m_2+1},0)$ with
$\psi((H_2,0))=(H_2,0)$ exists such that
\begin{eqnarray}
f_1(x_0,...,x_{m_1})+x_{m_1+1}^2+...+x_{m_2}^2 
= f_2(x_0,...,x_{m_2})\circ \psi
\end{eqnarray}
and an isomorphism $(\phi,\phi_{base})$ of the unfoldings 
$F_1+x_{m_1+1}^2+...+x_{m_2}^2$ 
and $F_2\circ \psi$ of boundary singularities exists with 
\begin{eqnarray}
F_1+x_{m_1+1}^2+...+x_{m_2}^2 
= F_2\circ \psi \circ \phi \mbox{ \ and \ }
\phi_{base}=\varphi\ .
\end{eqnarray}
\end{theorem}

{\bf Proof:}
a) In \cite{NN} (Proposition 1) an unfolding 
$F:(\C^{m+1}\times M,0) \to (\C,0)$ with nondegeneracy condition (17.11) 
of a boundary singularity is constructed such that $F$ is a generating 
family for $L^{(1)}\subset T^{*}M$ and $F|H\times M$ is a generating 
family for $L^{(2)}$.

One can show that there are canonical maps $C^{(i)}\to L^{(i)}$ which 
combine to an isomorphism $\qqq:C\to L$ with $\aaa_C={\bf \hat{q}}\circ \aaa$
(as in Theorem 15.1). Then the Kodaira-Spencer map 
$\aaa_C:\TT_{M,0}\to \OO_{C,0}$ is an isomorphism and $F$ is a semiuniversal
unfolding of a boundary singularity. 
(Implicitly this is also contained in \cite{NN} (Th\'eor\`eme)).
Because of 
$\aaa_C={\bf \hat{q}}\circ \aaa$ its base is $(M,0)$ as F-manifold.

b) \cite{NN} (Proposition 3).
\hfill $\qed$

\begin{remarks}
i) Two boundary singularities $f_i:(\C^{m_i+1},0)\to (\C,0)$ with 
$m_1\leq m_2$ are {\it stably right equivalent} if a 
coordinate change $\psi$ as in 
Theorem 17.6 b) exists. They are right equivalent if furthermore
$m_1=m_2$. A splitting lemma for boundary singularities is formulated
below in Lemma 17.8.

\smallskip
ii) Theorem 17.6 gives a 1-1 correspondence between boundary singularities
up to stable right equivalence and irreducible germs of massive
F-manifolds with analytic spectrum 
$(L,\lambda)\cong (\{(x,y)\in \C^2\ |\ xy=0\},0)\times (\C^2,0)$
and ordered components $(L^{(1)},\lambda)\cup (L^{(2)},\lambda)
=(L,\lambda )$.

\smallskip
iii) Interchanging the two components of $(L,\lambda )$ corresponds to a
duality for boundary singularities which goes much further and has been
studied by I. Shcherbak, A. Szpirglas
\cite{Sz}\cite{SS1}\cite{SS2}, and others.
\end{remarks}

\begin{lemma} (Splitting lemma for boundary singularities) \\
A boundary singularity $(f,H)$ with $f:(\C^{m+1},0)\to (\C,0)$ and 
$H=\{x\ |\ x_0=0\}$ is stably right equivalent to a boundary singularity
$g:(\C^{r+1},0)\to (\C,0)$ in 
\begin{eqnarray}
r+1=\max (2;m+1-\rank \bigl( 
{\partial^2 f (x_0^2,x_1,...,x_m) \over 
\partial x_i \partial x_j }
\bigr) (0)) 
\end{eqnarray}
coordinates. 
The boundary singularity $g$ is unique up to right equivalence.
\end{lemma}

{\bf Proof:}
Existence of $g$: The group $G=\Z_2$ acts on $(\C^{m+1},0)$ by 
$(x_0,x_1,...,x_m)\mapsto (\pm x_0,x_1,...,x_m)$. 
Boundary singularities on the quotient $(\C^{m+1},0)$ correspond 
to $G$-invariant singularities on the double cover, branched along 
$H$ (\cite{AGV} 17.4). 

One applies an equivariant splitting lemma of Slodowy \cite{Sl}(4.2 Satz)
to the $G$-invariant singularity $f(x_0^2,x_1,...,x_m)$.
The nondegenerate quadratic part of the $G$-invariant singularity in 
splitted form does not contain $x_0^2$ because $f$ is not smooth.

Uniqueness of $g$: This follows with Theorem 17.6 b).
\hfill $\qed$

\bigskip

The following two observations give some information on generating families
in the sense of Definition 15.4 for the F-manifolds of boundary 
singularities. The first one gives a distinguished generating family
and is essentially well known.
The second one explains why $B_2=I_2(4)$ is missing in Example 15.5 b).

\begin{lemma}
a) Let $F:(Z,0)=(\C^{m+1}\times M,0)\to (\C,0)$ be a semiuniversal unfolding
of a boundary singularity $(f,H)$. 

Then the function
$\widetilde{F}:(\widetilde{Z},0)=(\C^{m+1}\times M,0)\to (\C,0)$ 
with
$\widetilde{F}(x,t)=F(x_0^2,x_1,...,x_m,t)$ is a generating family for the
germ $(M,0)$ of an F-manifold. 

The finite map $\widetilde{\qqq}:\widetilde{C}\to 
L = L^{(1)}\cup L^{(2)}$ from its critical set 
$\widetilde{C}$ to the analytic spectrum $L$ has degree 2 on $L^{(1)}$ and
degree 1 on $L^{(2)}$. The branched covering $\widetilde{C}\to M$ has
degree $2\mu(f)+\mu(f|_H)$.

b) Let $(M,0)$ be a germ of a massive F-manifold with analytic spectrum
$(L,\lambda)\cong (\{(x,y)\in \C^2\ |\ xy=0\},0)\times (\C^{n-1},0)$.
There does not exist a generating family $F:(Z,0)\to (\C,0)$ 
with critical set $C$ such that the canonical map 
$\qqq:C\to L$ 
is a homeomorphism.
\end{lemma}

{\bf Proof:}
a) Consider the branched covering $\pi_G:\widetilde{Z}\to Z$,
$(x_0,...,x_m,t)\mapsto (x_0^2,x_1,...,x_m,t)$ which is induced by the
action $(x_0,...,x_m,t)\mapsto (\pm x_0,x_1,...,x_m,t)$ 
of the group $G=\Z_2$ on 
$\widetilde{Z}$.
The composition $\widetilde{F}=F\circ \pi_G$ is an unfolding of the 
$G$-invariant singularity $\widetilde{F}|(\C^{m+1}\times \{0\},0)$,
in fact, semiuniversal within the $G$-invariant unfoldings 
(cf. \cite{Sl}(4.5)). The ideals
\begin{eqnarray}
J_{\widetilde{F}}= \bigl( 
2x_0\cdot ( { \partial F \over \partial x_0 }\circ \pi_G),
{\partial F \over \partial x_1} \circ \pi_G,...,
{\partial F \over \partial x_m} \circ \pi_G \bigl)
\end{eqnarray}
and $\pi_G^{*} J_{F,H}$ have the same zero sets
$\widetilde{C} = \pi_G^{-1}(C)$.

Comparison of $\alpha_{\widetilde{Z}} =d\widetilde{F}$ on $\widetilde{Z}$
and $\alpha_Z = -{ \partial F\over \partial x_0 }dx_0 +dF $ on $Z$ 
(formula (17.4)) 
shows that the map $\aaa_{\widetilde{C}}:\TT_{M,0} \to \OO_{\widetilde{C},0}$
factorizes into the Kodaira-Spencer map 
$\aaa_C:\TT_{M,0}\to \OO_{{C},0}$ and the map
$\pi_G^{*}|_{(C,0)} : \OO_{C,0} \to \OO_{\widetilde{C},0}$.

Therefore $\aaa_{\widetilde{C}}$ is injective with multiplication invariant
image and induces the correct multiplication on $\TT_{M,0}$. The rest 
is clear.

b) Assume that such a generating family $F$ exists.
$(L,\lambda )$ is its own maximalisation. Therefore the homeomorphism $\qqq$
is an isomorphism. 

Then $\aaa_C={\bf \hat{q}}\circ \aaa$ 
(cf. Definition 15.4 and Theorem 15.1) is an isomorphism. 
We are simultaneously in the special cases Definition 15.4 and 
Corollary 15.2 of Theorem 15.1. By Lemma 16.8
$(L,\lambda)$ is smooth, a contradiction.
\hfill $\qed$

\section{Coxeter groups and F-manifolds}
\setcounter{equation}{0}

The complex orbit space of an irreducible Coxeter group is equipped
with the discriminant, the image of the reflection hyperplanes,
and with a certain distinguished vector field (see below), which
is unique up to a scalar. Together they induce as in Corollary 11.6 
the structure of an F-manifold on the complex orbit space (Theorem 18.1).
This follows independently from \cite{Du1}\cite{Du2}(Lecture 4) and
from \cite{Gi2}(Theorem 14).

In fact, both give stronger results. Dubrovin established the structure
of a Frobenius manifold. This will be discussed in chapter 19.
Givental proved that these F-manifolds are distinguished by certain
geometric conditions (Theorem 18.4). 
With one additional argument we will show that
the germs of these F-manifolds and their products are the only germs
of simple F-manifolds whose tangent spaces are Frobenius algebras
(Theorem 18.3). This complements in a nice way the relation between
Coxeter groups and simple hypersurface and boundary singularities.

We will also present simple explicit formulas for these F-manifolds
which are new for $H_3$ and $H_4$ (Theorem 18.5).

\bigskip

A Coxeter group is a finite group $W$ of linear transformations of the
Euclidean space $\R^n$ generated by reflections in hyperplanes.
Each Coxeter group is the direct sum of irreducible Coxeter groups.
They are (\cite{Co}\cite{Bo}) $A_n$ ($n\geq 1$), $D_n$ ($n\geq 4$), $E_6$, 
$E_7$, $E_8$, $B_n$, ($n\geq 2$), $F_4$, $G_2$, $H_3$, $H_4$, $I_2(m)$
($m\geq 3$) with $A_2=I_2(3)$, $B_2=I_2(4)$, $H_2:=I_2(5)$, $G_2=I_2(6)$.

The Coxeter group $W$ acts on $\C^n=\R^n\otimes_\R \C$ and on 
$\C[x_1,...,x_n]$, where $x_1,...,x_n$ are the coordinates on $\C^n$.
The ring $\C[x_1,...,x_n]^W$ of invariant polynomials is generated by 
$n$ algebraically independent homogeneous polynomials $P_1,...,P_n$.
Their degrees $d_i:=\deg P_i$ are unique (up to ordering).
The quotient $\C^n/W$ is isomorphic to $\C^n$ as affine algebraic variety.
The $\C^{*}$-action and the vector field 
$\sum_i x_i{\partial \over \partial x_i}$ 
on the original $\C^n$ induce a $\C^{*}$-action and a vector field
$\sum d_i t_i{\partial \over \partial t_i}$ 
on the  orbit space $\C^n/W \cong \C^n$. 
The image in the orbit space of the union of the reflection hyperplanes
is the discriminant $\DD$ of the Coxeter group.

Suppose for a moment that $W$ is irreducible. Then there is precisely one 
highest degree, which is called the Coxeter number $h$. The degrees can be
ordered to satisfy
\begin{eqnarray}
d_1=h > d_2 \geq ... \geq d_{n-1}>d_n=2\ ,\\
d_i+d_{n+1-i}=h+2\ .
\end{eqnarray}
The vector field $e:= {\partial \over \partial t_1}$
is unique up to a scalar.

\begin{theorem}
The complex orbit space $M:= \C^n/W \cong \C^n$ of an irreducible 
Coxeter group $W$ carries a unique structure of a massive F-manifold
with the unit field 
$e={\partial \over \partial t_1}$ and the discriminant $\DD$. The discriminant 
$\DD$ corresponds to the Euler field
\begin{eqnarray}
E={1\over h}\sum_{i=1}^n d_it_i
{\partial \over \partial t_i}
\end{eqnarray}
of weight 1.
\end{theorem}

{\bf Proof:}
The uniqueness follows from Corollary 11.6. 

The existence follows from Dubrovin's result 
(\cite{Du1}\cite{Du2} (Lecture 4), cf. Theorem 19.1)
or Givental's result \cite{Gi2}(Theorem 14) together with Theorem 9.4.

Below in Theorem 18.5 we will follow Givental and reduce it to 
classical results on the appearance of discriminants in singularity theory
(\cite{Bri}\cite{Ar2}\cite{Ly2}\cite{ShO}).
\hfill $\qed$

\begin{remarks}
i) Corollary 11.6 gives probably the most elementary way how $e$ and $\DD$
determine the multiplication on the complex orbit space $M=\C^n/W\cong \C^n$,
at least at a generic point: the $e$-orbit of a generic point $p\in M$ 
intersects $\DD$ transversally in $n$ points. One shifts the tangent spaces
of $\DD$ at these points with the flow of $e$ to $T_pM$.
Then there exists a basis $e_1,...,e_n$ of $T_pM$ such that 
$\sum_{i=1}^n e_i=e$ and such that the hyperplanes 
$\bigoplus_{i\neq j} \C\cdot e_i$, $j=1,...,n$, are the shifted tangent spaces
of $\DD$. The multiplication on $T_pM$ is given by 
$e_i\circ e_j=\delta_{ij}e_i$.

\smallskip
ii) The unit field $e={\partial \over \partial t_1}$ 
is only unique up to a scalar. The flow of the Euler field respects the
discriminant $\DD$ and maps the unit field $e$ and the multiplication
to multiples, because of $\lie_E(e)=-e$ and $\lie_E(\circ)=\circ$.

Therefore the isomorphism class of the F-manifold $(M,\circ,e,E)$ 
is independent of the choice of the scalar.

\smallskip
iii) The complex orbit space of a reducible Coxeter group $W$ is 
isomorphic to the product of the complex orbit spaces of the irreducible
subgroups. The discriminant decomposes as in Remark 11.2 iv). Now any
sum of unit fields for the components yields a unit field for 
$\C^n/W$. The choices are parametrized by 
$(\C^{*})^{|irr. subgroups|}$. But the resulting F-manifold is unique 
up to isomorphism. It is the product of F-manifolds for the irreducible
subgroups. This F-manifold and its germ at 0 will be 
denoted by the same combination of letters as the Coxeter group.
\end{remarks}

\begin{theorem}
Let $((M,p),\circ,e)$ be a germ of a massive F-manifold.

$((M,p),\circ,e)$ is simple and $T_pM$ is a Frobenius algebra if and only if
$((M,p),\circ,e)$ is isomorphic to the germ at 0 of an F-manifold
of a Coxeter group.
\end{theorem}

This builds on the following result, which is a reformulation with chapter 9
of a theorem of Givental \cite{Gi2}(Theorem 14). 
The Theorems 18.3, 18.4, and 18.5 will be proved below in the opposite
order. Some arguments
on $H_3$ and $H_4$ in the proof of Theorem 18.4 will only be sketched.

\begin{theorem} {\rm (Givental)} \ 
a) The F-manifold of an irreducible Coxeter group is simple. 
The analytic spectrum $(L,\lambda)$ of its germ at 0 is isomorphic to
$(\{(x,y)\in \C^2\ |\ x^2=y^r\},0)\times (\C^{n-1},0)$  
with $r=1$ for $A_n, D_n, E_n$, $r=2$ for $B_n, F_4$, $r=3$ for 
$H_3,H_4$ and $r=m-2$ for $I_2(m)$.

b) An irreducible germ of a simple F-manifold with analytic spectrum
isomorphic to a product of germs of plane curves is isomorphic to the 
germ at 0 of an F-manifold of an irreducible Coxeter group.
\end{theorem}

Finally, we want to present the F-manifolds
of the irreducible Coxeter groups explicitly with data as in Corollary 15.2.
We will use the notations of Corollary 15.2. The following is
a consequence of results in \cite{Bri}\cite{Ar2}\cite{Ly2}\cite{ShO}
on the appearance of the discriminants of Coxeter groups in singularity theory.

\begin{theorem}
a) The germs at 0 of the F-manifolds of the Coxeter groups $A_n, D_n, E_n,
B_n, F_4$ are isomorphic to the base spaces of the semiuniversal 
unfoldings of the corresponding simple hypersurface singularities
$A_n, D_n, E_n$ and simple boundary singularities $B_n$ (or $C_n$) and $F_4$.

b) For the F-manifolds $(M,\circ,e)=(\C^n,\circ,e)$ of the irreducible 
Coxeter groups,
a space $Z$ with projection $\pi_Z:Z\to M$, a subspace $C\subset Z$ and 
a 1--form $\alpha_Z$ will be given such that the map
\begin{eqnarray}
\aaa_C:\tm \to (\pi_C)_{*}\OO_C,\ 
X\mapsto \alpha_Z(\widetilde{X})|_C
\end{eqnarray}
is welldefined and an isomorphism of $\OO_M$-algebras.
$C$ is isomorphic to the analytic spectrum of $(M,\circ,e)$.
The Euler field is always 
$E={1\over h}\sum_{i=1}^n d_it_i{\partial \over \partial t_i}$.
The discriminant $\DD\subset M$ is $\DD=\pi_C(\aaa_C(E)^{-1}(0))$.

\smallskip
i)
$A_n, B_n, H_3, I_2(m): \\ 
Z=\C\times M=\C\times \C^n$ with coordinates
$(x,t)=(x,t_1,...,t_n)$, \\
$\alpha_Z=dt_1+xdt_2+...+x^{n-1}dt_n$, \\ 
$\widetilde{t_2}(x,t):= t_2+2xt_3+...+(n-1)x^{n-2}t_n$,

\smallskip \noindent
$A_n: \ \ C=\{(x,t)\in Z\ |\ x^n-\widetilde{t_2}=0\} ,$\\
$B_n: \ \ C=\{(x,t)\in Z\ |\ x\cdot (x^{n-1}-\widetilde{t_2})=0\} ,$\\
$H_3: \ \ C=\{(x,t)\in Z\ |\ x^3-\widetilde{t_2}^2=0\} ,$\\
$I_2(m): \ C=\{(x,t)\in Z\ |\ x^2-\widetilde{t_2}^{m-2}=0\} .$

\smallskip
ii)
$D_4, F_4, H_4: \\  
Z=\C^2\times M=\C^2\times \C^4$ with coordinates
$(x,y,t)=(x,y,t_1,...,t_4)$, \\
$\alpha_Z=dt_1+xdt_2+ydt_3+xydt_4$, \\ 
$\widetilde{t_2}(x,y,t):= t_2+yt_4$,
$\widetilde{t_3}(x,y,t):= t_3+xt_4$,

\smallskip \noindent
$D_4: \ \ C=\{(x,y,t)\in Z\ |\ x^2+\widetilde{t_2}=0, 
                               y^2+\widetilde{t_3}=0\} ,$\\
$F_4: \ \ C=\{(x,y,t)\in Z\ |\ x^2+\widetilde{t_2}=0, 
                               y^2+\widetilde{t_3}^2=0 \} ,$\\
$H_4: \ \ C=\{(x,y,t)\in Z\ |\ x^2+\widetilde{t_2}=0, 
                               y^2+\widetilde{t_3}^3=0 \} .$

\smallskip
iii)
$D_n, E_6, E_7, E_8: \\ 
Z=\C^2\times M=\C^2\times \C^n$ with coordinates
$(x,y,t)=(x,y,t_1,...,t_n)$,  \\
$F:Z\to \C$ a semiuniversal unfolding of $F|\C^2\times \{0\}$, \\ 
$\alpha_Z= \sum_{i=1}^n 
{\partial F\over \partial t_i} dt_i$
(or $\alpha_Z=dF$), \\
$C=\{(x,y,t)\in Z\ |\ 
{\partial F \over \partial x} = {\partial F\over \partial y}=0\}$, 

\smallskip \noindent
$D_n: \ \ F= x^{n-1}+xy^2+\sum_{i=1}^{[n/2]} x^{i-1}t_i + yt_{[n/2]+1}
               +\sum_{i=[n/2]+2}^n x^{i-2}t_i$,\\
$E_6: \ \ F= x^4+y^3+t_1+xt_2+yt_3+x^2t_4+xyt_5+x^2yt_6$,\\
$E_7: \ \ F= x^3y+y^3+t_1+xt_2+yt_3+x^2t_4+xyt_5+x^3t_6+x^4t_7$,\\
$E_8: \ \ F= x^5+y^3+t_1+xt_2+yt_3+x^2t_4+xyt_5+x^3t_6+x^2yt_7+x^3yt_8$.
\end{theorem}

{\bf Proof of Theorem 18.5:}
a) One can choose a semiuniversal unfolding 
$F=f(x_1,...,x_m)+\sum_{i=1}^nm_it_i$ of the hypersurface or boundary 
singularity which is weighted homogeneous with positive degrees in 
all variables and parameters. There is an isomorphism from its base space
$\C^n$ to the complex orbit space of the corresponding Coxeter group
which respects the discriminant, the Euler field, and the unit field
(\cite{Bri}\cite{Ar2}). It respects also the F-manifold structure
(Corollary 11.6).

\smallskip
b) i) for $I_2(m)$ is Remark {15.5 a).} i) for $A_n, B_n$ and ii) for $D_4$ 
follow with a), with semiuniversal unfolding as in a) for the singularities
$-{1\over n+1}x^{n+1}$ ($A_n$), $(-{1\over n}x^n+y^2,H=\{x=0\})$ ($B_n$),
${1\over 3}x^3+{1\over 3}y^3$ ($D_4$). Also iii) follows with a).

The same procedure gives for the boundary singularity 
$({1\over 2}y^2+{1\over 3}x^3,H=\{y=0\})$ ($F_4$) the data in ii) with
critical set 
\begin{eqnarray}
C' = \{(x,y,t)\in Z\ |\ x^2+\widetilde{t_2}=0,y^2+y\widetilde{t_3}=0\} .
\end{eqnarray}
It is a nontrivial, but solvable exercise to find compatible automophisms
of $Z$ and $M$ which map $C'$ to $C$ and $\alpha_Z$ to $\alpha_Z$ 
modulo $I_C\Omega_Z^1$.
Independently of explicit calculations, the proof of Theorem 18.4 will show
that the data $(Z,\alpha_Z,C)$ in ii) correspond to $F_4$.

The data in i) for $H_3$ and in ii) for $H_4$ can be obtained from results
of O.P. Shcherbak (\cite{ShO} pp 162, 163; \cite{Gi2} Proposition 12)
(for $H_3$ one could use instead \cite{Ly2}). The unfoldings
\begin{eqnarray}
F_{H_3} = \int_0^y (u^2+x)^2du +t_1+xt_2+x^2t_3
\end{eqnarray}
of $D_6$ and 
\begin{eqnarray}
F_{H_4} = \int_0^y (u^2+t_3+xt_4)^2du +x^3+t_1+xt_2
\end{eqnarray}
of $E_8$ have only critical points with even Milnor number and are 
maximal with this property. Their discriminants are isomorphic to the 
discriminants of the Coxeter groups $H_3$ and $H_4$. The unfoldings are
generating families in the sense of Definition 15.4 for the F-manifolds 
of the Coxeter groups $H_3$ and $H_4$.
We will determine the data in ii) for $H_4$ from $F_{H_4}$; 
the case $H_3$ is similar.

Consider the map
\begin{eqnarray}
\phi:\C^2\times \C^4 \to Z=\C^2\times \C^4,\ 
(x,y,t)\mapsto (x,\widetilde{y},t), \\
\widetilde{y}(x,y,t):= \int_0^y 2(u^2+t_3+xt_4)du 
= {2\over 3}y^3+2(t_3+xt_4)y\ ,
\end{eqnarray}
and observe
\begin{eqnarray}
dF_{H_4} &=& (3x^2+t_2+\widetilde{y}t_4)dx + (y^2+t_3+xt_4)^2dy\nonumber \\
&& +dt_1+xdt_2+\widetilde{y}dt_3+x\widetilde{y}dt_4\ ,
\end{eqnarray}
\begin{eqnarray}
{9\over 4}\widetilde{y}^2 + 4(t_3+xt_4)^3 = 
(y^2+t_3+xt_4)^2(y^2+4(t_3+xt_4))\ .
\end{eqnarray}
Therefore
\begin{eqnarray}
\phi^{*}(\alpha_Z)=dF 
- {\partial F\over \partial x}dx 
- {\partial F \over \partial y}dy
\end{eqnarray}
and the image under $\phi$ of the reduced critical set $C_F$ of $F_{H_4}$ is
\begin{eqnarray}
\phi (C_F) = \{(x,y,t)\in Z &|& 3x^2+t_2+yt_4=0,\\
&&{9\over 16}y^2+t_3+xt_4=0\}\ .\nonumber
\end{eqnarray}
An automorphism 
\begin{eqnarray}
Z\to Z, \ (x,y,t_1,t_2,t_3,t_4)\mapsto
(r^{-1}x,s^{-1}y,t_1,rt_2,st_3,rst_4)
\end{eqnarray}
for suitable $r,s\in \C^{*}$ maps 
$\phi(C_F)$ to $C$ and respects $\pi_Z$ and $\alpha_Z$, together with the
induced automorphism $M\to M$.
\hfill $\qed$

\bigskip

{\bf Sketch of the proof of Theorem 18.4:}
a) Consider the data in Theorem 18.5 b). The Euler field on $M=\C^n$
is $E={1\over h}\sum_id_it_i{\partial \over \partial t_i}$.
The coefficients of the Lyashko-Looijenga map $\Lambda:M\to \C^n$
are up to a sign the symmetric polynomials in the eigenvalues of $E\circ$.
Because of $\lie_E(E\circ)=E\circ$, the coefficient $\Lambda_i$
is weighted homogeneous of degree $i$ with respect to the weights 
$({d_1\over h},...,{d_n\over h})$ for $(t_1,...,t_n)$.
The Lyashko-Looijenga map is a branched covering of degree
\begin{eqnarray}
n!\cdot (\prod {d_i\over h})^{-1} = n!h^n\cdot |W|^{-1}
\end{eqnarray}
and $(M,\circ,e)$ is simple (Corollary 14.4 b)).
The analytic spectrum is isomorphic to $C$.

\smallskip
b) The dimension $\dim (M_k,p)$ of an irreducible germ in the 
decomposition $(M,p)=\prod_{k=1}^l (M_k,p)$ of a germ of a massive
F-manifold is equal to the intersection multiplicity of 
$T^{*}_pM$ with the corresponding germ $(L,\lambda_k)$ of the 
analytic spectrum $L$. This number will be called the intersection
multiplicity $IM(\lambda_k)$.

$(S_\mu(q),q)$ denotes for any $q\in M$ the $\mu$-constant stratum
through $q$ (chapter 14), and $l(q)$ the number of irreducible components of 
$(M,q)$. For any subvariety $S\subset L$ we have the estimates 
\begin{eqnarray}
&&\max (l(q)|q\in \pi(S)) \leq n+1-\min (IM(\sigma )| \sigma \in S)\ ,\\
&&\max (\dim (S_\mu(q),q)\ |\ q\in S) \geq \dim S\ ,\\
&&\max (\modmu(M,q)\ |\ q\in S) \nonumber \\
&&\qquad + n+1-\min (IM(\sigma )\ |\ \sigma \in S)
\geq \dim S\ .
\end{eqnarray}
Therefore, if $M$ is simple then
\begin{eqnarray}
\min (IM(\sigma )\ |\ \sigma\in S) \leq n+1-\dim S
\end{eqnarray}
for any subvariety $S\subset L$.

Now suppose that $((M,p),\circ,e)$ is an irreducible germ of a 
simple F-manifold and that 
$\phi:(M,p)\to \prod_{i=1}^n (C_i,0)$ is an isomorphism to a product of germs
of plane curves (they are necessarily plane because of Proposition 7.4).

If at least two curve germs were not smooth, e.g. $(C_{n-1},0)$ and $(C_n,0)$,
then the intersection multiplicities $IM(p)$ for points $p$ in 
$S_1:= \phi^{-1}(\prod_{i=1}^{n-2}(C_i,0)\times \{0\})$
were at least 4; but $\dim (S_1,p)=n-2$, a contradiction to (18.19).
So, at most one curve, e.g. $(C_n,0)$, is not smooth.

The irreducible germs of F-manifolds which correspond to generic points of
$\pi (S_2)$ for $S_2:= \phi^{-1}(\prod_{i=1}^{n-1}(C_i,0)\times \{0\})$
are at most 2-dimensional because of (18.19).
Therefore
\begin{eqnarray}
(C_n,0)\cong (\{(x,y)\in \C^2\ |\ x^2=y^r\},0)
\end{eqnarray}
for some $r\in \N$. If $r\geq 4$ and $n\geq 3$
then the set of possible intersection 
multiplicities for points in $S_2$ has a gap at 3 and a subvariety 
$S_3\subset S_2$ exists with $\dim S_3=n-2$ and 
$\min (IM(\sigma)\ |\ \sigma\in S_3)\geq 4$ \cite{Gi2}(p 3266), 
a contradiction to (18.19).
Therefore $r\in \{ 1,2,3\}$ or $n\leq 2$.

If $r\in \{1,2\}$ then $(M,p)$ is the base space of a semiuniversal unfolding
of a hypersurface singularity ($r=1$, Theorem 16.7, \cite{AGV}(19.))
or boundary singularity ($r=2$, Theorem 17.7, \cite{NN}).
Simplicity of their F-manifolds corresponds to simplicity of the 
singularities. The simple hypersurface singularities are 
$A_n, D_n, E_6, E_7, E_8$ \cite{AGV}.
The simple boundary singularities are $B_n, C_n$, and 
$F_4$ \cite{Ar2}\cite{AGV}.
$B_n$ and $C_n$ are dual boundary singularities and have isomorphic
discriminants and F-manifolds.

The details of the case $r=3$ (\cite{Gi2} pp 3269--3271) are difficult and
will not be given here. In that case the set of possible intersection 
multiplicities for points in $S_2$ has a gap at 5. If $n\geq 6$ then
a subvariety $S_4\subset S_2$ exists with $\dim S_4=n-4$ and 
$\min (IM(\sigma)\ |\ \sigma\in S_4)\geq 6$, 
a contradiction to (18.19). The case $r=3$, $n=3$ 
corresponds to $H_3$,
the case $r=3$, $n=4$ corresponds to $H_4$.
\hfill $\qed$

\bigskip

{\bf Proof of Theorem 18.3:}
It is sufficient to consider an irreducible germ $((M,p),\circ, e)$.
If it corresponds to a Coxeter group then it is simple (Theorem 18.4 a))
and $T_pM$ is a Frobenius algebra (Theorem 18.5 b)).

Suppose that $(M,p)$ is simple and that $T_pM$ is a Frobenius algebra.
We will show by induction on the dimension $n=\dim M$ that the analytic
spectrum $(L,\lambda)$ is isomorphic to
$(\{(x,y)\in \C^2\ |\ x^2=y^r\},0)\times (\C^{n-1},0)$ 
for some $r\in \N$.

This is clear for $n=2$. Suppose that $n\geq 3$. 
The maximal ideal of $T_pM$ is called $\mmm$. The socle
$\Ann_{T_pM}(\mmm)$ of the Gorenstein ring $T_pM$ has dimension 1, therefore
$\Ann_{T_pM}(\mmm)\ { \genfrac{}{}{0pt}{}{\subset}{\neq}} 
\ \mmm$ and $\mmm^2\neq 0$.
In the equations for the analytic spectrum $(L,\lambda)\subset T^{*}_pM$
one can eliminate fiber coordinates which correspond to 
$\mmm^2\subset T_pM$: the embedding dimension of $(L,\lambda)$ is
\begin{eqnarray}
\embdim (L,\lambda) \leq n + \dim {\mmm \over \mmm^2} 
\leq 2n-2 
\end{eqnarray}
(Lemma 15.3). Then $(L,\lambda)\cong (\C^2,0)\times (L'',\lambda'')$
(Propositiom 7.4). There exists $\lambda_2\in L$ close to $\lambda$ such that
$(L,\lambda_2)\cong (L,\lambda)$ and $\pi(\lambda_2)$ is not in the
$e$-orbit of $p$.

Now for all $q$ near $p$, but outside of the $e$-orbit of $p$ the germ $(M,q)$
is reducible because of $\modmu (M,p)=0$.
For all $q$ near $p$ the germ $(M,q)$ is simple and $T_qM$ is a Frobenius 
algebra (Lemma 1.2). 

One can apply the induction hypothesis to the irreducible component
of $(M,\pi(\lambda_2))$ which corresponds to $\lambda_2$.
Its analytic spectrum $(L',\lambda')$ is isomorphic to a product of a smooth
germ and a curve as above. Now 
$(L,\lambda)\cong (L,\lambda_2)\cong (\C^{n-\dim L'},0)\times (L',\lambda')$.
One applies Theorem 18.4 b).
\hfill $\qed$

\section{Coxeter groups and Frobenius manifolds}
\setcounter{equation}{0}

K. Saito \cite{SaK2} introduced a flat metric on the complex orbit space
of an irreducible Coxeter group. Dubrovin \cite{Du1}\cite{Du2}(Lecture 4)
showed that this metric and the multiplication and the Euler field 
from Theorem 18.1 together 
yield the structure of a massive Frobenius manifold on the complex 
orbit space (Theorem 19.1).

The Euler field has positive degrees. Dubrovin \cite{Du1}\cite{Du2}(p 268)
conjectured that these Frobenius manifolds and products of them are the 
only massive Frobenius manifolds with an Euler field with positive degrees.

We will prove this conjecture (Theorem 19.3). Theorem 18.3,
which builds on Givental's result (Theorem 18.4, \cite{Gi2}(Theorem 14)),
will be crucial.

\bigskip
We use the same notations as in chapter 18.
A metric on a complex manifold is a nondegenerate complex bilinear form
on the tangent bundle. The flat standard metric on $\C^n$ is invariant
with respect to the Coxeter group $W$ and induces a flat metric 
$\check{g}$ on $M-\DD$. Dubrovin proved the following with differential
geometric tools \cite{Du1}\cite{Du2}(Lecture 4 and pp 191, 195).

\begin{theorem} {\rm (Dubrovin)}
Let $W$ be an irreducible Coxeter group with complex orbit space 
$M=\C^n/W$, Euler field $E$, a unit field $e$, and a multiplication $\circ$
on $M$ as in Theorem 18.1.

The metric $g$ on $M-\DD$ with 
\begin{eqnarray}
g(X,Y):=\check{g}(E\circ X,Y)
\end{eqnarray}
for any (local) vector fields $X$ and $Y$ extends to a flat metric on $M$
and coincides with K. Saito's flat metric.
$(M,\circ,e,E,g)$ is a Frobenius manifold. The Euler field satisfies
\begin{eqnarray}
\lie_E(g)=(1+{2\over h})g\ .
\end{eqnarray}
There exists a basis of flat coordinates $z_1,...,z_n$ on $M$ with 
$z_i(0)=0$ and 
$e={\partial \over \partial z_1}$ and
\begin{eqnarray}
E={1\over h}\sum d_i \cdot z_i{\partial \over \partial z_i}\ .
\end{eqnarray}
\end{theorem}

\begin{remarks}
i) K. Saito (and also Dubrovin) introduced the flat metric $g$ in a way
different from formula (19.2): The metrics $\check{g}$ and $g$ on 
$M-\DD$ induce two isomorphisms $T(M-\DD)\to T^{*}(M-\DD)$. One lifts
$\check{g}$ and $g$ with the respective isomorphism to metrics
${\check{g}}^{*}$ and $g^{*}$ on the cotangent bundle $T^{*}(M-\DD)$.
Then
\begin{eqnarray}
g^{*}=\lie_e({\check{g}}^{*})
\end{eqnarray}
(\cite{Du2} pp 191, 195).
(Here ${\check{g}}^{*}$ and $g^{*}$ are considered as $(0,2)$-tensors.)
K. Saito introduced $g$ with formula (19.4).

\smallskip
ii) Closely related to (19.2) and (19.4) is 
(\cite{Du2} pp 191, 270)
\begin{eqnarray}
\sum_{i=1}^n {\partial Q_1 \over \partial x_i }\cdot 
{\partial Q_2 \over \partial x_i} 
= {\check{g}}^{*} (dQ_1,dQ_2) = i_E(dQ_1\circ dQ_2)\ .
\end{eqnarray}
Here $Q_1,Q_2\in \C[x_1,...,x_n]^W$ are $W$-invariant polynomials; 
$dQ_1$ and $dQ_2$ are interpreted as sections in $T^{*}M$; the 
multiplication $\circ $ is lifted to $T^{*}M$ with the isomorphism 
$TM\to T^{*}M$ induced by $g$; $i_E$ is the contraction of a 1--form
with $E$. 
 
The first equality is trivial. 
(19.5) is related to Arnold's convolution of invariants 
(\cite{Ar3}\cite{Gi1}).

\smallskip
iii) A Frobenius manifold as in Theorem 19.1 for an irreducible Coxeter group
is not unique as the unit field and the multiplication are not unique. 
Contrary to the F-manifold, it is not even unique up to isomorphism.
There is one complex
parameter between $(M,\circ, e)$ and $(M,g)$ to be chosen:
$(M,\circ,e,E,c\cdot g)$ respectively $(M,c\cdot \circ, c^{-1}\cdot e,E,g)$
is a Frobenius manifold for any $c\in \C^{*}$.

\smallskip
iv) We consider only Frobenius manifolds with an Euler field which is
normalized by $\lie_E(\circ )=1\cdot \circ $ (compare Remark 5.5 c)).
The product $\prod M_i$ of Frobenius manifolds $(M_i,\circ_i,e_i,E_i,g_i)$
carries also the structure of a Frobenius manifold if 
$\lie_{E_i}(g_i)=D\cdot g_i$ holds with the same number $D\in \C$ for all
$i$.
This follows from Proposition 4.1, Theorem 5.3 and Remark 5.5 c)
(compare also \cite{Du2}(p 136)).

\smallskip
v) Especially, the complex orbit space $\C^n/W$ of a reducible 
Coxeter group can be provided with the structure of a Frobenius manifold
if the irreducible Coxeter subgroups have the same Coxeter number.

The Frobenius manifold is not unique. The different choices are 
parametrized by 
$(\C^{*})^{|irr. subgroups|}$ in the obvious way
(cf. the Remarks 18.2 iii) and 19.2 iii)).
\end{remarks}

The {\it spectrum} of a Frobenius manifold $(M,\circ,e,E,g)$ is 
defined as follows. $\nabla$ is the Levi-Civita connection of the 
metric $g$. The operator $\nabla E: \tm\to \tm$, $X\mapsto \nabla_XE$,
acts on the space of flat fields (\cite{Du2}(p 132), \cite{Man}(p 24))
and coincides there with $-{\rm ad\,}E$.
The set of its eigenvalues $\{w_1,...,w_n\}$ is the spectrum 
(\cite{Man}). If $-{\rm ad\,}E$ acts semisimple on the space of flat
fields then there exists locally a basis of flat coordinates 
$z_1,...,z_n$ with 
\begin{eqnarray}
E=\sum_i (w_iz_i+r_i){\partial \over \partial z_i}
\end{eqnarray}
for some $r_i\in \C$.

The following was conjectured by Dubrovin
(\cite{Du1}\cite{Du2}(p 268)).

\begin{theorem}
Let $((M,p),\circ,e,E,g)$ be the germ of a Frobenius manifold with the
following properties: \\
generically semisimple multiplication;\\
$\lie_E(\circ )=1\cdot \circ$ and $\lie_E(g)=D\cdot g$;
\begin{eqnarray}
E=\sum w_iz_i{\partial \over \partial z_i}
\end{eqnarray}
for a basis of flat coordinates $z_i$ with $z_i(p)=0$;\\
positive spectrum
$(w_1,...,w_n)$, that is, $w_i>0$ for all $i$.

Then $(M,p)$ decomposes uniquely into a product of germs at 0 of Frobenius
manifolds for certain irreducible Coxeter groups. The 
Coxeter groups have all the same Coxeter number 
$h={2\over D-1}$.
\end{theorem}

{\bf Proof:}
As in the proof of Theorem 18.4 a), the hypotheses on the Euler field 
show that the Lyashko-Looijenga map $\Lambda:(M_k,p)\to \C^n$ is finite and
that the F-manifold $(M,\circ,e)$ is simple.
One applies Theorem 18.3 and Theorem 19.4.
\hfill $\qed$

\begin{theorem}
Let $((M,p),\circ,e,E,g)$ be the germ of a Frobenius manifold such that
$((M,p),\circ,e,E)$ is isomorphic to the germ at 0 of the F-manifold
of a Coxeter group with the standard Euler field.

Then the irreducible Coxeter subgroups have the same Coxeter number and
$((M,p),\circ,e,E,g)$ is isomorphic to a product of 
germs at 0 of Frobenius manifolds for these Coxeter groups.
\end{theorem}

{\bf Proof:}
First we fix notations. $W$ is a Coxeter group which acts on $V=\C^n$ and 
respects the standard bilinear form. The decomposition of $W$ into $l$
irreducible Coxeter groups $W_1,...,W_l$ corresponds to an orthogonal
decomposition $V=\bigoplus_{k=1}^l V_k$.
The choice of $n$ algebraically independent homogeneous polynomials
$P_1,...,P_n\in \C[x_1,...,x_n]^W$ identifies the quotient $M=V/W$ with
$\C^n$. The quotient map $\psi:V\to M$ decomposes into a product of 
quotient maps $\psi_k:V_k\to V_k/W_k = M_k$. 
The F-manifold $M\cong \prod_{k=1}^lM_k$ is the product of the F-manifolds
$M_k$.

Setting $\varepsilon :=\sum_{i=1}^n x_i{\partial \over \partial x_i}$
on $V$ and $\varepsilon_k:= \varepsilon |_{V_k}$, the standard Euler field 
$E_k$ on $M_k$ is 
$E_k = {1\over h_k} d\psi_k(\varepsilon_k)$. Here $h_k$ is the Coxeter number of 
$W_k$. The Euler field on $M$ is 
\begin{eqnarray}
E=\sum_{k=1}^l E_k = 
\sum_{i=1}^n w_it_i{\partial \over \partial t_i}\ ,
\end{eqnarray}
$\{w_1,...,w_n\}$ is the union of the invariant degrees of $W_k$, 
divided by $h_k$.

Now suppose that $g$ is a flat metric on the germ $(M,0)$ such that 
$((M,0),\circ,e,E,g)$ is a germ of a Frobenius manifold with 
$\lie_E(g)=D\cdot g$.
Consider a system of flat coordinates $z_1,...,z_n$ of $(M,0)$ with
$z_i(0)=0$.
The space $\bigoplus_{i=1}^n \C \cdot {\partial \over \partial z_i}$ 
of flat fields is invariant with respect to $\ad E$ (\cite{Du2} p 132,
\cite{Man} p 24) and the space of affine linear functions
$\C\cdot 1 \oplus
\bigoplus_{i=1}^n \C\cdot z_i\subset \OO_{M,0}$ is invariant with respect to
$E$. 
Because $E$ vanishes at 0 
even the subspace $\bigoplus_{i=1}^n \C\cdot z_i$ is invariant with 
respect to $E$.

The weights $w_1,...,w_n$ of $E$ are positive. Therefore 
the coordinates $z_1,...,z_n$ 
can be chosen to be weighted homogeneous polynomials in $\C[t_1,...,t_n]$
of degree $w_1,...,w_n$.
Thus the spectrum of the Frobenius manifold is $\{w_1,...,w_n\}$.
It is symmetric with respect to ${D\over 2}$, because of 
$\lie_E(g)=D\cdot g$; hence $1+{2\over h_k}=D$ for all $k=1,...,l$.
The Coxeter numbers are all equal, $h:=h_1=...=h_l$, 
the Euler field $E$ is 
\begin{eqnarray}
E={1\over h}d\psi (\varepsilon)\ .
\end{eqnarray}
It rests to show that $g$ is induced as in Theorem 19.1 from a metric
on $V$ which it the orthogonal sum of multiples of the standard metrics
on the subspaces $V_k$.

The operator $\UU=E\circ : \tm \to \tm $ is invertible on $M-\DD$.
The metric $\check{g}$ on $M-\DD$ with 
\begin{eqnarray}
\check{g}(X,Y):=g(\UU^{-1}(X),Y)
\end{eqnarray}
is flat (\cite{Du2} pp 191, 194; \cite{Man}). It lifts to a flat
metric
$\widetilde{g}$ on $V-\psi^{-1}(\DD)$. We claim that $\widetilde{g}$ extends
to a flat metric on the union $\psi^{-1}(\DD)$ of the reflection 
hyperplanes. 

It is sufficient to consider a generic point $p$ in 
one reflection hyperplane.
Then the $e$-orbit of $\psi (p)$ intersects $\DD$ in $n$ points; there
exist canonical coordinates $u_1,...,u_n$ in a neighborhood of $\psi(p)$ with
$e_i\circ e_j=\delta_{ij}e_i$, $g(e_i,e_j)=0$ for $i\neq j$,
\begin{eqnarray*}
&& E=u_1e_1+\sum_{i=2}^n (u_i+r_i)e_i \mbox{ \ \  for some }r_i\in \C^{*}\ ,\\
&& (\DD,\psi(p))\cong (\{u\ |\ u_1=0\},0)\ . 
\end{eqnarray*}
The map germ $\psi:(V,p)\to (M,\psi(p))$ is a twofold covering, 
branched along $(\DD,\psi(p))$, and is given by
$(\widetilde{u_1},...,\widetilde{u_n})\mapsto ({\widetilde{u_1}}^2,
\widetilde{u_2},...,\widetilde{u_n})=(u_1,...,u_n)$
for some suitable local coordinates 
$\widetilde{u_1},...,\widetilde{u_n}$ on $(V,p)$. Then
\begin{eqnarray*}
\widetilde{g}({\partial \over \partial \widetilde{u_1}},
{\partial \over \partial \widetilde{u_1}})=
\check{g}(4u_1e_1,e_1)=
\check{g}(4E\circ e_1,e_1)=
4g(e_1,e_1)\ ,\\
\widetilde{g}({\partial \over \partial \widetilde{u_i}},
{\partial \over \partial \widetilde{u_i}})=
\check{g}(e_i,e_i)=
{1\over u_i+r_i}g(e_i,e_i)\mbox{ \ \ \  for } i\geq 2\ ,\\
\widetilde{g}({\partial \over \partial \widetilde{u_i}},
{\partial \over \partial \widetilde{u_j}})=0
\mbox{ \ \ \  for }i\neq j\ .
\end{eqnarray*}
So $\widetilde{g}$ extends to a nondegenerate (and then flat) metric on $V$. 

The Coxeter group $W$ acts
as group of isometries with respect to $\widetilde{g}$.
It rests to show that the vector space structure
on $V$ which is induced by $\widetilde{g}$ (and $0\in V$) coincides
with the old one.
Then $\widetilde{g}$ is an orthogonal sum of multiples of the standard
metrics on the subspaces $V_k$, because each $W$-invariant quadratic
form is a sum of $W_k$-invariant quadratic forms on the subspaces $V_k$
and they are unique up to scalars.

Let $\widetilde{\varepsilon}$ 
be the vector field on $V$ which corresponds to the
$\C^{*}$-action of the vector space structure induced by $\widetilde{g}$.
Then $\lie_{\widetilde{\varepsilon}} \widetilde{g} = 2\cdot \widetilde{g}$.
Because of 
$\lie_E(\UU)=\UU$, 
$\lie_E(g)=(1+{2\over h})g$, and 
$E={1\over h}d\psi(\varepsilon)$
we have also $\lie_{\varepsilon}{\widetilde{g}} = 2\cdot \widetilde{g}$
for the vector field $\varepsilon$, which corresponds to the
$\C^{*}$-action of the old vector space structure.
The difference $\varepsilon- \widetilde{\varepsilon}$
satisfies 
$\lie_{\varepsilon-\widetilde{\varepsilon}} (\widetilde{g})=0$
and is a generator of a 1-parameter group of isometries. 
As it is also tangent to the union of reflection hyperplanes, it 
vanishes. The vector field $\varepsilon =\widetilde{\varepsilon}$
determines a unique space of linear functions on $V$ and a unique
vector space structure.
\hfill $\qed$

\section{3-dimensional and other F-manifolds}
\setcounter{equation}{0}

The F-manifolds in chapters 16 -- 18 were special in several aspects:
the analytic spectrum was weighted homogeneous and a complete intersection.
Therefore always an Euler field of weight 1 existed, and the tangent 
spaces were Frobenius algebras.
Furthermore, the stratum of points with irreducible germs of 
dimension $\geq 3$ had codimension 2.

Here we want to present examples with different properties.
A partial classification of 3-dimensional germs of massive F-manifolds
will show that already in dimension 3 most germs are not simple
and do not even have an Euler field of weight 1.
Examples of germs $(M,p)$ of simple F-manifolds such that $T_pM$ is not
a Frobenius algebra will complement Theorem 18.3.

\bigskip

First, a construction which is behind the formulas for 
$A_n$, $B_n$, $H_3$, $I_2(m)$
in Theorem 18.5 b)i) provides many other examples.

\begin{proposition}
Fix the following data: $(M,0)=(\C^n,0)$,\\ 
$(Z,0)=(\C,0)\times (M,0)$ 
with coordinates $(x,t)=(x,t_1,...,t_n)$,\\
the projection $\pi_Z:(Z,0)\to (M,0)$,\\
the 1--form $\alpha_Z:= dt_1+xdt_2+...+x^{n-1}dt_n$  on $Z$, \\
the function 
$\widetilde{t_2}(x,t):= t_2+2xt_3+...+(n-1)x^{n-2}t_n$,\\
an isolated plane curve singularity (or a smooth germ)
$f:(\C^2,0)\to (\C,0)$ with
$f(x,0)=x^n\cdot unit\in \C\{x\}$,\\
the subvariety
$C:=\{(x,t)\in Z\ |\ f(x,\widetilde{t_2})=0\}\subset Z$.

\smallskip
a) The map
\begin{eqnarray}
\aaa_C:\TT_{M,0}\to \OO_{C,0},\ 
X\mapsto \alpha_Z(\widetilde{X})|_C
\end{eqnarray}
($\widetilde{X}$ is a lift of $X$ to $Z$) is welldefined and an 
isomorphism of $\OO_{M,0}$-modules.
$(M,0)$ with the induced multiplication on $\TT_{M,0}$ is an irreducible
germ of a massive F-manifold. Its analytic spectrum is isomorphic to 
$(C,0)\cong (\C^{n-1},0)\times (f^{-1}(0),0)$.
For each $t\in M$ the tangent space $T_tM$ is isomorphic to a product of 
algebras $\C\{x\}/(x^k)$ and is a Frobenius algebra.

\smallskip
b) An Euler field of weight 1 exists on $(M,p)$ iff the curve singularity
$f(x,y)$ is weighted homogeneous.

\smallskip
c) Suppose that $\mult f=n$. Then the caustic is 
$\KK = \{t\in M\ |\ t_2=0\}$.
The germ $(M,t)$ is irreducible for all $t\in \KK$, so the caustic 
 is equal to the 
$\mu$-constant stratum of $(M,0)$. The modality is $\modmu (M,0)=n-2$
(the maximal possible).
\end{proposition}

{\bf Proof:}
a) $\alpha_Z$ is exact on $C_{reg}$ because of $d\alpha_Z=dxd\widetilde{t_2}$. 
One can apply Corollary 15.2.

b) Corollary 7.5 b).

c) For $t_2=0$ fixed we have $f(x,\widetilde{t_2})=
x^n\cdot unit\in\C\{x,t_3,...,t_n\}$.
Thus the projection $\pi_C:C\to M$ is a branched covering of degree $n$,
with
$\pi_C^{-1}(\{t\ |\ t_2=0\})=
\{0\}\times \{t\ |\ t_2=0\}$
and unbranched outside of $\{t\ |\ t=2=0\}$.
The analytic spectrum is isomorphic to $C$.
\hfill $\qed$

\begin{remarks}
i) The function $\widetilde{t_2}$ is part of a coordinate system on
$T^{*}M$ for a different Lagrange fibration: the coordinates
\begin{eqnarray}
&& \widetilde{y_1}=y_1,\ \widetilde{y_2}=y_2,\ 
\widetilde{y_i}=y_i-y_2^{i-1}\mbox{ \ for }i\geq 3\ ,\\
&& \widetilde{t_1}=t_1,\ \widetilde{t_2}=t_2+2y_2t_3+...+(n-1)y_2^{n-2}t_n,\ 
\widetilde{t_i}=t_i\mbox{ \ for }i\geq 3 \nonumber
\end{eqnarray}
satisfy
\begin{eqnarray}
\sum_{i=1}^n d \widetilde{y_i}d \widetilde{t_i} = 
\sum_{i=1}^n dy_idt_i = d\alpha\ .
\end{eqnarray}
The analytic spectrum of an F-manifold as in Proposition 20.1 is
\begin{eqnarray}
L=\{(y,t)\in T^{*}M\ |\ \widetilde{y_1}=1, 
f(\widetilde{y_2},\widetilde{t_2})=0,\widetilde{y_i}=0\mbox{ for }i\geq3\}\ .
\end{eqnarray}
It is a product of Lagrange curves.

ii) Another different Lagrange fibration is behind the formulas for
$D_4, F_4, H_4$ in Theorem 18.5 b)ii). There are many possibilities to 
generalize the construction of examples above.
\end{remarks}

\medskip
In dimension 3, there exist up to isomorphism only two irreducible
commutative and associative algebras,
\begin{eqnarray}
Q^{(1)}&:=&\C\{x\}/(x^3) \mbox{ \ \ \ \ and} \\
Q^{(2)}&:=&\C\{x,y\}/(x^2,xy,y^2)\ .
\end{eqnarray}
$Q^{(1)}$ is a Frobenius algebra, $Q^{(2)}$ not.

\begin{theorem}
Let $(M,p)$ be an irreducible germ of a 3-dimensional massive F-manifold
with analytic spectrum $(L,\lambda)\subset T^{*}M$.

\smallskip
a) Suppose $T_pM\cong Q^{(1)}$. Then $(L,\lambda)$ has embedding dimension
3 or 4 and $(L,\lambda)\cong (\C^2,0)\times (C',0)$ for a
plane curve $(C',0)\subset (\C^2,0)$ with $\mult (C',0)\leq 3$.
An Euler field of weight 1 exists iff $(C',0)$ is 
weighted homogeneous.

\smallskip
b) Suppose $T_pM\cong Q^{(1)}$ and 
$(L,\lambda)\cong (\C^2,0)\times (C',0)$ 
with $\mult (C',0)<3$. 
Then $((M,p),\circ,e)$ is one of the germs $A_3, B_3, H_3$.

\smallskip
c) Suppose $T_pM\cong Q^{(1)}$ and 
$(L,\lambda)\cong (\C^2,0)\times (C',0)$ 
with $\mult (C',0)=3$. 
Then the caustic $\KK$ is a smooth surface and coincides with the 
$\mu$-constant stratum; that means, $T_qM\cong Q^{(1)}$ for each 
$q\in \KK$.
The modality is $\modmu (M,p)=1$ (the maximal possible).

\smallskip
d) Suppose $T_pM\cong Q^{(2)}$. Then $(L,\lambda)$ has embedding dimension
5 and $(L,\lambda)\cong (\C,0)\times (L^{(r)},0)$.
Here $(L^{(r)},0)$ is a Lagrange surface with embedding dimension 4.
$\OO_{L^{(r)},0}$ is a  Cohen--Macaulay ring, but not a Gorenstein ring.
\end{theorem}

{\bf Proof:}
a) One chooses coordinates $(t_1,t_2,t_3)$ for $(M,p)$ 
(as in the proof of Lemma 15.3) with 
$e={\partial \over \partial t_1}$ and
\begin{eqnarray}
\C\cdot [{\partial \over \partial t_2}] + \C\cdot [{\partial \over \partial t_3}] 
&=&\mmm \subset T_pM \ , \\
\C\cdot [{\partial \over \partial t_3}] &=&\mmm^2 \subset T_pM\ .
\end{eqnarray}
The dual coordinates on $(T^{*}M,T^{*}_pM)$ are $y_1,...,y_n$.
There exist functions
$a_0, a_1, a_2, b_0, b_1, b_2 \in \C\{t_2,t_3\}$
with
\begin{eqnarray}
L&=&\{(y,t)\ |\ y_1=1,
y_3=b_2y_2^2+b_1y_2+b_0,\nonumber \\
&&\qquad \qquad y_2^3=a_2y_2^2+a_1y_2+a_0\}\ .
\end{eqnarray}
The Hamilton fields of the smooth functions $y_1-1$ and 
$y_3-\sum_{i=0}^2b_iy_2^i$
are ${\partial \over \partial t_1}$ and 
${\partial \over \partial t_3}+...$ in $T^{*}M$.
They are tangent to $L$.
Therefor $(L,\lambda)\cong (\C^2,0)\times (C',0)$ with
$(C',0)\cong (L,\lambda)\cap T_p^{*}M$
(cf. Proposition 7.4).
The statement on the Euler field is contained in Corollary 7.5 b).

b) We have $\mult (C',0)\leq 2$; and the intersection multiplicity of 
$(C',0)$ with a suitable smooth curve is 3. So,
$(C',0)$ is either smooth or a double point or a cusp.
In the first two cases, one can apply the correspondence between
F-manifolds and hypersurface or boundary singularities (Theorem 16.7 and
Theorem 17.7) and the fact that $A_3, B_3$, and $C_3$ are the only
hypersurface or boundary singularities with Milnor number 3.

Suppose $(C',0)$ is a cusp and $((M,p),\circ, e)$ is not $H_3$.
Then it is not simple because of Theorem 18.3. 
The $\mu$-constant stratum $S_\mu=\{q\in M\ |\ T_qM\cong Q^{(1)}\}$
is more than the $e$-orbit of $p$.
It can only be the image in $M$ of the surface of cusp points of $L$,
because at other points $q$ close to $p$
the germ $(M,q)$ is $A_1^3$ or $A_1\times A_2$. 

So at each cusp point $\lambda'$ of $L$ the intersection multiplicity
of $T^{*}_{\pi(\lambda')}M$ and $(L,\lambda')$ is 3.
This property is not preserved by small changes of the Lagrange 
fibration (e.g. as in Remark 20.2 i)).
But Givental proved that a versal Lagrange map is stable
with respect to small changes of the Lagrange fibration
(\cite{Gi2} p 3251, Theorem 3 and its proof). This together with 
chapter 9 yields a contradiction.

c) In this case, at each point $\lambda'$ (close to $\lambda$) of
the surface of singular points of $L$ the intersection multiplicity
of $T^{*}_{\pi(\lambda')}M$ and $(L,\lambda')$ is 3 and 
the map germ $\pi:(L,\lambda')\to (M,\pi(\lambda'))$ is a branched
covering of degree 3.
This implies all the statements.

d) If the embedding dimension of $(L,\lambda)$ is $\leq 4$ then 
$(L,\lambda)\cong (\C^2,0)\times (C',0)$ for some plane curve $(C',0)$ by 
Proposition 7.4. Then
$(L,\lambda)$ is a complete intersection and the tangent spaces $T_pM$
are Frobenius algebras.

So if $T_pM\cong Q^{(2)}$ then the embedding dimension of $(L,\lambda)$ is 5.
The ring
$\OO_{L^{(r)},0}$ is a Cohen--Macaulay ring because the projection
$L^{(r)}\to M^{(r)}$ of  the restricted Lagrange map is finite and flat.
It is not a Gorenstein ring because $T_pM$ is not a Gorenstein ring.
\hfill $\qed$

\bigskip

The next result provides a complete classification and normal forms for those
irreducible germs $(M,p)$ of 3-dimensional massive F-manifolds which
satisfy $T_pM\cong Q^{(1)}$ and whose analytic spectrum consists of
3 components. Part a) gives an explicit construction of all those
F-manifolds.

\begin{theorem}
a) Choose two discrete parameters $p_2, p_3 \in \N$ with
$p_2\geq p_3\geq 2$ and choose $p_3-1$ holomorphic parameters
$(g_0,g_1,...,g_{p_3-2}) \in \C^{*}\times \C^{p_3-2}$
with $g_0\neq 1$ if $p_2=p_3$.

Define $(M,0):= (\C^3,0)$,      \\
$(Z,0):=(\C,0)\times (M,0)$ with coordinates
$(x,t_1,t_2,t_3)=(x,t)$,        \\
$g:=\sum_{i=0}^{p_3-2} g_it_2^i + t_2^{p_3-2}\cdot t_3$,\\
$f_1:=0,\ f_2:=t_2^{p_2},\ f_3:=t_2^{p_3}\cdot g$,\\
$C:=\bigcup_{i=1}^3 C_i := \bigcup_{i=1}^3 \{(x,t)\in Z\ |\ 
x={\partial f_i \over \partial t_2} \} \subset Z$, \\
$b_2(t_2,t_3):= (p_3\cdot g+t_2\cdot 
{\partial g\over \partial t_2})^{-1} \cdot
(p_3 \cdot g +t_2\cdot 
{\partial g\over \partial t_2}
-p_2\cdot t_2^{p_2-p_3})^{-1}$  \\
($b_2$ is a unit in $\C\{t_2,t_3\}$),\\
$b_1(t_2,t_3):= -b_2\cdot p_2\cdot t_2^{p_2-1}$, \\
the 1--form $\alpha_Z:= dt_1+xdt_2+(b_2x^2+b_1x)dt_3$ \ \ on $Z$.

\smallskip
i)
Then
\begin{eqnarray}
{\partial f_i \over \partial t_3} 
&=&  b_2 \bigl( {\partial f_i \over \partial t_2 }\bigr)^2 
+ b_1{\partial f_i \over \partial t_2} \mbox{ \ \ for }i=1,2,3, \\
\alpha_Z|_{C_i} &=& d(t_1+f_i)|_{C_i}
\mbox{ \ \ for }i=1,2,3.
\end{eqnarray}

ii)
The map 
\begin{eqnarray} 
\aaa_C:\TT_{M,0}\to \OO_{C,0},
\ X\mapsto \alpha_Z(\widetilde{X})|_C
\end{eqnarray}
($\widetilde{X}$ is a lift of $X$ to $Z$) is welldefined and an isomorphism
of $\OO_{M,0}$-modules. $(M,0)$ with the induced multiplication on 
$\TT_{M,0}$ is an irreducible germ of a massive F-manifold with
$T_0M\cong Q^{(1)}$.
Its analytic spectrum is isomorphic to 
$(C,0)\cong (\C^2,0)\times (C',0)$, with 
\begin{eqnarray}
(C',0)=(C,0)\cap (\{(x,t)\ |\ t_1=t_3=0\},0)\ .
\end{eqnarray}

iii)
The caustic is $\KK=\{t\in M\ |\ t_2=0\}$ and coincides with the bifurcation
diagram $\BB$ and with the $\mu$-constant stratum.

\smallskip
iv)
The functions $t_1+f_i|_{C_i}$, $i=1,2,3$, combine to a function 
$F:C\to \C$ which is continuous on $C$ and holomorphic on 
$C_{reg}=C\cap \{(x,t)\ |\ t_2\neq 0\}$.
The Euler field $E$ on $M-\KK$ with 
$\aaa_C|_{M-\KK} (E)=F|_{M-\KK}$ is 
\begin{eqnarray}
E=t_1{\partial \over \partial t_1} +
{1\over p_2}t_2{\partial \over \partial t_2} + 
{1\over p_2t_2^{p_3-2}} \bigl((p_2-p_3)g- 
t_2{\partial g \over \partial t_2}\bigr) \cdot
{\partial \over \partial t_3} \ .
\end{eqnarray}
The following conditions are equivalent:
\begin{list}{}{}
\item[$\alpha)$]
$F$ is holomorphic on $C$ and $E$ is holomorphic on $M$,
\item[$\beta)$]
$p_3=2$ or ($p_2=p_3\geq 3$ and $g_i=0$ for $1\leq i < p_3-2$),
\item[$\gamma)$] 
the curve $(C',0)$ is weighted homogeneous.
\end{list}

\smallskip
b) Each irreducible germ $(M,p)$ of a massive F-manifold such that 
$T_pM\cong Q^{(1)}$ and such that $(L,\lambda)$ has 3 components is 
isomorphic to a finite number of normal forms as in a).
The numbers $p_2$ and $p_3$ are determined by $(L,\lambda)$.
The number of isomorphic normal forms is
$\leq 2p_2$ if $p_2>p_3$ and $\leq 6p_2$ if $p_2=p_3$.
\end{theorem}

{\bf Proof:}
a) i) Direct calculation.

ii) $\aaa_C$ is an isomorphism because $b_2$ is a unit in $\C\{t_2,t_3\}$.
One can apply Corollary 15.2 because of (20.11). For the analytic spectrum
see Theorem 20.3 a).

iii) The branched covering $(C,0)\to (M,0)$ is branched along 
$\{(x,t)\ |\ x=t_2=0\}$.
Compare Theorem 20.3 c). 
The generating function $F:C\to \C$ has three different values on 
$\pi^{-1}_C(t)$
for $t\in M$ with $t_2\neq 0$ because of (20.11) and the definition of
$f_i$.
Therefore $\KK=\BB$.

iv) (20.14) can be checked by calculation. 
$\alpha)$$\iff$$ \beta)$ follows. Corollary 7.5 b) shows 
$\alpha)$$\iff$$ \gamma)$ (one can see $\beta)$$\iff$$ \gamma)$ also directly).

\smallskip
b) We start with coordinates $(t_1,t_2,t_3)$ for $(M,p)$ as in the proof
of Theorem 20.3 a). The proofs of Theorem 20.3 a) and Lemma 15.3 b) give
a unique construction of data 
$(Z,0)=(\C,0)\times (\C^3,0)$, 
$(C,0)\subset (Z,0)$ and 
\begin{eqnarray}
\alpha_Z = dt_1+xdt_2+(b_2(t_2,t_3)x^2+b_1(t_2,t_3)x+b_0(t_2,t_3))dt_3
\end{eqnarray}
as in Corollary 15.2 for the germ $(M,p)\cong (\C^3,0)$ of an F-manifold.

$(C,0)=\bigcup_{i=1}^3 (C_i,0)$ is the union of 3 smooth varieties,
which project isomorphically to $(\C^3,0)$, and is isomorphic to the
product of $(\C^2,0)$ and $(C,0)\cap (\{(x,t)\ |\ t_1=t_3=0\},0)$.

The components $(C_i,0)$ can be numbered such that the intersection numbers
of the curves $(C_i,0)\cap (\{(x,t)\ |\ t_1=t_3=0\},0)$ are 
$p_2-1$ for $i=1,2$ and $p_3-1$ for $i=1,3$ and for $i=2,3$.
The numbers $p_2$ and $p_3$ are defined hereby and satisfy 
$p_2\geq p_3\geq 2$.

The 1--form $\alpha_Z$ is exact on $C_{reg}$ and can be integrated to a 
continuous function $F:(C,0)\to (\C,0)$ with 
$F|_{C_i}=t_1+f_i$ for a unique function $f_i\in \C\{t_2,t_3\}$.
Then $C_i=\{(x,t)\ |\ 
x={\partial f_i \over \partial t_2}\}$ and 
\begin{eqnarray}
{\partial f_i\over \partial t_3}= 
b_2\cdot \bigl( {\partial f_i \over \partial t_2 }\bigr)^2 +
b_1\cdot {\partial f_i \over \partial t_2} + b_0\ .
\end{eqnarray}
We will refine $(t_1,t_2,t_3)$ in several steps and change 
$Z$, $C$, $\alpha_Z$, $f_1, f_2, f_3$ accordingly, without explicit mentioning.

{\bf 1st step:}
$(t_1,t_2,t_3)$ can be chosen such that $(C,0)\to (\C^3,0)$ is branched
precisely over $\{t\in \C^3\ |\ t_2=0\}\ .$

{\bf 2nd step:}
$t_1$ can be changed such that $f_1=0$.\\
Then $C_1=\{(x,t)\ |\ x=0\}$ and
\begin{eqnarray}
C_1\cap C_i = \{(x,t)\ |\ x=0={\partial f_i\over \partial t_2}\}
=\{(x,t)\ |\ x=0,t_2=0\}\ .
\end{eqnarray}
Because of $f_i|_{C_1\cap C_i}=f_1|_{C_1\cap C_i}=0$, the functions
$f_2$ and $f_3$ can be written uniquely as 
$f_2=t_2^{\widetilde{p_2}} \cdot \widetilde{g}$, 
$f_3=t_2^{\widetilde{p_3}} \cdot g$
with $\widetilde{p_2}, \widetilde{p_3} \geq 1 $,
$\widetilde{g}, g \in \C\{t_2,t_3\}-t_2\cdot \C\{t_2,t_3\}$.
Now (20.17) shows
${\partial f_i \over \partial t_2} = t_i^{\widetilde{p_i}-1}\cdot unit$
and $\widetilde{p_i}\geq 2$. 
Therefore $\widetilde{p_i}=p_i$ and 
$\widetilde{g}$ and $g$ are units,
with $g(0)\neq \widetilde{g}(0)$ if $p_2=p_3$.

{\bf 3rd step:}
$t_2$ can be chosen such that $f_2=t_2^{p_2}$. \\
(20.16) yields $b_0=0$ for $i=1$ and 
$b_1=-b_2\cdot p_2\cdot t_2^{p_2-1}$ for $i=2$ and
\begin{eqnarray}
t_2^{p_3}{\partial g \over \partial t_3} =
b_2\cdot t_2^{p_3-1}(p_3g+ 
t_2{\partial g \over \partial t_2}) t_2^{p_3-1}
(p_3g+t_2{\partial g \over \partial t_2}-p_2\cdot t_2^{p_2-p_3})&&
\end{eqnarray}
for $i=3$. The first, third and fifth factor on the right are units, therefore
we have ${\partial g \over \partial t_3}=t_2^{p_3-2}\cdot unit$.

{\bf 4th step:}
$t_3$ can be changed such that 
$g=\sum_{i=0}^{p_3-2} g_i\cdot t_2^i + t_2^{p_3-2}\cdot t_3$.

We have brought the germ $(M,p)$
to a normal form as in a).
The numbering of $C_1, C_2, C_3$ was unique up to permutation of 
$C_1$ and $C_2$ if $p_2>p_3$ and arbitrary if $p_2=p_3$. The choice of 
$t_2$ was unique up to a unit root of order $p_2$. Everything else
was unique.
\hfill $\qed$

\begin{remark}
Results of Givental motivate some expectations on the moduli of germs of 
F-manifolds, which are satisfied in the case of Theorem 20.4.

An irreducible germ $(M,p)$ of a massive F-manifold is determined by
its restricted Lagrange map 
$(L^{(r)},\lambda^{(r)})\hookrightarrow (T^{*}M^{(r)},T^{*}_pM^{(r)}) 
\to (M^{(r)},p^{(r)})$ (chapter 9).
Suppose that $(M,p)$ is 3-dimensional with $T_pM\cong Q^{(1)}$.
Then $(L^{(r)},\lambda^{(r)})$ decomposes into a product of two
Lagrange curves, a smooth one and a plane curve $(C',0)$
(Theorem 20.3 a), Proposition 7.4).

If we fix only the topological type of the curve $(C',0)$, we can 
divide the moduli for the possible germs $(M,p)$ into three pieces:
\begin{list}{}{}
\item[i)]
moduli for the complex structure of the germ $(C',0)$,
\item[ii)]
moduli for the Lagrange structure of $(C',0)$,
\item[iii)]
moduli for the Lagrange fibration in the restricted Lagrange map.
\end{list}
Within the $\mu$-constant stratum 
$S_\mu=\{q\in M\ |\ T_qM\cong Q^{(1)}\}$
of a representative $M$, the modules of type i) and ii) are not visible
because the Lagrange structure of the curve $(C',0)$ is constant along
$S_\mu$.

But the modules for the Lagrange fibration are precisely reflected by 
$S_\mu$ because of a result of Givental \cite{Gi2} (proof of Theorem 3):
as a miniversal Lagrange map, the restricted Lagrange map is stable 
with respect to small changes of the Lagrange fibration which 
preserve the symplectic structure; that means, the germ of the 
Lagrange map after such a small change is the restricted Lagrange map
of $(M,q)$ for a point $q\in S_\mu$ close to $p$.

In view of Theorem 14.2 and Theorem 20.3 b)+c) there is 1 module of type
iii) if $\mult (C',0)=3$ and no module of type iii) if 
$\mult (C',0)=1$ or 2. 

Fixing the complex structure of the plane curve $(C',0)$, the choice of
a Lagrange structure is equivalent to the choice of a volume form.
Equivalence classes of it are locally parametrized by
$H_{Giv}^1((C',0))$ (\cite{Gi2}(Theorem 1), \cite {Va}),
so the number of moduli of type ii) is $\mu-\tau$ (Theorem 7.2 b)).
It is equal to the number of moduli of right equivalence classes
of function germs $f:(\C^2,0)\to (\C,0)$ with 
$(f^{-1}(0),0)\cong (C',0)$. 

The $\mu$-constant stratum of a plane curve singularity in the 
semiuniversal unfolding is smooth by a result of Wahl (\cite{Wah}, cf. also
\cite{Mat}) and its 
dimension depends only on the topological type of the curve. 
So one may expect that the number of moduli of types i) and ii) together
depends only on the topological type of $(C',0)$ and is equal
to this dimension.

(But a canonical relation between the choice of a Lagrange structure and
the choice of a function germ for a plane curve $(C',0)$ is not known.)

In the case of Theorem 20.4 these expectations are met:
the topological type of the plane curve is given by the intersection
numbers $p_2-1$ and $p_3-1$; the last one $g_{p_3-2}$ of the complex modules
is of type iii), 
it is the module for the $\mu$-constant stratum and for the Lagrange fibration;
the other $p_3-2$ modules $(g_0,...,g_{p_3-3})$ are of type i) and ii).
One can check with \cite{Mat}(4.2.1) that $p_3-2$ is the dimension of the 
$\mu$-constant stratum for such a plane curve singularity.
\end{remark}

\smallskip
Finally, at least a few examples of germs $(M,p)$ of F-manifolds with
$T_0M\cong Q^{(2)}$ shall be presented.

\begin{proposition}
Consider $M=\C^3$ with coordinates $(t_1,t_2,t_3)$ and $T^{*}M$ with 
fiber coordinates $y_1,y_2,y_3$.
Choose $p_2, p_3\in\N_{\geq 2}$. Then the variety
\begin{eqnarray}
L=\{(y,t)\in T^{*}M &|& y_1=1, \ y_2(y_2-p_2t_2^{p_2-1})=y_2y_3\nonumber \\
&&=y_3(y_3-p_3t_3^{p_3-1})=0\}
\end{eqnarray}
has three smooth components
and is the analytic spectrum of the structure of a simple F-manifold on $M$
with $T_0M\cong Q^{(2)}$. The field
\begin{eqnarray}
E=t_1{\partial \over \partial t_1} 
+ {1\over p_2}t_2{\partial \over \partial t_2}
+ {1\over p_3}t_3{\partial \over \partial t_3}
\end{eqnarray}
is an Euler field of weight 1.
\end{proposition}

{\bf Proof:}
One checks easily that $\alpha=y_1dt_1+y_2dt_2+y_3dt_3$ is exact on the
three components of $L$, that the map
$\aaa_C:\tm \to \pi_{*}\OO_L$ is an isomorphism of $\OO_M$-modules,
and that $E$ in (20.20) is an Euler field.

The weights of $E$ are positive. This shows via the Lyashko-Looijenga map
that $(M,\circ,e)$ is a simple F-manifold (cf. the proof of 
Theorem 18.4 a)).
\hfill $\qed$

\begin{remark}
In \cite{Gi2}(Theorem 15) the restricted Lagrange maps of two other 
series of simple F-manifolds with $M\cong \C^n$ and 
$T_0M$ not a Frobenius manifold are given,
the series $\Xi_n$ ($n\geq 3$) and $\Omega_n$ ($n\geq 4$) (also
$\Xi_1=A_1$, $\Xi_2=H_2$, $\Omega_2=A_2$, $\Omega_3=H_3$).

They have Euler fields of weight 1 with positive weights. 
The analytic spectra of $\Xi_n$ and $\Omega_n$ are isomorphic to
$\C\times \Sigma_{n-1}(2n-1)$ and $\C^2\times \Sigma_{n-2}(2n-3)$, 
respectively. Here $\Sigma_k(2k+1)$ is the {\it open swallowtail},
the subset of polynomials in the set of polynomials
$\{z^{2k+1}+a_2z^{2k-1}+...+a_{2k+1}\ |\ a_2,...,a_{2k+1}\in \C\}$
which have a root of multiplicity $\geq k+1$ (\cite{Gi2} p 3256).
It has embedding dimension $2k$.

The germs $((M,0),\circ,e)$ are irreducible for $\Xi_n$ and 
$\Omega_n$, the socle $\Ann_{T_0M}(\mmm)$ of $T_0M$ is the maximal ideal
$\mmm \subset T_0M$ itself in the case of $\Xi_n$ and has 
dimension $n-2$ in the case of $\Omega_n$. 

Givental \cite{Gi2}(Theorem 15) proved that the germs $(M,0)$
for $\Xi_n$ and $\Omega_n$ are the only irreducible germs of simple
F-manifolds whose analytic spectra are products of smooth germs and 
open swallowtails. Generating functions in the sense of 
Definition 15.4 are due to O.P. Shcherbak and are given in 
\cite{Gi2}(Proposition 12).
\end{remark}


\end{document}